\documentclass[11pt]{article}

\usepackage[letterpaper,margin=1in]{geometry}
\usepackage[T1]{fontenc}
\usepackage{lmodern}
\usepackage{graphicx}
\usepackage{multirow}
\usepackage{amsmath,amssymb,amsfonts}
\usepackage{amsthm}
\usepackage{mathrsfs}
\usepackage{booktabs}
\usepackage{float}
\usepackage{array}
\usepackage{mathtools}
\usepackage{xcolor}
\usepackage{adjustbox}
\usepackage[font=small,labelfont=bf,labelsep=space,skip=6pt]{caption}
\usepackage[numbers,sort&compress]{natbib}
\usepackage[title]{appendix}
\usepackage[hyphens]{url}
\usepackage{hyperref}
\hypersetup{colorlinks=true,linkcolor=blue!60!black,citecolor=blue!60!black,urlcolor=blue!60!black}

\theoremstyle{plain}
\newtheorem{theorem}{Theorem}
\newtheorem{proposition}[theorem]{Proposition}
\newtheorem{lemma}[theorem]{Lemma}
\newtheorem{corollary}[theorem]{Corollary}
\theoremstyle{definition}

\theoremstyle{remark}
\newtheorem{remark}[theorem]{Remark}

\newenvironment{plaintable}[1][tbp]{\begin{table}[#1]\centering\small
  \renewcommand{\resizebox}[3]{\adjustbox{max width=##1}{##3}}}{\end{table}}
\newcommand{\backmatter}{}

\allowdisplaybreaks[4]
\begin{document}

\begin{center}
{\LARGE\bfseries Exact Spinal Optima and Numerical Bounds for the\\[2pt] Single-Root Maximum Overhang Problem\par}
\vspace{1.2em}
{\large Aolong Zha\textsuperscript{1,*}, Qianwei Zhuang\textsuperscript{2}, Qiong Chang\textsuperscript{3}\par}
\vspace{0.8em}
{\small
\textsuperscript{1}Faculty of Systems Engineering, Wakayama University, Japan\\
\textsuperscript{2}Research Center for Advanced Science and Technology, The University of Tokyo, Japan\\
\textsuperscript{3}School of Computing, Institute of Science Tokyo, Japan\\[4pt]
\textsuperscript{*}Corresponding author: \href{mailto:cyouryuuryuu@gmail.com}{cyouryuuryuu@gmail.com}\\
Contributing authors: qweizhuang@gmail.com; changqiong30@gmail.com\par}
\end{center}

\begin{abstract}
The maximum overhang problem seeks the largest horizontal extension of a balanced stack of $N$ identical frictionless unit blocks on a horizontal table. For general stacks, optimal overhangs lie between order-$N^{1/3}$ constructions (Paterson and Zwick, 2009) and a matching $6N^{1/3}$ upper bound (Paterson et al., 2009), while candidate optima for $N \le 30$ were computed numerically without optimality proofs. Formulating single-root stacking as a nonconvex mixed-integer bilinear program, we develop a global optimization framework combining structural reductions, a loaded-spine equilibrium bound, and layer-profile decomposition. We prove a base-singleton reduction, a base-shift lemma, cut-set difference linearizations, and a second-moment energy criterion that prunes infeasible profiles analytically. We also establish a loaded-spine chain bound for arbitrary stability margins, matching Paterson and Zwick's loaded-spine optimum at zero margin, and convert it into sparse, globally valid big-$M$ spinal cuts. With exact rational and algebraic stacks checked without any solver, we prove, for $N \in \{4, 6, 8, \dots, 19\}$, Paterson and Zwick's conjecture that spinal stacks attain the loaded-spine bound, and we determine the spinal optima at the exception sizes $N \in \{5, 7\}$ in closed form. A spatial branch-and-bound solver reports every subproblem of the framework closed for $N \in \{5, \dots, 16\}$; read as exact, these numerical certificates bound the optimal single-root overhangs to relative gaps below $7.8 \times 10^{-6}$ and show that every single-root stack whose overhang exceeds $99\%$ of the loaded-spine bound, rounded upward at the seventh decimal, is spinal. If these solver closures are correct, the exact spinal optima are the optimal single-root overhangs for $N \in \{5, \dots, 16\}$.

\end{abstract}

\noindent{\small\textbf{Keywords:} Global optimization, Mixed-integer bilinear programming, Spatial branch-and-bound, Maximum overhang problem, Exact rational verification, Numerical optimality certificates\par
\par}

\section{Introduction}
\label{sec:Intro}

The maximum overhang problem is a classical question in discrete geometry and rigid-body statics~\cite{coffin1923,johnson1955}. Given $N$ identical, homogeneous, frictionless rectangular blocks of unit length and unit weight resting on a horizontal table, the objective is to determine the maximum horizontal distance from the table edge to the rightmost edge of any block in a statically balanced stack. When each block supports at most one block above it, the optimum is attained by the classical harmonic stack, in which the $i$-th block from the top extends by $1 / (2i)$ beyond the block beneath it and yields total overhang $\frac{1}{2} H_N = \frac{1}{2} \sum_{i=1}^{N} 1/i \sim \frac{1}{2} \ln N$~\cite{johnson1955}.

Removing the one-on-one restriction lets several blocks share a layer as counterweights and force splitters, so that the supports form a directed acyclic graph (DAG)~\cite{ainley1979,hall2005}. For general stacks, the maximum overhang is of order $N^{1/3}$: parabolic constructions attain this order~\cite[Theorem~4.4]{paterson2009overhang}, and Paterson, Peres, Thorup, Winkler, and Zwick proved the matching upper bound $6N^{1/3}$~\cite[Theorem~2.3]{paterson2009maximum}. Exact optima at finite $N$, with proofs of global optimality, are a different matter. For $N \le 30$, Paterson and Zwick~\cite{paterson2009overhang} enumerated combinatorially distinct arrangements and optimized block positions numerically, displaying optimal stacks for $2 \le N \le 19$ in their Figures~3 and~4 and for $N \in \{20, 30\}$ in Figure~5; they state that they are confident of the optimality of these stacks but have no formal proofs, since numerical techniques were used~\cite[Section~1]{paterson2009overhang}. Their study also indicated that optimal stacks are \emph{spinal}, with a single support block in every layer up to the principal block, for $N \le 19$, and that non-spinal stacks take over at $N = 20$~\cite{paterson2006overhang,paterson2009overhang}.

Optimizing over both support topologies and block coordinates gives a nonconvex mixed-integer bilinear program (MIBLP)~\cite{adams1993,gupte2013}. Checking equilibrium for fixed positions is a linear program (LP)~\cite{zwick2002,paterson2009overhang}, and initial collapse is a convex quadratic program~\cite{rote2011}, but optimizing positions and forces together couples each force $F_{i,j}$ with its unknown contact point $P_{i,j}$ in the moments $\tau_{i,j} = F_{i,j} P_{i,j}$. Spatial branch-and-bound (B\&B) then faces three obstacles. First, McCormick relaxations~\cite{mccormick1976,alkhayyal1983} over global coordinate boxes have an $O(N^2)$ error per contact pair~\cite{luedtke2012,boland2017}, and big-$M$ contact-interval constraints are weak at fractional indicators~\cite{ceria1999,vielma2015}. Second, integer layer assignments create combinatorial symmetry~\cite{atamturk2003,dey2024}. Third, subproblems that contain near-optimal spinal stacks stall when relative gaps of order $10^{-5}$ are certified, and double-precision incumbents at the usual feasibility tolerance of $10^{-6}$ can exceed the loaded-spine optimum $S^*(N)$.

We develop a global optimization framework for the single-root maximum overhang problem that addresses these obstacles. All optimal stacks displayed by Paterson and Zwick~\cite[Figures~3--5]{paterson2009overhang} rest on a single block on the table. Our numerical certificates apply within this class of single-root stacks and hold up to the numerical tolerances of the solver; we do not claim to exclude stacks with several blocks on the table. Within this class, the exact values of our rational stacks agree, after rounding, with the numerical configurations reported by Paterson and Zwick for all $N \in \{4, \dots, 16\}$.

\subsection{Related Work}
\label{subsec:RelatedWork}

\paragraph{Block-Stacking Equilibrium and Overhang Theory.}
Ainley~\cite{ainley1979} showed that a counterweight in the second layer raises the overhang for $N = 4$ to $(15 - 4\sqrt{2}) / 8 \approx 1.16789$ (see also~\cite[Figure~1]{paterson2009maximum}), and Hall~\cite{hall2005} raised the problem of maximizing the overhang without the one-on-one restriction and proposed counterbalanced stacks that he claimed, without proof, to be optimal. Paterson and Zwick~\cite{paterson2006overhang,paterson2009overhang} decomposed every balanced stack into a support set and a balancing set, which they abstracted as point weights on the support set. For spinal stacks of weight $w$ with $k$ spine blocks, they proved that the loads of an optimal loaded spine with overhang $S_k^*(w)$ satisfy the recurrence $t_i^2 = (t_{i-1} + 1/2) t_{i+1}$ below some index, above which all point weights vanish, and that $\ln w - 1.313 < S^*(w) < \ln w + 1$~\cite[Lemmas~3.3--3.5 and Theorems~3.6 and~3.8]{paterson2009overhang}. With a shield-and-tower heuristic, they realized optimal loaded spines with ordinary blocks for every $n \le 1000$ except $n \in \{3, 5, 7\}$, which led to their Conjecture~3.9 that $S(n) = S^*(n)$ for all $n \neq 3, 5, 7$, where $S(n)$ is the maximum overhang of an $n$-block spinal stack. They also observed that the stacks of Hall are spinal, so that their overhang grows only as $\ln n + O(1)$, and that they are optimal only for $n \le 19$. Paterson, Peres, Thorup, Winkler, and Zwick~\cite{paterson2009maximum} proved the upper bound $6N^{1/3}$, refinable to $4.5N^{1/3}$ for large $N$, by reducing force transmission across horizontal slices to mass-movement problems. Treeby~\cite{treeby2018} studied blocks of different lengths, and Gmeiner and Schulz~\cite{gmeiner2026} proved that counterbalanced stacking of such blocks is NP-hard and established right-alignment properties of optimal stacks. Equilibrium for fixed coordinates is an LP with $O(N)$ variables and constraints when interface resultants are concentrated at the ends of the contact intervals~\cite{livesley1978,zwick2002,paterson2009overhang}; when positions and topologies are optimized jointly, the resulting nonconvex programs have so far been handled by enumerating topologies and applying local solvers~\cite[Section~5]{paterson2009overhang}, without global dual bounds.

\paragraph{Nonconvex MIBLP, Polyhedral Strengthening, and Global Certification.}
Spatial B\&B solvers relax bilinear equations $\tau = F P$ over boxes by McCormick envelopes~\cite{mccormick1976,ryoo1996,tawarmalani2005}, the convex hull of a single bilinear graph~\cite{alkhayyal1983}, whose error grows with the product of the interval widths~\cite{luedtke2012,boland2017}. Strengthening relies on the reformulation-linearization technique (RLT)~\cite{sherali1990}, convex-hull results for bilinear products over structured polytopes and bipartite graphs~\cite{davarnia2017,dey2019,santana2020,anstreicher2021,gu2023,gu2024}, disjunctive and perspective formulations against big-$M$ weakness~\cite{ceria1999,frangioni2006,gunluk2010,vielma2015,wei2024}, bounds tightening~\cite{belotti2009,gleixner2017,xia2020}, and sparse constraint matrices~\cite{atamturk2003,dey2024}. We combine problem-specific structural reductions and an exact layer-profile decomposition with sparse big-$M$ cuts derived from a loaded-spine equilibrium bound, which separate the analytically understood spinal regime from the rest of each subproblem.

\subsection{Main Contributions and Paper Organization}
\label{subsec:Contributions}

Our contributions are fivefold:
\begin{enumerate}
\item \textbf{Structural reductions and energy pruning:} We prove a base-singleton reduction, base right-corner dominance, and a base-shift lemma that extends these conditions without loss of generality to every feasible single-root stack with two second-layer blocks, together with cut-set linearizations and a second-moment energy criterion that discards infeasible layer profiles before branching.
\item \textbf{Spinal cuts:} We prove a loaded-spine bound for every unique-support chain at every stability margin, which coincides at zero margin with the loaded-spine optimum of Paterson and Zwick~\cite[Lemmas~3.4--3.5]{paterson2009overhang}, and turn it into a sparse family of valid big-$M$ cuts inside an exact layer-profile decomposition.
\item \textbf{Exact lower bounds and exact spinal optima:} We give rational stacks for all $N \in \{4, \dots, 16\}$ that are verified in exact rational arithmetic. With the exact algebraic loads of the loaded-spine optimum, we obtain stacks of ordinary blocks that attain $S^*(N)$, which proves the equality $S(N) = S^*(N)$ of Conjecture~3.9 of~\cite{paterson2009overhang} for $N \in \{4, 6, 8, \dots, 19\}$, and we determine $S(5) = (12 - \sqrt{30})/5$ and $S(7)$ in closed form, excluding every other spine height and every face with a vanishing added load at these exception sizes. Paterson and Zwick computed the optimal loads and numerical configurations; the new results are exact realizations and upper bounds that can be checked independently, collected in Theorem~\ref{thm:exact_spinal_main}. They use no optimization solver.
\item \textbf{Numerical upper bounds for $N \in \{5, \dots, 16\}$:} Within the single-root class, a spatial branch-and-bound solver reports every subproblem of one decomposition closed for every $N \in \{5, \dots, 16\}$. Read as exact, these numerical certificates give relative gaps between $4.67 \times 10^{-6}$ and $7.68 \times 10^{-6}$ and show that no single-root stack attains $S^*(N)$ at the exception sizes $N \in \{5, 7\}$.
\item \textbf{Numerical exclusion of non-spinal stacks:} The solver also reports every non-spinal branch closed, so that, read as exact, every single-root stack whose overhang exceeds $99\%$ of $S^*(N)$, rounded upward at the seventh decimal, is spinal for $N \in \{5, \dots, 16\}$; to our knowledge, this is the first computational certificate, within the single-root class, of the numerical finding of Paterson and Zwick that optimal stacks are spinal over this range. If these solver closures are correct, the exact spinal optima are the optimal single-root overhangs for every $N \in \{5, \dots, 16\}$. Theorem~\ref{thm:exact_spinal_main} collects the results verified in exact arithmetic, and Theorem~\ref{thm:main_certified_optimality} those that rest on the solver.
\end{enumerate}

Section~\ref{sec:Formulation} states the model and the baseline formulation, Section~\ref{sec:Structural} the structural reductions, and Section~\ref{sec:SpinalCert} the layer-profile decomposition, the spinal cuts, the non-spinal exclusion, and a classification of all constraints by validity. Section~\ref{sec:Algorithm} describes the certification pipeline, the exact rational stacks, and the theorem on the exact spinal optima, Section~\ref{sec:Computational} reports the computations, and Section~\ref{sec:Conclusion} concludes. Appendix~\ref{app:RationalSeeds} lists the rational stacks, Appendix~\ref{app:proofs} contains the proofs not given in the main text and the solver settings, and Appendix~\ref{app:ExactSpinal} proves the exact spinal optima.

\section{Problem Formulation}
\label{sec:Formulation}

This section states the physical model with a contact stability margin $\eta \in [0, 1)$, the class of single-root stacks, the baseline mixed-integer bilinear formulation $\mathcal{M}_0$, two elementary reductions of it, and the error of its McCormick relaxation.

\subsection{Physical Model, Stability Margin, and Single-Root Stacks}
\label{subsec:PhysicalModel}

We consider a two-dimensional stack of $N$ identical, homogeneous, frictionless rectangular blocks of length $2W = 1$, height $h > 0$, and unit weight on a horizontal table~\cite{paterson2009overhang,paterson2009maximum}. Equilibrium does not depend on $h$, so only horizontal coordinates and layer indices enter the problem. The table, block~$B_0$, occupies $(-\infty, 0] \times (-\infty, 0]$.

Let $\mathcal{N} := \{1, \dots, N\}$ and $\mathcal{N}_0 := \mathcal{N} \cup \{0\}$. Block~$i \in \mathcal{N}$ has center $X_i \in \mathbb{R}$ and layer $L_i \in \{1, \dots, N\}$ and occupies $[X_i - W, X_i + W] \times [(L_i - 1)h, L_i h]$; the table has $X_0 := 0$ and $L_0 := 0$. Blocks in the same layer have disjoint interiors:
\begin{equation}
\label{eq:same_layer_sep}
|X_i - X_k| \ge 2W = 1 \qquad (\forall i, k \in \mathcal{N},\; i \neq k,\; L_i = L_k).
\end{equation}
The overhang of the stack is the horizontal coordinate of the rightmost physical edge among all blocks:
\begin{equation}
\label{eq:overhang_def}
z := \max_{i \in \mathcal{N}} (X_i + W).
\end{equation}

Optimal stacks in the classical model are precarious, with forces acting at block corners and at the table edge; stronger notions of stability require every contact resultant to lie inside its contact interval~\cite[Section~6]{paterson2009overhang}. We therefore introduce a relative stability margin $\eta \in [0, 1)$, the edge clearance $\epsilon := \eta W = \eta / 2$, and the effective half-width
\begin{equation}
\label{eq:effective_halfwidth}
w := W - \epsilon = \frac{1 - \eta}{2} \in (0, 1/2],
\end{equation}
so that $2w = 1 - \eta$. The margin shrinks only the support core of each block to $[X_i - w, X_i + w]$ and the support region of the table to $(-\infty, -\epsilon]$; the non-overlap width $2W$ in~\eqref{eq:same_layer_sep} and the overhang edge $X_i + W$ in~\eqref{eq:overhang_def} are unchanged. For $\eta = 0$, we obtain the classical frictionless model of~\cite{paterson2009overhang} with one relaxation: our contact intervals are closed, so two blocks whose cores touch at a single point may exchange a force there, whereas in~\cite{paterson2009overhang} only contacts of positive length carry force. Upper bounds proved in our model therefore hold for theirs, and the rational stacks of Appendix~\ref{app:RationalSeeds} use only contacts of positive length; Appendix~\ref{app:pz_relation} states what our results imply for their model.

Block~$i \in \mathcal{N}$ can be supported by block~$j \in \mathcal{N}$ only if $L_i = L_j + 1$ and their effective cores intersect, and by the table only if $L_i = 1$ and $[X_i - w, X_i + w] \cap (-\infty, -\epsilon] \neq \emptyset$. Binary variables $T_{i,j}$ for $i \in \mathcal{N}$ and $j \in \mathcal{N}_0 \setminus \{i\}$ encode contacts; a contact with $T_{i,j} = 1$ admits a force but may carry none, so a stack generally has several contact representations $T$ (Lemma~\ref{lem:maximal_contacts}). Since blocks and table are frictionless and cohesionless, all forces are vertical and compressive~\cite{livesley1978,paterson2009overhang}: for $T_{i,j} = 1$, block~$j$ exerts an upward force $F_{i,j} \ge 0$ on block~$i$, and block~$i$ an equal downward force on block~$j$, at a contact point in the $\eta$-shrunk contact interval
\begin{equation}
\label{eq:contact_interval_eta}
P_{i,j} \in I_{i,j}(\eta) := \begin{cases}
[\max\{X_i - w, X_j - w\},\; \min\{X_i + w, X_j + w\}], & \text{if } j \in \mathcal{N}, \\
[X_i - w,\; \min\{X_i + w, -\epsilon\}], & \text{if } j = 0.
\end{cases}
\end{equation}
With the moment $\tau_{i,j} := F_{i,j} P_{i,j}$ about the origin, a stack is \emph{$\eta$-stable}, or \emph{balanced} when $\eta = 0$, if there exist $F_{i,j} \ge 0$ and $P_{i,j} \in I_{i,j}(\eta)$, with $F_{i,j} = P_{i,j} = \tau_{i,j} = 0$ whenever $T_{i,j} = 0$, that satisfy force and moment balance for every block:
\begin{align}
  \sum_{j \in \mathcal{N}_0 \setminus \{i\}} F_{i,j} &= 1 + \sum_{k \in \mathcal{N} \setminus \{i\}} F_{k,i}, && \forall i \in \mathcal{N}, \label{eq:force_bal} \\
  \sum_{j \in \mathcal{N}_0 \setminus \{i\}} \tau_{i,j} &= X_i + \sum_{k \in \mathcal{N} \setminus \{i\}} \tau_{k,i}, && \forall i \in \mathcal{N}, \label{eq:moment_bal} \\
  \tau_{i,j} &= F_{i,j} P_{i,j}, && \forall i \in \mathcal{N},\, j \in \mathcal{N}_0 \setminus \{i\}. \label{eq:bilinear_torque}
\end{align}

In general, several blocks may rest on the table, either in components to the left of the origin or as a foundation of several blocks. The harmonic stacks~\cite{johnson1955}, the counterbalanced stacks of Hall~\cite{hall2005}, and every optimal stack displayed by Paterson and Zwick~\cite[Figures~3--5]{paterson2009overhang}, for $2 \le N \le 20$ and $N = 30$, rest on a single block on the table, but no proof rules out optima with several. We therefore work within the class of \emph{single-root} stacks and leave the general case open (Section~\ref{sec:Conclusion}). We denote by $z^*(N, \eta)$ the maximum overhang of a feasible single-root $N$-block $\eta$-stable stack and write $z^*(N) := z^*(N, 0)$.

\begin{proposition}[Single-Root Canonical Indexing]
\label{prop:single_root}
Let $N \ge 2$ and $\eta \in [0, 1)$. Any feasible single-root $N$-block $\eta$-stable stack (that is, a stack in which exactly one block lies in layer~$1$ resting on the table $B_0$) can be indexed uniquely in canonical layer-and-coordinate order such that:
\begin{enumerate}
\item[(i)] block~$1$ is the unique base block on the table, with $L_1 = 1$, $T_{1,0} = 1$, and $L_i \ge 2$, $T_{i,0} = 0$ for every $i \in \{2, \dots, N\}$;
\item[(ii)] layer indices are nondecreasing ($1 = L_1 \le L_2 \le \dots \le L_N$) and blocks within each layer are ordered strictly from right to left ($X_i - X_{i+1} \ge 1$ whenever $L_i = L_{i+1}$).
\end{enumerate}
Moreover, if a general stack contains a disconnected component resting independently on the table that transmits zero force to the support set of the principal overhanging block, removing that component preserves $\eta$-stability and overhang with fewer blocks.
\end{proposition}

\begin{proof}
Sorting the blocks lexicographically by $(L_i, -X_i)$ gives the indexing: the base block is the only block in layer~$1$, and~\eqref{eq:same_layer_sep} makes the order within each layer strict and unique. A component resting independently on the table exerts no force or moment on the support set of the principal block, so deleting it preserves~\eqref{eq:force_bal}--\eqref{eq:bilinear_torque} for the remaining blocks and the overhang.
\end{proof}

\subsection{Baseline MIBLP Formulation}
\label{subsec:BaselineM0}

The formulation enforces the canonical indexing of Proposition~\ref{prop:single_root} with binary layer-jump indicators $b_i$, $i \in \{1, \dots, N - 1\}$, where $b_i = 1$ if and only if $L_{i+1} > L_i$.

Restricting the centers to $[-N, N]$ loses no optimality, and no generality for any overhang threshold. Translating a feasible stack rightward by $-\epsilon - P_{1,0} \ge 0$ preserves feasibility, decreases no overhang, and gives $P_{1,0} = -\epsilon$, hence $X_1 \in [-W, W - 2\epsilon]$; since every block is linked to block~$1$ by a support path of at most $N - 1$ contacts, each with $|X_i - X_j| \le 2w \le 1$, every block satisfies $|X_i| \le N - 1/2$. We therefore set $M_X := N$, $M_F := N$, $M_\tau := M_F M_X = N^2$, and $M := 2 M_X + 1 = 2N + 1$. With the equilibrium equations~\eqref{eq:force_bal}--\eqref{eq:bilinear_torque}, the baseline single-root formulation $\mathcal{M}_0(N, \eta)$ is

\begingroup\small
\begin{align}
  \max \quad & z \label{eq:M0_obj} \\
  \text{s.t.} \quad & z = \max_{i \in \mathcal{N}} R_i, \quad R_i = X_i + W, && \forall i \in \mathcal{N}, \label{eq:M0_right_edge} \\
  & X_0 = 0, \quad L_0 = 0, \quad L_1 = 1, \quad T_{1,0} = 1, && \label{eq:M0_base_root} \\
  & L_i \ge 2, \quad T_{i,0} = 0, \quad T_{i,i} = 0, && \forall i \in \mathcal{N} \setminus \{1\}, \label{eq:M0_upper_root} \\
  & \sum_{j \in \mathcal{N}_0 \setminus \{i\}} T_{i,j} \ge 1, && \forall i \in \mathcal{N}, \label{eq:M0_support_min} \\
  & 0 \le F_{i,j} \le M_F T_{i,j}, && \forall i \in \mathcal{N},\, j \in \mathcal{N}_0 \setminus \{i\}, \label{eq:M0_force_cap} \\
  & -M_X T_{i,j} \le P_{i,j} \le M_X T_{i,j}, && \forall i \in \mathcal{N},\, j \in \mathcal{N}_0 \setminus \{i\}, \label{eq:M0_p_cap} \\
  & -M_\tau T_{i,j} \le \tau_{i,j} \le M_\tau T_{i,j}, && \forall i \in \mathcal{N},\, j \in \mathcal{N}_0 \setminus \{i\}, \label{eq:M0_tau_cap} \\
  & X_i - w - M(1 - T_{i,j}) \le P_{i,j} \le X_i + w + M(1 - T_{i,j}), && \forall i \in \mathcal{N},\, j \in \mathcal{N}_0 \setminus \{i\}, \label{eq:M0_bigm_i} \\
  & X_j - w - M(1 - T_{i,j}) \le P_{i,j} \le X_j + w + M(1 - T_{i,j}), && \forall i \in \mathcal{N},\, j \in \mathcal{N} \setminus \{i\}, \label{eq:M0_bigm_j} \\
  & P_{i,0} \le -\epsilon + M(1 - T_{i,0}), && \forall i \in \mathcal{N}, \label{eq:M0_bigm_table} \\
  & 1 - N(1 - T_{i,j}) \le L_i - L_j \le 1 + N(1 - T_{i,j}), && \forall i \in \mathcal{N},\, j \in \mathcal{N}_0 \setminus \{i\}, \label{eq:M0_layer_link} \\
  & b_i \le L_{i+1} - L_i \le N b_i, && \forall i \in \{1, \dots, N - 1\}, \label{eq:M0_layer_jump} \\
  & X_i - X_{i+1} \ge 1 - M b_i, && \forall i \in \{1, \dots, N - 1\}, \label{eq:M0_non_overlap} \\
  & \sum_{i=1}^N X_i \le -N\epsilon, \label{eq:M0_global_com} \\
  & \mathrlap{X_i \in [-M_X, M_X],\; L_i \in \{1, \dots, N\},\; T_{i,j} \in \{0, 1\},\; b_i \in \{0, 1\}.} \label{eq:M0_domains}
\end{align}
\endgroup
The maximum in~\eqref{eq:M0_right_edge} is modeled by binary selection constraints or a general maximum constraint of the solver, since its linear relaxation alone only gives $z \ge R_i$; it simplifies once the layers are fixed (Section~\ref{sec:SpinalCert}). Inequality~\eqref{eq:M0_global_com} follows from summing~\eqref{eq:moment_bal} and~\eqref{eq:force_bal} over all blocks together with $P_{1,0} \le -\epsilon$ and is included as a valid center-of-mass constraint.

\subsection{Upper-Triangular DAG Reduction and Exact Base Linearization}
\label{subsec:UpperDAG}

Two elementary consequences of the canonical order eliminate more than half of the contact variables and linearize the base reaction exactly.

\begin{lemma}[Upper-Triangular Variable Elimination and Exact Base Linearization]
\label{lem:upper_dag_base}
Let $N \ge 2$ and $\eta \in [0, 1)$. Every feasible solution of the baseline single-root formulation $\mathcal{M}_0(N, \eta)$ satisfies the following conditions unconditionally:
\begin{enumerate}
\item[(i)] \textbf{Upper-triangular support structure:} $T_{i,j} = F_{i,j} = P_{i,j} = \tau_{i,j} = 0$ for every pair $(i, j) \in \mathcal{N} \times \mathcal{N}$ with $j \ge i$, and $T_{i,0} = F_{i,0} = P_{i,0} = \tau_{i,0} = 0$ for every $i \in \{2, \dots, N\}$.
\item[(ii)] \textbf{Layer-index bounds and upper force capacities:} Every block $i \in \{2, \dots, N\}$ satisfies $2 \le L_i \le i$, $\sum_{j=1}^{i-1} T_{i,j} \ge 1$, and
\begin{equation}
\label{eq:upper_force_cap}
0 \le F_{i,j} \le u_i^F T_{i,j}, \qquad u_i^F := N - i + 1, \qquad \forall 1 \le j < i \le N.
\end{equation}
Moreover, the layer-linking inequalities~\eqref{eq:M0_layer_link} tighten to $1 - j(1 - T_{i,j}) \le L_i - L_j \le 1 + (i - 1)(1 - T_{i,j})$ for all $1 \le j < i \le N$.
\item[(iii)] \textbf{Exact base reaction linearization and global torque identity:} At the table contact $(1, 0)$, the vertical reaction force is fixed identically to $F_{1,0} = N$. Consequently, the bilinear equation $\tau_{1,0} = F_{1,0} P_{1,0}$ reduces to the exact linear identity
\begin{equation}
\label{eq:base_linearization}
\tau_{1,0} = N P_{1,0}, \qquad X_1 - w \le P_{1,0} \le X_1 + w, \qquad P_{1,0} \le -\epsilon,
\end{equation}
and global rotational equilibrium reduces to the linear identity
\begin{equation}
\label{eq:global_torque_identity}
\sum_{i=1}^N X_i = N P_{1,0}.
\end{equation}
\end{enumerate}
\end{lemma}

The proof, which uses only the canonical layer order and summation of the equilibrium equations, is given in Appendix~\ref{app:upper_dag_proof}.

Lemma~\ref{lem:upper_dag_base} reduces the candidate contacts from $N^2$ to the $N(N - 1)/2$ internal pairs $(i, j)$ with $1 \le j < i \le N$ and removes the bilinear product at the base.

\subsection{McCormick Relaxation Gap}
\label{subsec:McCormickGap}

Every constraint of $\mathcal{M}_0$ is linear except the bilinear equations $\tau_{i,j} = F_{i,j} P_{i,j}$ in~\eqref{eq:bilinear_torque}. Over a box $[\underline{F}_{i,j}, \overline{F}_{i,j}] \times [\underline{P}_{i,j}, \overline{P}_{i,j}]$ with widths $\Delta F_{i,j}$ and $\Delta P_{i,j}$, spatial branch-and-bound relaxes each of them by the McCormick envelope~\cite{mccormick1976}, the convex hull of the saddle surface $\tau_{i,j} = F_{i,j} P_{i,j}$ over the box~\cite{alkhayyal1983}. Its largest vertical thickness, $\frac{1}{2} \Delta F_{i,j} \Delta P_{i,j}$, is attained at the center of the box, where each of the two envelopes deviates from the bilinear term by $\frac{1}{4} \Delta F_{i,j} \Delta P_{i,j}$~\cite{luedtke2012,boland2017}. On the box $[0, N] \times [-N, N]$ of $\mathcal{M}_0$, this thickness is $N^2$ at the root node, and after Lemma~\ref{lem:upper_dag_base}, which removes the table contact, the internal pairs with $F_{i,j} \in [0, N - i + 1]$ retain $(N - i + 1) N = O(N^2)$. In addition, the big-$M$ rows~\eqref{eq:M0_bigm_i}--\eqref{eq:M0_bigm_j} let $P_{i,j}$ leave $[X_i - w, X_i + w]$ by up to $(2N + 1)(1 - T_{i,j})$ at fractional $T_{i,j}$~\cite{ceria1999,vielma2015}, and the integer layer variables $L_i$ create branching symmetry. Sections~\ref{sec:Structural} and~\ref{sec:SpinalCert} address these weaknesses by structural reductions, by an exact decomposition over layer profiles, which eliminates $L_i$ and $b_i$ and contracts the boxes, and by spinal cuts. A different remedy would measure each contact moment about the centers of the two blocks, so that the only bilinear factor is the center offset in $[-2w, 2w]$; we do not use this reformulation here.

\section{Structural Results for Single-Root Stacks}
\label{sec:Structural}

This section establishes structural properties of feasible and of optimal single-root stacks in the notation of Section~\ref{sec:Formulation}, with contact variables $(T_{ij}, F_{ij}, P_{ij}, \tau_{ij})$ for $1 \le j < i \le N$; Lemma~\ref{lem:upper_dag_base} gives the upper-triangular structure and the capacities $F_{ij} \le (N - i + 1) T_{ij}$.

\subsection{Support Degrees and Envelopes}
\label{subsec:support_bound}

The following elementary bounds hold for every feasible stack, optimal or not.

\begin{lemma}[Support Degrees and Envelopes]
\label{lem:support_bound}
Let $\eta \in [0, 1)$ and consider any feasible single-root $N$-block $\eta$-stable stack.
\begin{enumerate}
\item[(i)] \textbf{Support degrees:} every block $i \ge 2$ rests on at most three blocks, $\sum_{j < i} T_{ij} \le 3$, and at most three blocks rest on any block $j$, $\sum_{i > j} T_{ij} \le 3$. When $\eta > 0$, both bounds improve to two.
\item[(ii)] \textbf{Base envelope:} $X_1 \le \frac{N - 1}{N} w - \epsilon$.
\item[(iii)] \textbf{Coordinate envelope:} $|X_i - X_1| \le 2w (L_i - 1)$ for every block $i$.
\end{enumerate}
\end{lemma}

The proof, which counts pairwise unit-spaced support centers in an interval of length $4w$, applies rotational equilibrium to block~$1$, and follows support paths downward, is given in Appendix~\ref{app:support_bound_proof}.

\subsection{Optimal Table-Edge Alignment}
\label{subsec:translation}

We first remove the continuous horizontal translation invariance inherent in the table-support constraint $P_{1,0} \le -\epsilon$.

\begin{lemma}[Optimal Table-Edge Alignment]
\label{lem:s1_translation}
For every $N \ge 1$ and every $\eta \in [0, 1)$, every globally optimal single-root $N$-block $\eta$-stable stack satisfies
\begin{equation}
\label{eq:s1_alignment}
P_{1,0} = -\epsilon, \qquad \tau_{1,0} = -N\epsilon, \qquad \sum_{i=1}^{N} X_i = -N\epsilon.
\end{equation}
\end{lemma}

The proof, which translates an optimal stack rightward until its table contact reaches $-\epsilon$, is given in Appendix~\ref{app:s1_translation_proof}.

\subsection{Layer-Width Doubling and Structural Variable Elimination}
\label{subsec:layer_doubling}

Let $H := L_N$ be the total number of occupied layers of a feasible single-root stack, and let $\mathcal{B}_\ell := \{i \in \{1, \dots, N\} : L_i = \ell\}$ and $n_\ell := |\mathcal{B}_\ell|$ denote the block set and width of layer $\ell \in \{1, \dots, H\}$, with $n_1 = 1$. The vector $\mathbf{n} = (n_1, \dots, n_H)$ is called the \emph{layer profile} of the stack. The next theorem establishes unconditional upper bounds on the growth of layer widths that hold for every feasible single-root stack without requiring optimality or bounded support degree.

\begin{theorem}[Layer-Width Doubling and Expansion]
\label{thm:layer_doubling}
Let $\eta \in [0, 1)$, and consider any feasible single-root $N$-block $\eta$-stable stack with layer profile $(n_1, \dots, n_H)$. For each $\ell \in \{1, \dots, H\}$, define the extreme block centers $y_{\min}(\ell) := \min_{j \in \mathcal{B}_\ell} X_j$ and $y_{\max}(\ell) := \max_{j \in \mathcal{B}_\ell} X_j$. Then the following properties hold:
\begin{enumerate}
\item[(i)] \textbf{Strict endpoint exclusion:} For every $\ell \in \{1, \dots, H-1\}$,
\begin{equation}
\label{eq:endpoint_exclusion}
y_{\min}(\ell + 1) > y_{\min}(\ell) - 2w, \qquad y_{\max}(\ell + 1) < y_{\max}(\ell) + 2w.
\end{equation}
\item[(ii)] \textbf{Width doubling and linear expansion:}
\begin{equation}
\label{eq:layer_doubling}
n_{\ell + 1} \le 2 n_\ell \quad (\forall \ell \in \{1, \dots, H - 1\}),
\qquad
n_\ell \le 2\ell - 2 \quad (\forall \ell \in \{2, \dots, H\}).
\end{equation}
\item[(iii)] \textbf{Singleton second-layer refinement:} If $n_2 = 1$, then $n_\ell \le 2\ell - 4$ for all $\ell \in \{3, \dots, H\}$.
\end{enumerate}
\end{theorem}

The proof, which excludes a point contact at the end of a layer by moment balance, splits each layer into clusters of centers separated by gaps larger than $4w$, bounds the blocks of the next layer supported by each cluster, and inducts over the layers, is given in Appendix~\ref{app:layer_doubling_proof}.

For $\eta \ge 1/2$, part~(i) leaves one block in every layer: two blocks of one layer have centers at least $2W = 1$ apart, whereas by~\eqref{eq:endpoint_exclusion} and induction from $n_1 = 1$ the centers of every layer lie in an open interval of length $4w = 2(1 - \eta) \le 1$. Single-root stacks are then columns, whose maximum overhang is the harmonic value $w H_N$ with $H_N = \sum_{k=1}^{N} 1/k$; this case serves as an analytical check of the model, and our computations concern $\eta = 0$ only.

Theorem~\ref{thm:layer_doubling} immediately restricts the second layer to $n_2 \in \{1, 2\}$. When $n_2 = 2$, blocks~$2$ and~$3$ both rest only on block~$1$, so $1 \le F_{2,1}, F_{3,1} \le N - 2$, and strict endpoint exclusion~\eqref{eq:endpoint_exclusion} together with $X_2 - X_3 \ge 1$ gives $X_2 > X_1 > X_3$.

\subsection{Base-Singleton Reduction and Inductive Elimination of $n_2 = 1$}
\label{subsec:base_singleton}

A stack with $n_2 = 1$ has overhang at most the optimal overhang of $N - 1$ blocks plus $w / N$, so a certified bound for $N - 1$ eliminates the case $n_2 = 1$ inductively.

\begin{theorem}[Base-Singleton Reduction]
\label{thm:base_singleton}
Let $N \ge 2$ and $\eta \in [0, 1)$, and let $z$ be the overhang of any feasible single-root $N$-block $\eta$-stable stack satisfying $n_2 = 1$. Then
\begin{equation}
\label{eq:base_singleton_bound}
z \le z^*(N - 1, \eta) + \frac{w}{N} = z^*(N - 1, \eta) + \frac{1 - \eta}{2N}.
\end{equation}
\end{theorem}

The proof, given in Appendix~\ref{app:base_singleton_proof}, removes block~$1$, shifts the upper substack $\{2, \dots, N\}$ leftward by $P_{2,1} + \epsilon \le w / N$ onto the table, and compares the overhang of block~$1$ itself with the single-column harmonic stack of $N - 1$ blocks.

\begin{corollary}[Inductive Elimination of $n_2 = 1$, Strengthening \textsf{S2+}]
\label{cor:s2p_elimination}
Let $N \ge 3$ and $\eta \in [0, 1)$. Suppose that $\overline{z}(N - 1, \eta)$ is a certified upper bound satisfying $z^*(N - 1, \eta) \le \overline{z}(N - 1, \eta)$ and $\underline{z}(N, \eta)$ is the overhang of a known feasible single-root $N$-block $\eta$-stable stack such that
\begin{equation}
\label{eq:s2p_condition}
\overline{z}(N - 1, \eta) + \frac{w}{N} < \underline{z}(N, \eta).
\end{equation}
Then every globally optimal single-root $N$-block $\eta$-stable stack satisfies $n_2 = 2$, and consequently:
\begin{enumerate}
\item[(i)] $b_2 = 0$, $L_2 = L_3 = 2$, and $L_j \le j - 1$ for all $j \in \{3, \dots, N\}$;
\item[(ii)] $T_{i,3} = F_{i,3} = 0$ for all $i \ge 8$;
\item[(iii)] $T_{2,1} = T_{3,1} = 1$, $T_{3,2} = 0$, $1 \le F_{2,1} \le N - 2$, $1 \le F_{3,1} \le N - 2$, and $X_2 > X_1 > X_3$.
\end{enumerate}
\end{corollary}

The proof, which combines Theorem~\ref{thm:base_singleton} and~\eqref{eq:s2p_condition} with the layer-width bounds and the endpoint exclusion of Theorem~\ref{thm:layer_doubling}, is given in Appendix~\ref{app:s2p_elimination_proof}.

The reduction uses a certified bound $\overline{z}(N - 1, \eta)$ and no spinal assumption, and the certificate for $N - 1$ never depends on $N$. For $\eta = 0$, Table~\ref{tab:constants} verifies~\eqref{eq:s2p_condition} for every $N \in \{5, \dots, 16\}$ along one chain of our certificates, with margins of at least $1.79 \times 10^{-2}$.

\begin{remark}
\label{rem:s2p_scope}
Unlike Theorem~\ref{thm:layer_doubling}, which restricts every feasible single-root stack, Corollary~\ref{cor:s2p_elimination} restricts only \emph{globally optimal} stacks. Consequently, strengthening \textsf{S2+} is valid for maximum-overhang certification under~\eqref{eq:s2p_condition}, whereas feasibility investigations over all balanced stacks (such as testing the existence of essential three-support configurations) must retain profiles with $n_2 = 1$. A non-spinal analogue of this reduction, valid for non-optimal stacks, is given in Corollary~\ref{cor:nonspinal_exclusion}.
\end{remark}

\subsection{Base-Block Shift and Base Right-Corner Tightness}
\label{subsec:base_right_corner}

The moment $\tau_{2,1} = F_{2,1} P_{2,1}$ carries the heaviest load in the stack, so a narrow domain for $P_{2,1}$ reduces the McCormick error at the root node. We first normalize every feasible stack, without loss of overhang outside the base, so that the table reaction acts at $-\epsilon$ and block~$2$ rests on the right corner of block~$1$.

\begin{lemma}[Base-Block Shift and Global Translation]
\label{lem:wlog_base_shift}
Let $N \ge 2$ and $\eta \in [0, 1)$, and let $X$ be any feasible single-root $N$-block $\eta$-stable stack (not necessarily optimal) with layer profile $\mathbf{n} = (n_1, \dots, n_H)$, indexed canonically, so that $\mathcal{B}_1 = \{1\}$ and, by Theorem~\ref{thm:layer_doubling}, either $\mathcal{B}_2 = \{2\}$ or $\mathcal{B}_2 = \{2, 3\}$ with $X_2 - X_3 \ge 1$. Let $s_l$ denote the rightmost block of layer~$l$, and set $\Delta_0 := -\epsilon - P_{1,0} \ge 0$ and $\delta_1 := X_1 + w - P_{2,1} \ge 0$. Then there exists a feasible single-root $N$-block $\eta$-stable stack $\tilde{X}$, with contact points $\tilde{P}$, such that:
\begin{enumerate}
\item[(i)] $\tilde{X}$ has the same support topology $T$, contact forces $F$, and layer profile as $X$; every block $i \in \{2, \dots, N\}$ is shifted rightward by $\Delta := \Delta_0 + \delta_1 / N$, that is, $\tilde{X}_i = X_i + \Delta$, while $\tilde{X}_1 = X_1 + \Delta - \delta_1$;
\item[(ii)] $\tilde{X}$ satisfies \textsf{S1} ($\tilde{P}_{1,0} = -\epsilon$, hence $\sum_{i=1}^N \tilde{X}_i = -N\epsilon$) and \textsf{S4} ($\tilde{P}_{2,1} = \tilde{X}_1 + w$);
\item[(iii)] $\tilde{X}_1 \ge -W$, and, if $n_2 = 2$ and $N \ge 3$, $\tilde{X}_1 \le \frac{N - 3}{N} w - \epsilon$; moreover, $\tilde{X}$ satisfies the hypothesis $P_{1,0} = -\epsilon$ of Corollary~\ref{cor:energy_layer_ub};
\item[(iv)] every block $i \in \{2, \dots, N\}$ has overhang $\tilde{X}_i + W = X_i + W + \Delta \ge X_i + W$. Whenever block~$1$ does not attain the maximum overhang of $X$, the set of blocks attaining the maximum overhang is the same in $X$ and $\tilde{X}$, $z(\tilde{X}) = z(X) + \Delta$, the principal block and the spinal or non-spinal status are preserved, and for every maximum-overhang layer $\ell^* \ge 2$ the dominance inequalities $\tilde{X}_{s_l} \le \tilde{X}_{s_{\ell^*}}$ hold for all $l \in \{1, \dots, H\}$.
\end{enumerate}
\end{lemma}

The proof constructs $\tilde{X}$ by a global translation, a leftward shift of the base block, and a final translation, all of which keep the support topology and the contact forces unchanged; it is given in Appendix~\ref{app:wlog_base_shift_proof}.

The next theorem is the optimality statement behind strengthening~\textsf{S4}.

\begin{theorem}[Base Right-Corner Tightness, Strengthening \textsf{S4}]
\label{thm:base_right_corner}
Let $\eta \in [0, 1)$. For every $N \ge 3$, every globally optimal single-root $N$-block $\eta$-stable stack satisfies
\begin{equation}
\label{eq:s4_right_corner}
P_{2,1} = X_1 + w.
\end{equation}
For $N = 2$, at least one globally optimal single-root $\eta$-stable stack satisfies~\eqref{eq:s4_right_corner}.
\end{theorem}

The proof, which checks $N = 2$ directly and, for $N \ge 3$, shows that the optimal overhang is attained above the base block and then applies Lemma~\ref{lem:wlog_base_shift} to any optimal stack with $P_{2,1} < X_1 + w$, is given in Appendix~\ref{app:base_right_corner_proof}.

Combining Theorem~\ref{thm:base_right_corner} with $n_2 = 2$ tightens the upper bound on the base coordinate $X_1$ by $2w / N$ and contracts the domain of $P_{2,1}$ to an $N$-independent interval of width strictly less than one.

\begin{corollary}[Tightened Base Coordinate and Contact-Point Bounds]
\label{cor:base_ub_tight}
Let $N \ge 3$ and $\eta \in [0, 1)$. In any globally optimal single-root $N$-block $\eta$-stable stack with $n_2 = 2$, the base coordinate satisfies
\begin{equation}
\label{eq:s4_x1_ub}
X_1 \le \frac{N - 3}{N} w - \epsilon,
\end{equation}
and the base right contact point satisfies
\begin{equation}
\label{eq:s4_p21_box}
-\epsilon \le P_{2,1} \le \frac{2N - 3}{N} w - \epsilon.
\end{equation}
\end{corollary}

\begin{proof}
For an optimal stack, Lemma~\ref{lem:s1_translation} gives $P_{1,0} = -\epsilon$ and Theorem~\ref{thm:base_right_corner} gives $P_{2,1} = X_1 + w$, so $\Delta_0 = \delta_1 = 0$ in Lemma~\ref{lem:wlog_base_shift} and the normalized stack $\tilde{X}$ coincides with $X$. Lemma~\ref{lem:wlog_base_shift}(iii) therefore yields~\eqref{eq:s4_x1_ub} together with $X_1 \ge -W$. Substituting $-W \le X_1 \le \frac{N - 3}{N} w - \epsilon$ into $P_{2,1} = X_1 + w$ and using $-W + w = -\epsilon$ gives~\eqref{eq:s4_p21_box}, whose interval width is $\frac{2N - 3}{N} w = \big(1 - \frac{3}{2N}\big)(1 - \eta) < 1$.
\end{proof}

\begin{remark}[Optimality-Based Reductions for Non-Optimal Stacks]
\label{rem:s4_wlog_upgrade}
Lemma~\ref{lem:s1_translation}, Theorem~\ref{thm:base_right_corner}, and Corollary~\ref{cor:base_ub_tight} are statements about optimal stacks. Lemma~\ref{lem:wlog_base_shift} shows that the corresponding constraints (\textsf{S1}, \textsf{S4}, and for $n_2 = 2$ the bound $X_1 \le \frac{N - 3}{N} w - \epsilon$) can nevertheless be imposed on every feasible stack without loss of generality: the normalized stack has the same support topology, forces, and layer profile, and no block outside the base loses overhang. This is what allows these reductions inside the non-spinal exclusion check of Corollary~\ref{cor:nonspinal_exclusion}, whose subproblems must represent non-optimal stacks.
\end{remark}

Strengthening \textsf{S4} is imposed only at the base block; Appendix~\ref{app:s4_scope} discusses the left corner and the intermediate singleton layers.

\subsection{Cut-Set Difference-RLT Inequalities}
\label{subsec:c2_cuts}

Pairwise McCormick envelopes ignore that the forces crossing a horizontal cut between adjacent layers sum to the weight above the cut. Combining force and moment conservation across such cuts with same-layer non-overlap aggregates several bilinear products into one product with a narrow continuous factor, which yields sparse valid inequalities in the spirit of the reformulation-linearization technique~\cite{sherali1990} without auxiliary variables: family~\textsf{C2} for singleton layers that support two blocks, and family~\textsf{C2-Layer} for two-block layers.

\begin{proposition}[Singleton-Layer Difference-RLT Cuts, Family \textsf{C2}]
\label{prop:c2_singleton}
Let $\eta \in [0, 1)$, and consider any feasible single-root $N$-block $\eta$-stable stack with layer profile $(n_1, \dots, n_H)$. Suppose that some layer $\ell \in \{1, \dots, H - 1\}$ is a singleton $\mathcal{B}_\ell = \{u\}$ whose immediately higher layer contains two blocks $\mathcal{B}_{\ell + 1} = \{v_1, v_2\}$ with $v_1 < v_2$. Let $M_\ell := \sum_{t = \ell + 1}^{H} n_t \ge 2$ denote the total weight above layer~$\ell$, define the total downward moment leaving block~$u$ by
\[
\tau_u^{\mathrm{out}} :=
\begin{cases}
\tau_{1,0}, & \text{if } \ell = 1 \text{ (so } u = 1\text{)}, \\
\displaystyle\sum_{j \in \mathcal{B}_{\ell - 1}} \tau_{u, j}, & \text{if } \ell \ge 2,
\end{cases}
\]
and set $Q_u := \tau_u^{\mathrm{out}} - X_u - M_\ell P_{v_2, u}$, $\delta_u := P_{v_1, u} - P_{v_2, u}$, $d_L := 1 - 2w = \eta$, and $d_U := 2w = 1 - \eta$. Then $d_L < d_U$, and the following four linear inequalities hold unconditionally:
\begin{subequations}
\label{eq:c2_facets}
\begin{align}
Q_u &\ge \delta_u + d_L F_{v_1, u} - d_L, \label{eq:c2_a} \\
Q_u &\ge (M_\ell - 1) \delta_u + d_U F_{v_1, u} - (M_\ell - 1) d_U, \label{eq:c2_b} \\
Q_u &\le (M_\ell - 1) \delta_u + d_L F_{v_1, u} - (M_\ell - 1) d_L, \label{eq:c2_c} \\
Q_u &\le \delta_u + d_U F_{v_1, u} - d_U. \label{eq:c2_d}
\end{align}
\end{subequations}
\end{proposition}

The proof is given in Appendix~\ref{app:c2_proofs}.

At the base under \textsf{S2+} ($n_2 = 2$, $u = 1$, $v_1 = 2$, $v_2 = 3$, $M_1 = N - 1$), Proposition~\ref{prop:c2_singleton} applies inside the monolithic formulation. The width $1 - 2\eta$ of the domain of $\delta_1 = P_{2,1} - P_{3,1}$ does not depend on $N$, whereas $P_{2,1}$ and $P_{3,1}$ themselves range over $O(N)$; with $\tau_{1,0} = 0$ (Lemma~\ref{lem:s1_translation}) and $P_{2,1} = X_1 + w$ (Theorem~\ref{thm:base_right_corner}), system~\eqref{eq:c2_facets} reduces to four sparse rows in $X_1$, $P_{3,1}$, and $F_{2,1}$.

In the hardest layer profiles for $N \in \{10, \dots, 13\}$, most intermediate layers contain two blocks rather than one, so Proposition~\ref{prop:c2_singleton} does not trigger above the base. The following proposition extends the cut-set difference-RLT construction to every non-top two-block layer.

\begin{proposition}[Two-Block Layer Cut-Set RLT Cuts, Family \textsf{C2-Layer}]
\label{prop:c2_layer}
Let $\eta \in [0, 1)$, and consider any feasible single-root $N$-block $\eta$-stable stack with fixed layer profile $(n_1, \dots, n_H)$. Let $\ell \in \{2, \dots, H - 1\}$ be any non-top layer with $n_\ell = 2$, write $\mathcal{B}_\ell = \{u_1, u_2\}$ with $u_1 < u_2$, and set $\mathcal{B}_{>\ell} := \bigcup_{t = \ell + 1}^{H} \mathcal{B}_t$ and $M_\ell := |\mathcal{B}_{>\ell}| = \sum_{t = \ell + 1}^{H} n_t \ge 1$. Define the aggregated vertical forces crossing into layer~$\ell$, the layer spacing, and the upper-cut relative moment by
\[
F_{\ell, 1} := \sum_{k \in \mathcal{B}_{\ell + 1}} F_{k, u_1},
\qquad
\Delta_\ell := X_{u_1} - X_{u_2},
\qquad
W_\ell := \sum_{b \in \mathcal{B}_{>\ell}} X_b - M_\ell X_{u_2}.
\]
Let $s_0(\ell) := \max\{t \in \{1, \dots, \ell - 1\} : n_t = 1\}$ be the highest singleton layer below~$\ell$, and let $\overline{\Delta}_\ell := \min\big(4w(\ell - s_0(\ell)),\; x_{u_1}^{\mathrm{ub}} - x_{u_2}^{\mathrm{lb}}\big)$ for any valid coordinate bounds $X_{u_1} \le x_{u_1}^{\mathrm{ub}}$ and $X_{u_2} \ge x_{u_2}^{\mathrm{lb}}$. Then $1 \le \Delta_\ell \le \overline{\Delta}_\ell$ and the following four linear inequalities hold unconditionally:
\begin{subequations}
\label{eq:c2layer_facets}
\begin{align}
W_\ell &\ge F_{\ell, 1} - M_\ell w, \label{eq:c2l_a} \\
W_\ell &\ge M_\ell \Delta_\ell + \overline{\Delta}_\ell F_{\ell, 1} - M_\ell(\overline{\Delta}_\ell + w), \label{eq:c2l_b} \\
W_\ell &\le M_\ell \Delta_\ell + F_{\ell, 1} - M_\ell(1 - w), \label{eq:c2l_c} \\
W_\ell &\le \overline{\Delta}_\ell F_{\ell, 1} + M_\ell w. \label{eq:c2l_d}
\end{align}
\end{subequations}
\end{proposition}

The proof is given in Appendix~\ref{app:c2_proofs}.

Remark~\ref{rem:c2_discussion} in Appendix~\ref{app:c2_proofs} discusses the structure and scope of both families.

\subsection{Second-Moment Energy Inequality and Analytical Profile Screening}
\label{subsec:energy_bound}

The layer-width bounds are purely geometric and ignore the equilibrium cost of holding wide layers at low height. We now derive a second-moment identity from the equilibrium equations, which yields a quadratic energy bound, an analytical infeasibility criterion for wide low profiles, and a closed-form bound on the overhang of each layer.

\begin{theorem}[Global Second-Moment Energy Inequality]
\label{thm:second_moment_energy}
Let $\eta \in [0, 1)$, and consider any feasible single-root $N$-block $\eta$-stable stack with layer profile $\mathbf{n} = (n_1, \dots, n_H)$ satisfying $P_{1,0} = -\epsilon$. Define the shifted coordinates $x_b := X_b + \epsilon$ for $b \in \{1, \dots, N\}$. Then the exact second-moment identity
\begin{equation}
\label{eq:second_moment_identity}
\sum_{b=1}^{N} x_b^2 + N x_1^2 + \sum_{1 \le j < i \le N} F_{ij} (P_{ij} - X_i)^2
= \sum_{1 \le j < i \le N} F_{ij} (P_{ij} - X_j)^2
\end{equation}
holds, and consequently the block coordinates satisfy the convex quadratic inequality
\begin{equation}
\label{eq:energy_quad}
\sum_{b=1}^{N} (X_b + \epsilon)^2 + N (X_1 + \epsilon)^2
\le w^2 \sum_{b=1}^{N} (L_b - 1)
= w^2 \sum_{\ell=2}^{H} (\ell - 1) n_\ell
=: E(\mathbf{n}).
\end{equation}
\end{theorem}

The proof, which sums a moment identity for the squared contact points over all blocks and bounds every contact offset by $w$, is given in Appendix~\ref{app:second_moment_proof}.

Note that any feasible single-root stack (not necessarily optimal) can be translated rightward by $-\epsilon - P_{1,0} \ge 0$ to satisfy $P_{1,0} = -\epsilon$ without altering its layer profile or equilibrium. Combining Theorem~\ref{thm:second_moment_energy} with the minimum internal variance enforced by same-layer non-overlap immediately yields an unconditional feasibility criterion.

\begin{corollary}[Analytical Profile Infeasibility Criterion]
\label{cor:energy_infeasibility}
For any positive integer $m$, let $V_{\min}(m) := m(m^2 - 1) / 12$, and for any layer profile $\mathbf{n} = (n_1, \dots, n_H)$, define its energy margin
\begin{equation}
\label{eq:energy_margin}
\Delta E(\mathbf{n}) := E(\mathbf{n}) - \sum_{\ell=1}^{H} V_{\min}(n_\ell)
= w^2 \sum_{\ell=2}^{H} (\ell - 1) n_\ell - \sum_{\ell=1}^{H} \frac{n_\ell (n_\ell^2 - 1)}{12}.
\end{equation}
If $\Delta E(\mathbf{n}) < 0$, then no feasible single-root $N$-block $\eta$-stable stack with layer profile $\mathbf{n}$ exists.
\end{corollary}

The proof, which bounds the left-hand side of~\eqref{eq:energy_quad} from below by Lagrange's identity and the unit spacing within each layer, is given in Appendix~\ref{app:energy_infeasibility_proof}.

When $\Delta E(\mathbf{n}) \ge 0$, combining the global torque identity~\eqref{eq:global_torque_identity} under $P_{1,0} = -\epsilon$ with the layer-mean energy budget $N x_1^2 + \sum_{\ell=1}^{H} n_\ell \bar{x}_\ell^2 \le \Delta E(\mathbf{n})$, where $x_b := X_b + \epsilon$ and $\bar{x}_\ell$ is the mean of $x_b$ over layer~$\ell$ (Appendix~\ref{app:energy_infeasibility_proof}), provides a closed-form upper bound on the overhang achievable by each individual layer.

\begin{corollary}[Closed-Form Energy Upper Bound per Layer]
\label{cor:energy_layer_ub}
Let $\eta \in [0, 1)$, and consider any feasible single-root $N$-block $\eta$-stable stack with layer profile $\mathbf{n} = (n_1, \dots, n_H)$ satisfying $P_{1,0} = -\epsilon$ (in particular, any optimal stack by Lemma~\ref{lem:s1_translation}, or any stack normalized by Lemma~\ref{lem:wlog_base_shift}); hence $\Delta E(\mathbf{n}) \ge 0$ and $\sum_{b=1}^{N} (X_b + \epsilon) = 0$ by~\eqref{eq:global_torque_identity}. For any layer $\ell \in \{2, \dots, H\}$, set $m := n_\ell$ and
\[
A_m := m + \frac{m^2}{(N - m - 1) + \dfrac{1}{N + 1}}.
\]
Then the rightmost block $i_{\ell, 0}$ of layer~$\ell$ satisfies $X_{i_{\ell, 0}} + W \le \overline{z}^{\mathrm{E}}(\mathbf{n}, \ell)$, where
\begin{equation}
\label{eq:energy_layer_ub}
\overline{z}^{\mathrm{E}}(\mathbf{n}, \ell) :=
\begin{cases}
w + \dfrac{m - 1}{2} + \sqrt{\Delta E(\mathbf{n}) / A_m},
& \text{if } m = 1 \text{ or } 4 A_m \Delta E(\mathbf{n}) \le m^2, \\[6pt]
w + \sqrt{\big(\tfrac{1}{A_m} + \tfrac{m - 1}{m}\big)\big(\Delta E(\mathbf{n}) + \tfrac{m(m - 1)}{4}\big)},
& \text{otherwise}.
\end{cases}
\end{equation}
\end{corollary}

The proof, which combines Cauchy--Schwarz and Chebyshev sum inequalities with a planar elliptic bound, is given in Appendix~\ref{app:energy_proof}.

\begin{remark}[Screening Power and Role in the Formulation]
\label{rem:energy_filtering}
Corollary~\ref{cor:energy_infeasibility} removes between $5\%$ and $18\%$ of the profiles in $\mathcal{P}_2(N)$ for $N \in \{9, \dots, 16\}$ (Table~\ref{tab:profile_counts}), at negligible cost. Exhaustive solver checks over all profiles in $\mathcal{P}(N)$ for $N \in \{6, \dots, 9\}$ confirmed that no feasible profile is removed and that~\eqref{eq:energy_layer_ub} holds on all $437$ feasible profile and layer pairs. On tall near-optimal spinal profiles the energy margin is large (for example, $\Delta E(\mathbf{n}) = 6.0$ for $(1, 2, 2, 2, 3, 2, 1)$ at $N = 13$), so the layer bounds rarely prune inside the hardest subproblems. We therefore use Corollaries~\ref{cor:energy_infeasibility} and~\ref{cor:energy_layer_ub} as profile and candidate-layer pre-filters, and do not add the quadratic inequality~\eqref{eq:energy_quad} to the subproblems.
\end{remark}

\section{Layer-Profile Decomposition and Spinal Certificates}
\label{sec:SpinalCert}

This section turns the structural results of Section~\ref{sec:Structural} into a certification scheme. We first decompose the problem exactly over layer profiles, then bound every unique-support chain by the loaded-spine optimum of Paterson and Zwick, convert this bound into cuts, and finally derive a certified exclusion test for non-spinal stacks.

\subsection{Exact Layer-Profile Decomposition}
\label{subsec:profile_decomp}

In $\mathcal{M}_0$, every pair $(i, j)$ with $j < i$ is a potential contact, which gives $O(N^2)$ bilinear terms and big-$M$ links. Fixing the layer profile $\mathbf{n} = (n_1, \dots, n_H)$ eliminates all integer layer variables, fixes every non-adjacent contact pair to zero, reduces the number of bilinear pairs to $\sum_{\ell=1}^{H-1} n_\ell n_{\ell+1}$, and admits profile-specific valid inequalities.

Formally, define the set of structurally admissible layer profiles of size $N \ge 2$ by
\begin{align*}
\mathcal{P}(N) := \Big\{ \mathbf{n} \in \mathbb{Z}_{>0}^H :\;
& n_1 = 1,\; \textstyle\sum_{\ell=1}^{H} n_\ell = N,\; n_{\ell+1} \le 2 n_\ell \; (1 \le \ell < H), \\
& n_\ell \le 2\ell - 2 \; (2 \le \ell \le H), \\
& n_\ell \le 2\ell - 4 \; (3 \le \ell \le H) \text{ if } n_2 = 1 \Big\},
\end{align*}
its two-block base-layer subset $\mathcal{P}_2(N) := \{\mathbf{n} \in \mathcal{P}(N) : n_2 = 2\}$, and its energy-admissible subset
\[
\mathcal{P}_2^{\mathrm{E}}(N) := \big\{\mathbf{n} \in \mathcal{P}_2(N) : \Delta E(\mathbf{n}) \ge 0\big\}.
\]
Theorem~\ref{thm:layer_doubling} and Corollaries~\ref{cor:s2p_elimination} and~\ref{cor:energy_infeasibility} shrink the enumeration: for $N = 9, \dots, 16$, the number $2^{N-2}$ of unrestricted profiles with $n_1 = 1$ grows from $128$ to $16{,}384$, whereas $|\mathcal{P}_2^{\mathrm{E}}(N)|$ grows only from $18$ to $1{,}226$ (Table~\ref{tab:profile_counts}).

\begin{table}[t]
\centering
\caption{Progressive reduction in the number of layer profiles for $\eta = 0$ and $N \in \{9, \dots, 16\}$ from unrestricted single-root partitions ($2^{N-2}$) via layer-width doubling ($|\mathcal{P}(N)|$), inductive base-singleton elimination ($|\mathcal{P}_2(N)|$), and second-moment energy screening ($|\mathcal{P}_2^{\mathrm{E}}(N)|$).}
\label{tab:profile_counts}
\begin{tabular}{@{}rrrrrrr@{}}
\toprule
$N$ & $2^{N-2}$ & $|\mathcal{P}(N)|$ & $|\mathcal{P}_2(N)|$ & Energy-pruned & $|\mathcal{P}_2^{\mathrm{E}}(N)|$ & Reduction factor \\
\midrule
 9 &  128 &  50 &  22 &  4 &  18 &  $7.1\times$ \\
10 &  256 &  89 &  39 &  5 &  34 &  $7.5\times$ \\
11 &  512 & 159 &  70 &  5 &  65 &  $7.9\times$ \\
12 & 1024 & 285 & 126 & 12 & 114 &  $9.0\times$ \\
13 & 2048 & 510 & 225 & 13 & 212 &  $9.7\times$ \\
14 & 4096 & 913 & 403 & 25 & 378 & $10.8\times$ \\
15 & 8192 & 1636 & 723 & 41 & 682 & $12.0\times$ \\
16 & 16384 & 2933 & 1296 & 70 & 1226 & $13.4\times$ \\
\bottomrule
\end{tabular}
\end{table}

For any fixed profile $\mathbf{n} = (n_1, \dots, n_H) \in \mathcal{P}(N)$, canonical indexing uniquely determines the layer assignment $L_i$ of every block $i \in \{1, \dots, N\}$, the layer-jump vector $b_i = \mathbb{I}[L_{i+1} > L_i]$, and the ordered block set $\mathcal{B}_\ell = \{i_{\ell, 0} < i_{\ell, 1} < \dots < i_{\ell, n_\ell - 1}\}$ of each layer $\ell \in \{1, \dots, H\}$. Every non-adjacent pair $(i, j)$ with $L_i \neq L_j + 1$ satisfies $T_{ij} = F_{ij} = P_{ij} = \tau_{ij} = 0$. The next proposition establishes three profile-specific valid tightenings.

\begin{proposition}[Profile-Specific Bounds, Non-Crossing Cuts, and Singleton Bottlenecks]
\label{prop:profile_tightening}
Let $\eta \in [0, 1)$ and fix a layer profile $\mathbf{n} = (n_1, \dots, n_H) \in \mathcal{P}(N)$. Let $x_1^{\mathrm{lb}} \le X_1 \le x_1^{\mathrm{ub}}$ be valid bounds on the base block center. For optimal stacks with $N \ge 3$ and $n_2 = 2$, and without loss of generality for any feasible stack with $n_2 = 2$ by Lemma~\ref{lem:wlog_base_shift}, one may take $x_1^{\mathrm{lb}} = -W$ and $x_1^{\mathrm{ub}} = \frac{N - 3}{N} w - \epsilon$. Then every feasible stack with layer profile $\mathbf{n}$ satisfies:
\begin{enumerate}
\item[(i)] \textbf{Layer-rank coordinate bounds:} For every layer $\ell \in \{2, \dots, H\}$ and every block $i = i_{\ell, r} \in \mathcal{B}_\ell$ of within-layer rank $r \in \{0, \dots, n_\ell - 1\}$,
\begin{equation}
\label{eq:rank_box_rel}
X_1 - 2w(\ell - 1) + (n_\ell - 1 - r) \le X_i \le X_1 + 2w(\ell - 1) - r,
\end{equation}
and consequently $x_i^{\mathrm{lb}} \le X_i \le x_i^{\mathrm{ub}}$, where
\begin{equation}
\label{eq:rank_box_abs}
x_i^{\mathrm{lb}} := x_1^{\mathrm{lb}} - 2w(\ell - 1) + (n_\ell - 1 - r),
\qquad
x_i^{\mathrm{ub}} := x_1^{\mathrm{ub}} + 2w(\ell - 1) - r.
\end{equation}
Moreover, for every adjacent pair $(i, j) \in \mathcal{B}_{\ell+1} \times \mathcal{B}_\ell$, any active contact ($T_{ij} = 1$) satisfies $P_{ij} \in [p_{ij}^{\mathrm{lb}}, p_{ij}^{\mathrm{ub}}]$ with $p_{ij}^{\mathrm{lb}} := \max(x_i^{\mathrm{lb}}, x_j^{\mathrm{lb}}) - w$ and $p_{ij}^{\mathrm{ub}} := \min(x_i^{\mathrm{ub}}, x_j^{\mathrm{ub}}) + w$; when $T_{ij} = 0$ is permitted, $P_{ij}$ lies in $[\min(0, p_{ij}^{\mathrm{lb}}), \max(0, p_{ij}^{\mathrm{ub}})]$, inducing valid interval bounds on $\tau_{ij} = F_{ij} P_{ij}$ from the extremal products of the bounds on $F_{ij}$ and $P_{ij}$.
\item[(ii)] \textbf{Planar non-crossing inequalities:} For every $\ell \in \{1, \dots, H - 1\}$, any $i_1 < i_2$ in $\mathcal{B}_{\ell+1}$, and any $j_1 < j_2$ in $\mathcal{B}_\ell$,
\begin{equation}
\label{eq:planar_noncrossing}
T_{i_1, j_2} + T_{i_2, j_1} \le 1
\qquad\text{whenever } (i_2 - i_1) + (j_2 - j_1) \ge \theta(\eta),
\end{equation}
where $\theta(0) = 3$ and $\theta(\eta) = 2$ for $\eta > 0$.
\item[(iii)] \textbf{Singleton bottleneck reductions:} For every singleton layer $\mathcal{B}_\ell = \{u\}$ ($\ell < H$) with total upper weight $M_\ell = \sum_{t = \ell + 1}^{H} n_t$, every block $v \in \mathcal{B}_{\ell+1}$ satisfies $T_{v, u} = 1$. If $n_{\ell+1} = 1$ with $\mathcal{B}_{\ell+1} = \{v\}$, then $F_{v, u} = M_\ell$ is constant (linearizing $\tau_{v, u}$). If $n_{\ell+1} = 2$ with $\mathcal{B}_{\ell+1} = \{v_1, v_2\}$ ($v_1 < v_2$), then $X_{v_1} \ge X_u \ge X_{v_2}$ and the four inequalities~\eqref{eq:c2_facets} of Proposition~\ref{prop:c2_singleton} hold.
\end{enumerate}
\end{proposition}

The proof, which combines Theorem~\ref{thm:layer_doubling}(i) with same-layer spacing and contact-core overlap, is given in Appendix~\ref{app:profile_tightening_proof}.

We now state the global certification theorem for the cutoff-driven layer-profile decomposition. Given a known feasible overhang $\underline{z} > 0$ and a target cutoff $\gamma \ge \underline{z}$, such as the tight cutoff $\gamma_N$ used in Section~\ref{sec:Computational}, define for each energy-admissible profile $\mathbf{n} \in \mathcal{P}_2^{\mathrm{E}}(N)$ the candidate rightmost-block index set
\begin{equation}
\label{eq:candidate_set}
\mathcal{C}(\mathbf{n}, \gamma)
:= \Big\{ i_{\ell, 0} :\; \ell \in \{2, \dots, H\},\;
\min\!\big(x_1^{\mathrm{ub}} + W + 2w(\ell - 1),\; \overline{z}^{\mathrm{E}}(\mathbf{n}, \ell)\big) > \gamma \Big\},
\end{equation}
where $x_1^{\mathrm{ub}} = \frac{N - 3}{N} w - \epsilon$ by Corollary~\ref{cor:base_ub_tight} and $\overline{z}^{\mathrm{E}}(\mathbf{n}, \ell)$ is the closed-form energy upper bound~\eqref{eq:energy_layer_ub} of Corollary~\ref{cor:energy_layer_ub}. Let $\mathcal{F}(\mathbf{n})$ denote the feasible set of the strengthened profile subproblem obtained by combining the unconditional base formulation (Lemma~\ref{lem:support_bound}) with table-edge alignment (Lemma~\ref{lem:s1_translation}), base-layer fixation $n_2 = 2$ (Corollary~\ref{cor:s2p_elimination}), base right-corner tightness (Theorem~\ref{thm:base_right_corner} and Corollary~\ref{cor:base_ub_tight}), cut-set difference-RLT inequalities (Propositions~\ref{prop:c2_singleton} and~\ref{prop:c2_layer}), and profile-specific layer-rank boxes, planar non-crossing cuts, and singleton bottleneck reductions (Proposition~\ref{prop:profile_tightening}). This is the constraint set solved in Section~\ref{sec:Computational}. Whenever $\mathcal{C}(\mathbf{n}, \gamma) \neq \emptyset$, define the profile subproblem
\begin{equation}
\label{eq:profile_subproblem}
\mathcal{S}(\mathbf{n}, \gamma): \qquad
\max \left\{ \max_{i \in \mathcal{C}(\mathbf{n}, \gamma)} (X_i + W) :\; \text{constraints of } \mathcal{F}(\mathbf{n}) \right\},
\end{equation}
solved with objective cutoff $\gamma$. We say that $\mathcal{S}(\mathbf{n}, \gamma)$ is \emph{closed} if it is proved that no feasible point has objective strictly greater than $\gamma$. A solver reports such a proof either as termination by the cutoff or, because the cutoff may enter presolve as a constraint on the objective, as infeasibility; in both cases it says nothing about the feasible points with objective at most $\gamma$, so a closed subproblem contributes the bound $\gamma$, never $-\infty$. Only the energy test of Corollary~\ref{cor:energy_infeasibility} proves a profile empty.

\begin{theorem}[Exactness of Cutoff Layer-Profile Decomposition]
\label{thm:profile_decomp}
Let $N \ge 3$ and $\eta \in [0, 1)$, and suppose that the inductive condition~\eqref{eq:s2p_condition} of Corollary~\ref{cor:s2p_elimination} holds. Let $\underline{z} \le z^*(N, \eta)$ be the overhang of a known feasible stack, and fix any cutoff $\gamma \ge \underline{z}$. For each energy-admissible profile $\mathbf{n} \in \mathcal{P}_2^{\mathrm{E}}(N)$, define a subproblem upper bound $\overline{u}(\mathbf{n})$ as follows:
\begin{enumerate}
\item[(a)] if $\mathcal{C}(\mathbf{n}, \gamma) = \emptyset$, set $\overline{u}(\mathbf{n}) := \gamma$;
\item[(b)] if $\mathcal{C}(\mathbf{n}, \gamma) \neq \emptyset$ and $\mathcal{S}(\mathbf{n}, \gamma)$ is closed, set $\overline{u}(\mathbf{n}) := \gamma$;
\item[(c)] otherwise, let $\overline{z}(\mathbf{n})$ be any valid upper bound on the objective of the feasible points of $\mathcal{S}(\mathbf{n}, \gamma)$ with objective above $\gamma$, such as a dual bound returned for it, and set $\overline{u}(\mathbf{n}) := \max(\gamma, \overline{z}(\mathbf{n}))$.
\end{enumerate}
Then
\begin{equation}
\label{eq:pd_global_bound}
z^*(N, \eta) \le \max_{\mathbf{n} \in \mathcal{P}_2^{\mathrm{E}}(N)} \overline{u}(\mathbf{n}).
\end{equation}
In particular, if every profile $\mathbf{n} \in \mathcal{P}_2^{\mathrm{E}}(N)$ falls into case~(a) or~(b), or into case~(c) with $\overline{z}(\mathbf{n}) \le \gamma$, then $z^*(N, \eta) \le \gamma$.
\end{theorem}

The proof, which places a globally optimal stack in the subproblem of its own profile, is given in Appendix~\ref{app:profile_decomp_proof}.

\subsection{Spinal Bound, Big-$M$ Spinal Cuts, and Certificate Decomposition}
\label{subsec:spinal_bound}

Profiles that contain a near-optimal spinal stack still generate large branch-and-bound trees when $\gamma$ is close to the optimum (Section~\ref{subsec:method_progression}). We therefore bound every unique-support chain of an arbitrary single-root $\eta$-stable stack by the loaded-spine optimum of Paterson and Zwick~\cite[Lemmas~3.4--3.5]{paterson2009overhang}, scaled by $2w = 1 - \eta$, and turn this bound into a sparse family of big-$M$ inequalities (family~\textsf{SP}): when the rightmost chain up to a candidate block rests on no non-rightmost block, it is a unique-support chain and the bound applies, so every configuration above the cutoff must contain such a non-rightmost support.

For any integer $K \ge 1$ and any real load $M > 0$, define the unconstrained loaded-spine value
\begin{equation}
\label{eq:pz_unconstrained}
\hat{S}_K(M) := \max_{u_1, \dots, u_{K-1} > 0}
\sum_{s=1}^{K} \left(1 - \frac{u_{s-1} + 1/2}{u_s}\right)
\qquad (u_0 := 0,\; u_K := M),
\end{equation}
with $\hat{S}_1(M) := 1 - 1 / (2M)$. For $M \ge K$, define the corresponding load-constrained value
\begin{equation}
\label{eq:pz_constrained}
S_K^*(M) := \max_{\substack{u_0 = 0,\; u_K = M \\ u_s - u_{s-1} \ge 1 \; (1 \le s \le K)}}
\; \sum_{s=1}^{K} \left(1 - \frac{u_{s-1} + 1/2}{u_s}\right).
\end{equation}
Clearly $S_K^*(M) \le \hat{S}_K(M)$ for all $M \ge K$, and $S_K^*(M)$ coincides with the maximum overhang of a $K$-block loaded spinal stack of total weight $M$ established in~\cite[Lemmas~3.4--3.5]{paterson2009overhang}. Set $S^*(N) := \max_{1 \le K \le N} S_K^*(N)$ and $\hat{S}(N) := \max_{1 \le K \le N} \hat{S}_K(N)$.

\begin{proposition}[Loaded-Spinal Chain Bound]
\label{prop:spinal_bound}
Let $\eta \in [0, 1)$, so $w = W - \epsilon = (1 - \eta)/2$ and $2w = 1 - \eta$.
\begin{enumerate}
\item[(i)] \textbf{Uniqueness and exact evaluation of $\hat{S}_K(M)$:} For every $K \ge 2$ and every $M > 0$, problem~\eqref{eq:pz_unconstrained} admits a unique global maximizer $(u_1^*, \dots, u_{K-1}^*)$, characterized by the recurrence of Paterson and Zwick~\cite[Lemma~3.5]{paterson2009overhang},
\begin{equation}
\label{eq:pz_recurrence}
u_{s+1} = \frac{u_s^2}{u_{s-1} + 1/2} \qquad (s = 1, \dots, K - 1, \quad u_0 = 0,\; u_K = M).
\end{equation}
Moreover, the terminal value $u_K(u_1)$ generated by~\eqref{eq:pz_recurrence} is continuous and strictly increasing on $(0, \infty)$ with range $(0, \infty)$, so $\hat{S}_K(M)$ is determined uniquely by one-dimensional bisection on $u_1$. Whenever $(u_1^*, \dots, u_{K-1}^*)$ satisfies $u_s^* - u_{s-1}^* \ge 1$ for all $s \in \{1, \dots, K\}$, one has $S_K^*(M) = \hat{S}_K(M)$.
\item[(ii)] \textbf{Inequality for arbitrary unique-support chains:} Consider any feasible single-root $N$-block $\eta$-stable stack (not necessarily spinal and not necessarily optimal). Suppose that blocks $b_p, b_{p+1}, \dots, b_\ell$ lie in consecutive layers $p, p+1, \dots, \ell$ ($2 \le p \le \ell \le H$), and that each $b_t$ ($t \in \{p, \dots, \ell\}$) is supported exclusively by $b_{t-1} \in \mathcal{B}_{t-1}$ (that is, $T_{b_t, j} = 0$ for all $j \in \mathcal{B}_{t-1} \setminus \{b_{t-1}\}$). Let $r := \ell - p + 1$ and $F_p := F_{b_p, b_{p-1}}$. Then
\begin{equation}
\label{eq:spinal_suffix_bound}
X_{b_\ell} - X_{b_{p-1}} \le 2w \, S_r^*(F_p) \le 2w \, \hat{S}_r(F_p).
\end{equation}
In particular, if $p = 2$ (so $b_1 = 1$ is the unique base block) or if $\ell = 1$ (with $b_1 = 1$), then the total overhang of block $b_\ell$ satisfies
\begin{equation}
\label{eq:spinal_full_bound}
X_{b_\ell} + W \le 2w \, S_\ell^*(N) \le 2w \, \hat{S}_\ell(N).
\end{equation}
\end{enumerate}
\end{proposition}

Part~(i) is a self-contained version of Lemma~3.5 of~\cite{paterson2009overhang}. Part~(ii) follows by reverse induction along the chain $b_p, \dots, b_\ell$ from the force and moment balance of its blocks, followed by Abel's summation by parts and the definitions~\eqref{eq:pz_unconstrained} and~\eqref{eq:pz_constrained}. Both proofs are given in Appendix~\ref{app:spinal_i_proof}.

Proposition~\ref{prop:spinal_bound}(ii) needs neither the table-edge equality $P_{1,0} = -\epsilon$ nor layer-rightmost blocks. Applied to the chain $(s_1, \dots, s_\ell)$ of layer-rightmost blocks of a profile, it yields the following cuts.

\begin{proposition}[Validity of Spinal-Bound Big-$M$ Cuts, Family \textsf{SP}]
\label{prop:sp_cuts}
Let $\eta \in [0, 1)$ and fix any layer profile $\mathbf{n} = (n_1, \dots, n_H) \in \mathcal{P}(N)$. For each layer $t \in \{1, \dots, H\}$, denote its rightmost block by $s_t := i_{t, 0}$. For any layer $\ell \in \{2, \dots, H\}$, define the set of intermediate layers having a multi-block predecessor by
\[
J_\ell := \{t \in \{3, \dots, \ell\} : n_{t-1} \ge 2\},
\]
let $\overline{x}_{s_\ell}$ be an upper bound on $X_{s_\ell}$, and set
\begin{equation}
\label{eq:sp_bigm}
M_\ell := \max\!\big(0,\; \overline{x}_{s_\ell} + W - 2w \, \hat{S}_\ell(N)\big).
\end{equation}
Then the linear inequality
\begin{equation}
\label{eq:sp_cut}
X_{s_\ell} + W \le 2w \, \hat{S}_\ell(N)
+ M_\ell \sum_{t \in J_\ell} \sum_{q = 1}^{n_{t-1} - 1} T_{s_t,\, i_{t-1, q}}
\end{equation}
holds for every feasible single-root $N$-block $\eta$-stable stack with layer profile $\mathbf{n}$ that satisfies $X_{s_\ell} \le \overline{x}_{s_\ell}$, independently of any objective cutoff. In particular, \eqref{eq:sp_cut} holds for every feasible stack when $\overline{x}_{s_\ell}$ is taken from the envelopes of Lemma~\ref{lem:support_bound}, and, by Lemma~\ref{lem:wlog_base_shift}, without loss of generality for every feasible stack with $n_2 = 2$ when $\overline{x}_{s_\ell}$ is the layer-rank box bound $x_{s_\ell}^{\mathrm{ub}}$ of Proposition~\ref{prop:profile_tightening}(i). Our implementation takes $\overline{x}_{s_\ell}$ from this box intersected with the envelopes of Lemma~\ref{lem:support_bound} and replaces $2w\,\hat{S}_\ell(N)$ by a slightly larger computed value (Section~\ref{sec:Algorithm}); the energy layer bounds of Corollary~\ref{cor:energy_layer_ub} are used only to discard candidate layers.
\end{proposition}

If the double sum in~\eqref{eq:sp_cut} vanishes, then $(s_1, \dots, s_\ell)$ is a base-rooted unique-support chain and~\eqref{eq:spinal_full_bound} gives the cut; otherwise the big-$M$ term covers every $X_{s_\ell} \le \overline{x}_{s_\ell}$. The details are given in Appendix~\ref{app:sp_cuts_proof}.

\paragraph{Alignment with Paterson and Zwick's Definition of Spinal Stacks.}
Following Paterson and Zwick~\cite[Definitions~2.4 and~3.1]{paterson2009overhang}, the \emph{principal block} of a feasible stack is the block attaining its maximum overhang $z(X) = \max_{i \in \mathcal{N}}(X_i + W)$; if multiple blocks tie for the maximum overhang, the \emph{lowest} tied block is chosen. Because distinct blocks within any single layer $l$ are separated by at least unit distance ($|X_i - X_k| \ge 1 > 0$), each layer has a unique rightmost block $s_l = i_{l,0}$, so ties for the maximum overhang can occur only across distinct layers and the principal block is uniquely determined as $s_{\ell^*}$, where $\ell^* := \min\{l \in \{1, \dots, H\} : X_{s_l} + W = z(X)\}$. We say that block~$i$ \emph{rests on} block~$j$ if $L_i = L_j + 1$ and their effective cores intersect, $|X_i - X_j| \le 2w$, which is the geometric condition of Section~\ref{sec:Formulation} under which $j$ can support~$i$. The \emph{support set} of the stack is defined recursively as the set containing the principal block $s_{\ell^*}$ together with every block on which a block of the support set rests, and the stack is called \emph{spinal} if its support set contains exactly one block in each layer $1, \dots, \ell^*$, and \emph{non-spinal} otherwise. These notions depend only on the positions $X$, not on the contact representation~$T$.

For $\eta = 0$, Paterson and Zwick~\cite[Section~2]{paterson2009overhang} let a block rest on another when $|X_i - X_j| < 1$; our relation, with $2w = 1$, adds the point contacts $|X_i - X_j| = 1$. Every stack balanced in their sense is feasible in our model, and if it is spinal under our relation, it is spinal in their sense (Appendix~\ref{app:pz_relation}). Hence excluding non-spinal stacks above a threshold under our relation excludes them under their definition as well.

\begin{figure}[t]
\centering
\includegraphics[width=\textwidth]{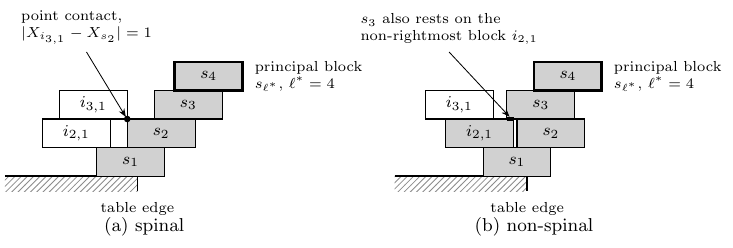}
\caption{Schematic layered configurations, not drawn as balanced stacks, illustrating the definitions of Section~\ref{subsec:spinal_bound}. Shaded blocks form the support set of the principal block (thick outline); $s_t$ is the rightmost block of layer~$t$, and $i_{t,q}$ with $q \ge 1$ are the other blocks of that layer. (a)~Each $s_t$ rests only on $s_{t-1}$, so the stack is spinal. Block $i_{3,1}$ touches $s_2$ at a single point; our relation counts this as resting and that of Paterson and Zwick does not, which does not matter here because $i_{3,1}$ is not in the support set. (b)~Block $s_3$ also rests on the non-rightmost block $i_{2,1}$, so the support set has two blocks in layer~$2$ and the stack is non-spinal; in its maximal contact representation, the covering constraint in~\eqref{eq:ns_branch_constraints} holds through $T_{s_3, i_{2,1}} = 1$, and the lowest impure layer of Proposition~\ref{prop:t0_split} is $t_0 = 3$.}
\label{fig:spinal_schematic}
\end{figure}

Figure~\ref{fig:spinal_schematic} illustrates these notions. We relate spinality to the contact variables in Lemma~\ref{lem:maximal_contacts} below.

For each energy-admissible profile $\mathbf{n} \in \mathcal{P}_2^{\mathrm{E}}(N)$ and cutoff $\gamma \ge \underline{z}$, let $\mathcal{F}_{\mathrm{SP}}(\mathbf{n})$ denote the feasible set $\mathcal{F}(\mathbf{n})$ augmented with the spinal-bound cut~\eqref{eq:sp_cut} for every candidate layer $\ell$ with $s_\ell \in \mathcal{C}(\mathbf{n}, \gamma)$, and let $\mathcal{S}_{\mathrm{SP}}(\mathbf{n}, \gamma)$ denote the corresponding cutoff subproblem~\eqref{eq:profile_subproblem} over $\mathcal{F}_{\mathrm{SP}}(\mathbf{n})$.

\begin{theorem}[Spinal Certificate Decomposition]
\label{thm:spinal_cert_decomp}
Let $N \ge 3$ and $\eta \in [0, 1)$, and suppose that condition~\eqref{eq:s2p_condition} of Corollary~\ref{cor:s2p_elimination} holds. Let $\underline{z} \le z^*(N, \eta)$ be the overhang of a known feasible stack, and fix any cutoff $\gamma \ge \underline{z}$.
\begin{enumerate}
\item[(i)] \textbf{Spinal versus non-spinal decomposition:} For each $\mathbf{n} \in \mathcal{P}_2^{\mathrm{E}}(N)$, let $\overline{u}_{\mathrm{ns}}(\mathbf{n})$ be any valid upper bound on the maximum overhang achievable by non-spinal stacks in $\mathcal{F}(\mathbf{n})$; in particular, if $\mathcal{S}_{\mathrm{SP}}(\mathbf{n}, \gamma)$ returns a closed or dual bound $\overline{u}_{\mathrm{SP}}(\mathbf{n})$ according to cases~(a) to~(c) of Theorem~\ref{thm:profile_decomp}, then $\overline{u}_{\mathrm{ns}}(\mathbf{n}) := \overline{u}_{\mathrm{SP}}(\mathbf{n})$ is valid. Then
\begin{equation}
\label{eq:spinal_decomp_bound}
z^*(N, \eta) \le \max\!\left\{
2w \max_{1 \le \ell \le N} \hat{S}_\ell(N),\;
\max_{\mathbf{n} \in \mathcal{P}_2^{\mathrm{E}}(N)} \overline{u}_{\mathrm{ns}}(\mathbf{n})
\right\}.
\end{equation}
\item[(ii)] \textbf{Exact cutoff certification (including exception sizes):} Independently of whether $\gamma$ exceeds $2w \max_{1 \le \ell \le N} \hat{S}_\ell(N)$, if every profile $\mathbf{n} \in \mathcal{P}_2^{\mathrm{E}}(N)$ satisfies $\overline{u}_{\mathrm{SP}}(\mathbf{n}) \le \gamma$, then $z^*(N, \eta) \le \gamma$. In particular, when $\eta = 0$ and $\gamma = S^*(N) + \delta$ with $\delta > 0$, closing every subproblem $\mathcal{S}_{\mathrm{SP}}(\mathbf{n}, \gamma)$ certifies $z^*(N, 0) \le S^*(N) + \delta$, and, whenever $\hat{S}(N) = S^*(N)$, every feasible stack exceeding $\gamma$ is forced to satisfy $\Sigma_\ell \ge 1$ in~\eqref{eq:sp_cut}.
\end{enumerate}
\end{theorem}

The proof is given in Appendix~\ref{app:spinal_cert_proof}.

\begin{remark}[Exception Sizes $N \in \{5, 7\}$ and Big-$M$ Conditioning]
\label{rem:sp_exceptions_numerics}
\begin{enumerate}
\item[(i)] \emph{Exception sizes.} For $\eta = 0$ and $N \in \{5, 7\}$, the exception sizes of Conjecture~3.9 of Paterson and Zwick~\cite{paterson2009overhang}, the loaded-spine maximum $S^*(5) = S_3^*(5) \approx 1.31716$ or $S^*(7) = S_4^*(7) \approx 1.53187$ needs a spine of height $k$ with $2k > N$, and no spinal stack attains it (Proposition~\ref{prop:exception_sizes}). The first term of~\eqref{eq:spinal_decomp_bound} is then not tight, but part~(ii) of Theorem~\ref{thm:spinal_cert_decomp} applies with the cutoff $\gamma = \underline{z}(N) + \delta$, and the solver closes every profile, including those with spines of the unrealizable height.
\item[(ii)] \emph{Big-$M$ magnitude and integrality tolerance.} Because of the layer-rank boxes of Proposition~\ref{prop:profile_tightening}(i), the constants~\eqref{eq:sp_bigm} stay small: at $\eta = 0$, their maximum over all profiles and candidate layers grows from $1.24$ at $N = 5$ to $6.10$ at $N = 16$. With integer feasibility tolerance $10^{-9}$, each indicator in~\eqref{eq:sp_cut} relaxes the cut by at most $M_\ell \cdot 10^{-9}$, so the total leakage stays below $10^{-7}$ for $N \le 16$, two orders of magnitude below the certification margin $\delta = 10^{-5}$.
\end{enumerate}
\end{remark}

\begin{lemma}[Maximal Contact Representation]
\label{lem:maximal_contacts}
Let $\eta \in [0, 1)$ and let $(X, T, F, P, \tau)$ be a feasible point of the profile subproblem $\mathcal{S}_{\mathrm{SP}}(\mathbf{n}, \gamma)$ or of a branch $\mathcal{S}_{\mathrm{ns}}(\mathbf{n}, \ell^*, \gamma_{\mathrm{ns}})$ with the non-spinal support constraint of~\eqref{eq:ns_branch_constraints} omitted. Set $\hat{T}_{i,j} := 1$ whenever block~$i$ rests on block~$j$, and $\hat{T}_{i,j} := T_{i,j}$ otherwise; for every pair with $\hat{T}_{i,j} = 1 > T_{i,j}$ set $\hat{F}_{i,j} := 0$, $\hat{\tau}_{i,j} := 0$, and $\hat{P}_{i,j} := (X_i + X_j)/2$, and keep all other variables. Then $(X, \hat{T}, \hat{F}, \hat{P}, \hat{\tau})$ is feasible for the same problem, and in $\hat{T}$ every block is linked to every block on which it rests.
\end{lemma}

The proof checks the constraint families one at a time: the new contacts carry zero force and zero moment, the contact rows hold at $\hat{P}_{i,j}$, and every remaining row that involves $T$ either holds for every contact representation or is nondecreasing in~$T$. It is given in Appendix~\ref{app:maximal_contacts_proof}.

We call $\hat{T}$ the \emph{maximal contact representation} of~$X$. Let $\Sigma_{\ell^*} := \sum_{t \in J_{\ell^*}} \sum_{q = 1}^{n_{t-1} - 1} T_{s_t,\, i_{t-1, q}}$ denote the rightmost-chain non-spinal indicator sum. If $\Sigma_{\ell^*} = 0$ in the maximal contact representation, the stack is spinal: every block $s_t$ with $t \ge 2$ is linked to, and hence rests on, at least one block of layer $t - 1$, and it rests on no non-rightmost block of that layer, trivially when $n_{t-1} = 1$ and in particular for $t = 2$; so $s_t$ rests solely on $s_{t-1}$, and the support set of $s_{\ell^*}$ is $\{s_1, \dots, s_{\ell^*}\}$. Hence every non-spinal stack satisfies $\Sigma_{\ell^*} \ge 1$ at its principal layer $\ell^*$ in its maximal contact representation. The converse fails: a spinal stack whose spine passes through a non-rightmost block also has $\Sigma_{\ell^*} \ge 1$, which only makes the non-spinal branches below conservative. Finally, the principal block of a non-spinal stack is not block~$1$, whose support set $\{1\}$ is spinal, and $s_2$ rests solely on block~$1$ because $n_1 = 1$; if $s_3$ rests on a non-rightmost block of layer~$2$, then $X_{s_3} \le X_3 + 2w \le X_2 - 1 + 2w = X_{s_2} - \eta \le X_{s_2}$, contradicting the minimality of $\ell^*$, so every non-spinal stack has principal layer $\ell^* \ge 4$ and satisfies $J_{\ell^*} \neq \emptyset$.

\begin{corollary}[Certified Non-Spinal Exclusion and Structural Spinal Optimality]
\label{cor:nonspinal_exclusion}
Let $N \ge 4$ and $\eta \in [0, 1)$, with $w = (1 - \eta)/2$, and let $z_{\mathrm{ns}}^*(N, \eta)$ denote the supremum overhang over all feasible non-spinal single-root $N$-block $\eta$-stable stacks (with $z_{\mathrm{ns}}^*(N, \eta) := -\infty$ if none exists). Fix a certification threshold $\gamma_{\mathrm{ns}} > 0$. Suppose that a previously certified upper bound $\overline{z}_{\mathrm{ns}}(N - 1, \eta)$ on non-spinal $(N - 1)$-block stacks, or more conservatively a general upper bound $\overline{z}(N - 1, \eta) \ge z^*(N - 1, \eta)$, satisfies
\begin{equation}
\label{eq:s2p_nonspinal_cond}
\overline{z}_{\mathrm{ns}}(N - 1, \eta) + \frac{w}{N} \le \gamma_{\mathrm{ns}}.
\end{equation}
For each energy-admissible profile $\mathbf{n} = (n_1, \dots, n_H) \in \mathcal{P}_2^{\mathrm{E}}(N)$ and each candidate principal layer $\ell^* \in \{2, \dots, H\}$, define the non-spinal branch subproblem $\mathcal{S}_{\mathrm{ns}}(\mathbf{n}, \ell^*, \gamma_{\mathrm{ns}})$ as the problem of maximizing $X_{s_{\ell^*}} + W$ with objective cutoff $\gamma_{\mathrm{ns}}$ over the profile constraint set $\mathcal{F}(\mathbf{n})$, whose optimality-based constraints Lemma~\ref{lem:wlog_base_shift} justifies for arbitrary feasible stacks, augmented with the spinal cuts~\eqref{eq:sp_cut}, the downward reachability bounds $X_{s_l} \ge \gamma_{\mathrm{ns}} - W - 2w(\ell^* - l)$ ($1 \le l < \ell^*$), and the rightmost-layer dominance and non-spinal support constraints
\begin{equation}
\label{eq:ns_branch_constraints}
X_{s_l} \le X_{s_{\ell^*}} \quad (\forall l \in \{1, \dots, H\}),
\qquad
\sum_{t \in J_{\ell^*}} \sum_{q = 1}^{n_{t-1} - 1} T_{s_t,\, i_{t-1, q}} \ge 1,
\end{equation}
where $J_{\ell^*} = \{t \in \{3, \dots, \ell^*\} : n_{t-1} \ge 2\}$ and any branch with $J_{\ell^*} = \emptyset$ is immediately closed as infeasible. The second constraint in~\eqref{eq:ns_branch_constraints} is a covering condition: every non-spinal stack above the threshold satisfies it in its maximal contact representation, while some spinal stacks may satisfy it as well, which only enlarges the branches. A branch \emph{closes} if it is proved that none of its feasible points has objective strictly greater than $\gamma_{\mathrm{ns}}$; since the downward reachability bounds hold only for stacks whose overhang is at least $\gamma_{\mathrm{ns}}$, a report of infeasibility for a branch is such a proof as well. If every branch $\mathcal{S}_{\mathrm{ns}}(\mathbf{n}, \ell^*, \gamma_{\mathrm{ns}})$ closes, then:
\begin{enumerate}
\item[(i)] $z_{\mathrm{ns}}^*(N, \eta) \le \gamma_{\mathrm{ns}}$, so every feasible single-root $N$-block $\eta$-stable stack with overhang strictly exceeding $\gamma_{\mathrm{ns}}$ is spinal;
\item[(ii)] in particular, when $\eta = 0$ and an exact feasible spinal stack attains $\underline{z}_{\mathbb{Q}}(N) > \gamma_{\mathrm{ns}}$, for instance with $\gamma_{\mathrm{ns}} = \lceil (1 - \rho) S^*(N) \rceil_7$ and $\rho \in (0, 1)$, every optimal single-root $N$-block stack (and every stack strictly exceeding $\gamma_{\mathrm{ns}}$) is spinal.
\end{enumerate}
\end{corollary}

The proof, given in Appendix~\ref{app:nonspinal_exclusion_proof}, removes block~$1$ when $n_2 = 1$ and uses condition~\eqref{eq:s2p_nonspinal_cond}; when $n_2 = 2$, it normalizes a non-spinal stack above $\gamma_{\mathrm{ns}}$ by Lemma~\ref{lem:wlog_base_shift}, passes to its maximal contact representation, and places it in the branch of its principal layer.

Each branch can be split further by the lowest layer at which the rightmost chain becomes impure. We use this refinement for the three branches at $N = 16$ that do not close within the per-branch time limit (Section~\ref{subsec:unconditional_certs}).

\begin{proposition}[Lowest-Impure-Layer Refinement of Non-Spinal Branches]
\label{prop:t0_split}
Under the hypotheses of Corollary~\ref{cor:nonspinal_exclusion}, fix a branch $(\mathbf{n}, \ell^*)$ with $J_{\ell^*} \neq \emptyset$. For each $t_0 \in J_{\ell^*}$, let the sub-branch $\mathcal{S}_{\mathrm{ns}}(\mathbf{n}, \ell^*, t_0, \gamma_{\mathrm{ns}})$ be $\mathcal{S}_{\mathrm{ns}}(\mathbf{n}, \ell^*, \gamma_{\mathrm{ns}})$ augmented with
\begin{gather}
T_{s_t,\, i_{t-1, q}} = 0 \quad (t \in J_{\ell^*},\; t < t_0,\; 1 \le q \le n_{t-1} - 1),
\qquad
\sum_{q = 1}^{n_{t_0 - 1} - 1} T_{s_{t_0},\, i_{t_0 - 1, q}} \ge 1, \label{eq:t0_partition} \\
X_{s_{t_0}} \le X_{s_{t_0 - 1}} - 1 + 2w,
\qquad
X_{s_{t_0 - 1}} + W \le 2w \, \hat{S}_{t_0 - 1}(N). \label{eq:t0_cuts}
\end{gather}
If every sub-branch with $t_0 \in J_{\ell^*} \setminus \{\ell^*\}$ closes in the sense of Corollary~\ref{cor:nonspinal_exclusion}, then no feasible non-spinal single-root $N$-block $\eta$-stable stack has profile $\mathbf{n}$, principal layer $\ell^*$, and overhang exceeding $\gamma_{\mathrm{ns}}$. In particular, the branch $(\mathbf{n}, \ell^*)$ may be replaced by these sub-branches in Corollary~\ref{cor:nonspinal_exclusion}.
\end{proposition}

The proof, which repeats the argument of Corollary~\ref{cor:nonspinal_exclusion} and places the normalized stack in the sub-branch of its lowest impure layer, is given in Appendix~\ref{app:t0_split_proof}.

The sub-branch $t_0 = \ell^*$ therefore needs no computation. Because the dominance inequalities in~\eqref{eq:ns_branch_constraints} are weak, this sub-branch may contain stacks in which layer $\ell^* - 1$ ties with layer $\ell^*$. Such a stack has its principal layer strictly below $\ell^*$, so if it is non-spinal, it is covered by the branch of that layer.

\begin{remark}[Validity Classes of Subproblem Constraints]
\label{rem:constraint_classification}
Every constraint of the subproblems $\mathcal{S}_{\mathrm{SP}}(\mathbf{n}, \gamma)$ and the branches $\mathcal{S}_{\mathrm{ns}}(\mathbf{n}, \ell^*, \gamma_{\mathrm{ns}})$ for $\mathbf{n} \in \mathcal{P}_2(N)$ belongs to one of three validity classes, listed in Table~\ref{tab:constraint_taxonomy}. Class~I constraints hold for every feasible single-root stack with profile $\mathbf{n}$. Class~II constraints hold without loss of generality for every feasible single-root stack with $n_2 = 2$, because Lemma~\ref{lem:wlog_base_shift} maps such a stack to one that satisfies them, has the same layer profile and non-spinal status, and in which no block outside the base loses overhang. Two families need a word of qualification. The cuts \textsf{C2-Layer} are valid for every stack satisfying the coordinate bounds that enter the width $\overline{\Delta}_\ell$, and the big-$M$ spinal cuts for every stack satisfying the bounds from which the constants $M_\ell$ are computed (Propositions~\ref{prop:c2_layer} and~\ref{prop:sp_cuts}); in our implementation both sets of bounds are of Class~II. Class~III constraints define the non-spinal branches and are not valid for every feasible stack: they hold for every stack that a branch must exclude, namely a non-spinal stack with principal layer $\ell^*$ and overhang above $\gamma_{\mathrm{ns}}$, taken in its maximal contact representation and, for a sub-branch, with lowest impure layer $t_0$. For example, the downward reachability inequalities hold in the relative form $X_{s_{\ell^*}} - X_{s_l} \le 2w(\ell^* - l)$ for every feasible stack, but in the cutoff form only for stacks whose overhang exceeds $\gamma_{\mathrm{ns}}$. Consequently, Classes~I and~II are valid for certifying optimal overhangs (Theorem~\ref{thm:spinal_cert_decomp}), and all three classes are valid for the non-spinal exclusion (Corollary~\ref{cor:nonspinal_exclusion} and Proposition~\ref{prop:t0_split}).
\end{remark}

\begin{plaintable}[t]
\caption{Validity taxonomy of all constraint families used in the layer-profile subproblems $\mathcal{S}_{\mathrm{SP}}(\mathbf{n}, \gamma)$ and non-spinal exclusion branches $\mathcal{S}_{\mathrm{ns}}(\mathbf{n}, \ell^*, \gamma_{\mathrm{ns}})$ for $\mathbf{n} \in \mathcal{P}_2(N)$.}
\label{tab:constraint_taxonomy}
\resizebox{\textwidth}{!}{%
\begin{tabular}{@{}lll@{}}
\toprule
Validity class & Constraint family & Theoretical basis \\
\midrule
Class~I & Equilibrium, contact intervals, non-overlap, base linearization & Section~\ref{sec:Formulation} \\
(all feasible & Support degrees, base and coordinate envelopes & Lemma~\ref{lem:support_bound} \\
single-root & Layer-width doubling ($n_t \le 2 n_{t-1}$, $n_t \le 2t - 2$) & Theorem~\ref{thm:layer_doubling} \\
stacks) & Singleton-layer fixing, non-crossing, relative rank & Proposition~\ref{prop:profile_tightening} \\
 & Cut-set cuts \textsf{C2} & Proposition~\ref{prop:c2_singleton} \\
 & Cut-set cuts \textsf{C2-Layer} (valid under the Class~II box bounds in $\overline{\Delta}_\ell$) & Proposition~\ref{prop:c2_layer}, Lemma~\ref{lem:wlog_base_shift} \\
 & Profile-level energy infeasibility criterion & Corollary~\ref{cor:energy_infeasibility} \\
 & Big-$M$ spinal cuts (family~\textsf{SP}; valid under the Class~II box bounds that define $M_\ell$) & Proposition~\ref{prop:sp_cuts}, Lemma~\ref{lem:wlog_base_shift} \\
\midrule
Class~II & Table-edge equality \textsf{S1} ($P_{1,0} = -\epsilon$, $\sum X_i = -N\epsilon$) & Lemma~\ref{lem:wlog_base_shift} \\
(WLOG for & Base right-corner equality \textsf{S4} ($P_{2,1} = X_1 + w$) & Theorem~\ref{thm:base_right_corner}, Lemma~\ref{lem:wlog_base_shift} \\
feasible stacks & Tightened base bound $X_1 \le \frac{N-3}{N}w - \epsilon$ & Corollary~\ref{cor:base_ub_tight}, Lemma~\ref{lem:wlog_base_shift} \\
with $n_2 = 2$) & Absolute coordinate and contact boxes from $X_1 \ge -W$ & Proposition~\ref{prop:profile_tightening}, Lemma~\ref{lem:wlog_base_shift} \\
 & Layer-wise energy upper bounds (candidate-layer filter) & Corollary~\ref{cor:energy_layer_ub}, Lemma~\ref{lem:wlog_base_shift} \\
\midrule
Class~III & Downward reachability (cutoff form) and dominance $X_{s_l} \le X_{s_{\ell^*}}$ & Corollary~\ref{cor:nonspinal_exclusion} \\
(stacks a branch & Non-spinal covering constraint in~\eqref{eq:ns_branch_constraints} (may also admit spinal stacks) & Corollary~\ref{cor:nonspinal_exclusion}, Lemma~\ref{lem:maximal_contacts} \\
must exclude) & Lowest-impure-layer constraints~\eqref{eq:t0_partition}, \eqref{eq:t0_cuts} & Proposition~\ref{prop:t0_split} \\
\bottomrule
\end{tabular}}
\end{plaintable}

\section{Algorithm and Exact Rational Construction}
\label{sec:Algorithm}

This section assembles the results of Sections~\ref{sec:Structural} and~\ref{sec:SpinalCert} into a certification algorithm and constructs the exact rational stacks that supply the lower bounds.

\subsection{Two-Stage Certification and Non-Spinal Exclusion Pipeline}
\label{subsec:TwoStagePipeline}

Given $N \ge 5$, a margin $\eta \in [0, 1)$, an exact lower bound $\underline{z}_{\mathbb{Q}}(N, \eta)$, and a cutoff $\gamma$, the procedure runs in two stages, both for optimality certification (Theorem~\ref{thm:spinal_cert_decomp}, with $\gamma = \gamma_N$) and for non-spinal exclusion (Corollary~\ref{cor:nonspinal_exclusion}, with $\gamma = \gamma_{\mathrm{ns}}$).

\paragraph{Stage~1: Base-Singleton Check, Profile Enumeration, and Energy Screening.}
The algorithm first checks the base-singleton condition with the certified bounds of size $N - 1$: condition~\eqref{eq:s2p_condition}, $\overline{z}(N - 1, \eta) + w / N < \underline{z}_{\mathbb{Q}}(N, \eta)$, for optimality certification, and condition~\eqref{eq:s2p_nonspinal_cond} for non-spinal exclusion. These checks use only bounds certified for $N - 1$ and never enter a model, so the dependency chain over $N$ is acyclic and starts from the anchor at $N = 4$. Only if the check holds does the algorithm enumerate $\mathcal{P}_2(N)$; it then discards every profile with $\Delta E(\mathbf{n}) < 0$, evaluated in closed form, before any model is built. Each remaining profile yields the subproblem $\mathcal{S}_{\mathrm{SP}}(\mathbf{n}, \gamma)$, which is closed analytically when its candidate set $\mathcal{C}(\mathbf{n}, \gamma)$ is empty, and, for non-spinal exclusion, one branch $\mathcal{S}_{\mathrm{ns}}(\mathbf{n}, \ell^*, \gamma_{\mathrm{ns}})$ for every $\ell^* \in \{2, \dots, H\}$, which is closed analytically when $J_{\ell^*} = \emptyset$, because the chain $s_1, \dots, s_{\ell^*}$ is then pure.

\paragraph{Stage~2: Subproblem Assembly.}
Fixing the profile removes the integer variables $L_i$ and $b_i$ of $\mathcal{M}_0(N, \eta)$ and restricts the contact quadruples $(T_{i,j}, F_{i,j}, P_{i,j}, \tau_{i,j})$ to pairs in adjacent layers. Each subproblem consists of the contact-interval rows and the bilinear equilibrium equations of $\mathcal{M}_0$ on these pairs together with the reductions of Sections~\ref{sec:Structural} and~\ref{sec:SpinalCert}, whose validity classes are listed in Table~\ref{tab:constraint_taxonomy}. For optimality certification, the objective is $z = \max_{b \in \mathcal{C}(\mathbf{n}, \gamma)} (X_b + W)$; for a non-spinal branch, it is $z = X_{s_{\ell^*}} + W$, and the branch adds the constraints of Corollary~\ref{cor:nonspinal_exclusion}. In both cases, $\gamma$ is passed to the solver as an objective cutoff.

\paragraph{Parallel Execution and Solver Settings.}
The subproblems and branches are independent and are dispatched to a pool of single-thread worker processes. Each worker calls a spatial branch-and-bound solver for nonconvex quadratic programs with four uniform settings: (i) nonconvex bilinear handling enabled; (ii) primal heuristics disabled, by design, since the cutoff comes from an exact lower bound and the runs only have to prove that no solution exceeds it (the setting can still change the course of the search); (iii) feasibility and optimality tolerances $10^{-6}$; and (iv) integer feasibility tolerance $10^{-9}$, so that the leakage of the indicators in the spinal cuts~\eqref{eq:sp_cut} stays below $10^{-7}$ (Remark~\ref{rem:sp_exceptions_numerics}(ii)). A subproblem counts as closed only with the status cutoff or infeasible, or by an analytic closure, so the relative gap parameter never decides a closure. When presolve reports the ambiguous status infeasible or unbounded, the worker disables dual reductions and solves again. Effort is reported in thread-seconds: the run time of each subproblem or branch multiplied by its threads, summed over all of them, with a run stopped by its time limit counted at its full time. The time limit per subproblem is $3{,}600$ to $259{,}200$~s for the certificates, $7{,}200$~s for the non-spinal branches, and $86{,}400$~s for the sub-branches at $N = 16$; Table~\ref{tab:run_settings} in Appendix~\ref{app:run_settings} lists the threads, parallel workers, time limits, seeds, and remaining parameters of every reported run.

\paragraph{Automated Certificate Auditing.}
After all workers finish, an auditor checks four invariants on the stored records: (a) every profile of $\mathcal{P}_2(N)$, or every branch $(\mathbf{n}, \ell^*)$ together with every sub-branch of a $t_0$ split, has a record, and no record lies outside this set; any repeated record is reported and the last one is used; (b) every record is closed analytically, by cutoff, or by proved infeasibility rather than by a time limit or numerical error; (c) every recorded bound and the cutoff of each run, or the solver threshold for non-spinal exclusion, do not exceed the certified value; and (d) the conditions that link consecutive sizes hold with the audited results of size $N - 1$, namely~\eqref{eq:s2p_condition}, starting from the anchor $\overline{z}(4)$, and~\eqref{eq:s2p_nonspinal_cond}. Analytic closures are checked against their conditions, recomputed by the auditor: the energy test, an empty candidate set, $J_{\ell^*} = \emptyset$, or $t_0 = \ell^*$. Repeated records with conflicting statuses are rejected, and where a record also stores the raw solver status, the objective bound, and the cutoff, these must agree with the closed status. The auditors thus verify coverage, statuses, and thresholds, while the upper bounds themselves rest on the closure decisions of the solver; they run without an optimization solver and are included in the companion repository. Positive reverse controls at lowered cutoffs (Section~\ref{sec:Computational}) recover known feasible spinal and non-spinal stacks, an empirical check that neither family~\textsf{SP} nor~\eqref{eq:ns_branch_constraints} excludes them. The primal lower bounds of Section~\ref{subsec:RationalConstructor} are verified with zero residual in exact rational arithmetic, whereas the dual upper bounds and the non-spinal exclusions rely on double-precision spatial branch-and-bound and therefore hold within the stated solver tolerances.

\paragraph{Analytic Constants and Rounding.}
The analytic quantities that discard profiles, candidate layers, or branches are computed in double precision during the runs and checked afterwards, so that rounding cannot discard a feasible stack. A separate program repeats every such decision of the reported runs in exact rational arithmetic, with the cutoff equal to the exact value of the double passed to the solver, and archives one record per decision with its criterion, threshold, and margin. At $\eta = 0$, $\Delta E(\mathbf{n})$ is a multiple of $1/4$, and every profile discarded by Corollary~\ref{cor:energy_infeasibility} has $\Delta E(\mathbf{n}) \le -1/4$. A candidate layer is dropped when its base envelope, its energy bound~\eqref{eq:energy_layer_ub}, or, for a pure chain, its loaded-spine bound does not exceed the cutoff. The program decides the energy bound without square roots: in the first case of~\eqref{eq:energy_layer_ub}, for instance, the bound does not exceed $\gamma$ if and only if $c := \gamma - w - (m - 1)/2 \ge 0$ and $c^2 \ge \Delta E(\mathbf{n}) / A_m$. It confirms all $3{,}119$ candidate layers dropped in the certificates, with smallest margin $0.0165$. The loaded-spine constants $2w\,\hat{S}_\ell(N)$ in~\eqref{eq:sp_cut} and~\eqref{eq:t0_cuts} are computed by bisection in double precision and increased by $10^{-9}$; they differ from a $60$-digit evaluation by at most $2 \times 10^{-15}$ for all $\ell \le N \le 20$, and interval enclosures with outward rounding, computed as in Appendix~\ref{app:ExactSpinal}, confirm that the constants used exceed $2w\,\hat{S}_\ell(N)$. At every $N \le 16$, the optimal loads of the maximizing height have increments of at least $1.19$, so $\hat{S}(N) = S^*(N)$ by Proposition~\ref{prop:spinal_bound}(i). For non-spinal exclusion, the solver receives the double-precision value of $0.99\, S^*(N)$, which does not exceed the reported threshold $\gamma_{\mathrm{ns}}(N)$. No branch is closed by a conflict of its bounds alone, as the program confirms, so every branch reported infeasible is closed by the solver in presolve. The remaining bounds of the subproblems, such as the coordinate and contact boxes, are evaluated in double precision from the closed forms of Sections~\ref{sec:Structural} and~\ref{sec:SpinalCert}; their rounding errors lie far below the feasibility tolerance $10^{-6}$, so, like the solver itself, they are valid within solver tolerances. We thus separate two kinds of numerical error: the analytic decisions that discard profiles, candidate layers, or branches are checked in exact or interval arithmetic, whereas the coefficients of the floating-point model and its solution are valid within solver tolerances only. For the hardest models, the $11$ sub-branches at $N = 16$ in Table~\ref{tab:n16_split}, the nonzero linear constraint coefficients lie in $[0.973, 240]$, all quadratic coefficients equal~$1$, the nonzero variable bounds lie in $[0.377, 35.2]$, and the nonzero right-hand sides lie in $[0.5, 59.5]$, ranges that the solver log rounds to one significant figure; a rerun of these sub-branches, and of all subproblems and branches at $N \in \{5, 7\}$, reproduces the archived node counts, and its solver logs, archived together with the models, report no numerical difficulties. The archived records of all runs store the closure status, the bound, the nodes, and the time of every subproblem and branch, together with the cutoff and the parameters of each run; the reruns at $N \in \{5, 7\}$ and of the sub-branches at $N = 16$ also store the raw solver status, the objective bound, and the cutoff of every record, with the solver log and the model of every subproblem solved.

\subsection{Closed-Form Exact Rational Stack Construction and Verification}
\label{subsec:RationalConstructor}

Relative gaps of order $5 \times 10^{-6}$ need lower bounds that do not depend on feasibility tolerances: at a feasibility tolerance of $10^{-6}$, floating-point incumbents at $N = 12$ and $N = 13$ exceeded $S^*(N)$ by $3.8 \times 10^{-6}$ and $4.4 \times 10^{-6}$ through small constraint violations. We therefore construct explicit rational stacks at $\eta = 0$ ($W = 1/2$) for every $N \in \{4, \dots, 16\}$ and verify every equilibrium and geometric constraint with zero residual in exact rational arithmetic.

\paragraph{Rational Loads of the Spine.}
Following Paterson and Zwick~\cite[Section~3]{paterson2009overhang}, index a $k$-block spine from top to bottom as $B_1, \dots, B_k$, so that $B_k$ is the base block resting on the table edge at $P_{B_k, 0} = 0$ and $B_r$ lies in layer $k - r + 1$. Let $t_r$ be the load transmitted from $B_r$ to $B_{r+1}$ ($t_0 := 0$, $t_k := N$), and let $\omega_r := t_r - t_{r-1} - 1$ be the point weight applied at the left edge $x_r$ of $B_r$ (their $w_r$, renamed to avoid a clash with the half-width $w$). By Proposition~\ref{prop:spinal_bound}(i), bisection on~\eqref{eq:pz_recurrence} gives the maximizer $(t_1^*, \dots, t_{k-1}^*)$ of~\eqref{eq:pz_unconstrained}; for every height used in Table~\ref{tab:rational_seeds} except at $N = 5$, its increments are at least $1.23$, so it attains $S_k^*(N) = \hat{S}_k(N)$. Since the objective is stationary there, rounding each $t_r^*$ to a rational $t_r$ with denominator at most $2{,}000$ changes the overhang only to second order, by $10^{-10}$ to $10^{-8}$, and keeps $\omega_r > 0$. Given such loads with $\omega_1 < 1$, the spine follows in $\mathbb{Q}$:
\begin{equation}
\label{eq:rat_spine_steps}
d_r = \frac{\omega_r + 1/2}{t_r} \quad (r = 1, \dots, k),
\qquad
x_k = d_k - 1,
\qquad
x_r = x_{r+1} + d_r \quad (r = k - 1, \dots, 1),
\end{equation}
with spine block centers $X(B_r) = x_r + 1/2$, spine-to-spine contacts $(F_{B_r, B_{r+1}}, \tau_{B_r, B_{r+1}}) = (t_r,\, t_r(x_{r+1} + 1))$ for $r \in \{1, \dots, k - 1\}$, base table contact $(F_{B_k, 0}, \tau_{B_k, 0}) = (N, 0)$, and exact rational stack overhang
\begin{equation}
\label{eq:rat_overhang}
\underline{z}_{\mathbb{Q}}(N) = x_1 + 1 = \sum_{r=1}^{k} d_r \in \mathbb{Q}.
\end{equation}

\paragraph{Shielded Spines and the Exception $N = 5$.}
For $N \in \{4, 6, 7, \dots, 16\}$, the point weights are realized by the shielded architecture of~\cite[Section~3]{paterson2009overhang}: $k - 1$ shield blocks $B_1', \dots, B_{k-1}'$, where $B_r'$ shares layer $k - r + 1$ with $B_r$, a bridge block $B_0'$ in layer $k + 1$, and a vertical tower of $v_1 := N - 2k \ge 0$ blocks on $B_1'$. Rational left edges $y_1, \dots, y_{k-1}$ of the shields and, when $v_1 \ge 1$, a rational contact seed $z_1$ of the bridge determine every remaining coordinate, force, and moment in closed form over $\mathbb{Q}$ by force and moment balance (Appendix~\ref{subsec:ReconIndex}). Every such size admits a strictly feasible seed, with overhang within $1.5 \times 10^{-9}$ to $3.1 \times 10^{-8}$ of $S^*(N)$ for $N \ne 7$ and within $7.1 \times 10^{-10}$ of $S_3^*(7)$ at $N = 7$, where $k = 3$ and $v_1 = 1$. At $N = 5$, a shielded spine of height $3$ would need $2k = 6 > 5$ blocks. We use instead a bridgeless spine with profile $(1, 2, 2)$ whose top block $B_1$ carries no external weight ($\omega_1 = 0$, $t_1 = 1$), while two counterweights, $C_2$ in layer~$2$ and $C_1$ in layer~$3$, supply $\omega_2 = t_2 - 2$ and $\omega_3 = 4 - t_2$. Maximizing over $t_2$ gives $t_2^* = \sqrt{15/2}$ and the overhang $\frac{12}{5} - 2\sqrt{3/10} \approx 1.3045548850$, which Proposition~\ref{prop:exception_sizes} identifies as the maximum overhang $S(5)$ of every $5$-block spinal stack; the numerical certificate brackets $z^*(5) \in [1.3045548840, 1.3045649]$ (Theorem~\ref{thm:main_certified_optimality}). Rational seeds $t_2$ and $X(C_2)$ determine the rest in closed form (Appendix~\ref{subsec:ReconIndex}).

\paragraph{Zero-Residual Verification.}
Each rational stack is passed to an independent verifier that uses exact rational arithmetic only and asserts with zero tolerance that (i) the layer profile is valid and single-rooted; (ii) blocks of the same layer satisfy $X_i - X_{i+1} \ge 1$; (iii) every contact force is positive and every contact point $P_{i,j} = \tau_{i,j} / F_{i,j}$ lies in its contact interval $I_{i,j}(0)$, with $P_{1,0} = 0$; (iv) the force balance~\eqref{eq:force_bal} and the moment balance~\eqref{eq:moment_bal} of every block hold with residual $0$; and (v) $\sum_{i=1}^N X_i = 0$. Displayed values of $\underline{z}_{\mathbb{Q}}(N)$ are truncated downward, so they remain lower bounds. Appendix~\ref{app:RationalSeeds} lists all seeds and the exact overhangs.

\paragraph{Exact Spinal Optima.}
Combined with the exact spinal stacks of Appendix~\ref{app:ExactSpinal}, these stacks give the results of this paper that use no optimization solver. Paterson and Zwick computed the optimal loads, the shielded construction, and numerical configurations for these sizes, and noted that $N \in \{5, 7\}$ cannot realize the optimal loads with the same number of ordinary blocks~\cite[Sections~1 and~3]{paterson2009overhang}; what is new here are exact realizations and upper bounds that can be checked independently, including the exclusion of every other spine height and every face with a vanishing added load at the exception sizes.

\begin{theorem}[Exact Spinal Optima]
\label{thm:exact_spinal_main}
Consider stacks at zero stability margin ($\eta = 0$), and let $S(N)$ denote the maximum overhang of an $N$-block spinal stack. The following statements are verified in exact rational arithmetic or in interval arithmetic with outward rounding, without an optimization solver.
\begin{enumerate}
\item[(i)] For every $N \in \{4, \dots, 16\}$, the rational stack of Table~\ref{tab:rational_seeds} is feasible and single-root, and it is spinal both under the resting relation of Section~\ref{sec:SpinalCert} and in the sense of Paterson and Zwick; hence $\underline{z}_{\mathbb{Q}}(N) \le z^*(N)$.
\item[(ii)] $S(N) = S^*(N)$ for every $N \in \{4, 6, 8, \dots, 19\}$, which proves Conjecture~3.9 of Paterson and Zwick~\cite{paterson2009overhang} for these sizes.
\item[(iii)] At the exception sizes, $S(5) = (12 - \sqrt{30})/5 \approx 1.304554884990 < S^*(5)$ and $S(7) = \hat{S}_3(7) \approx 1.530045237061 < S^*(7)$, where $\hat{S}_3(7)$ is given in closed form by the positive root of $8\xi^4 - 14\xi - 7 = 0$.
\item[(iv)] The values $S(N)$ of~(ii) and~(iii) are attained by single-root spinal stacks of ordinary blocks and bound every spinal stack of our model, so $\underline{z}_{\mathbb{Q}}(N) \le S(N) \le z^*(N)$ for every $N \in \{4, \dots, 16\}$.
\item[(v)] Every single-root stack of four blocks is spinal, so $z^*(4) = S(4) = S^*(4) = (15 - 4\sqrt{2})/8$.
\end{enumerate}
\end{theorem}

\begin{proof}
Part~(i) is the zero-residual verification above; every contact of the rational stacks that carries a force has positive length. Parts~(ii) to~(iv) are Propositions~\ref{prop:exact_spinal} and~\ref{prop:exception_sizes} of Appendix~\ref{app:ExactSpinal}, checked by two independent programs. For part~(v), a non-spinal single-root stack has principal layer $\ell^* \ge 4$ and a layer with two blocks below it (Section~\ref{subsec:spinal_bound}, after Lemma~\ref{lem:maximal_contacts}), hence at least five blocks; so every single-root stack of four blocks is spinal, and parts~(ii) and~(iv) give $z^*(4) = S(4) = S^*(4)$. The maximizing height is $2$ (Table~\ref{tab:exact_heights}); there, \eqref{eq:pz_recurrence} gives $2 t_1^2 = 4$, so $t_1 = \sqrt{2}$, and $S^*(4) = \hat{S}_2(4) = (15 - 4\sqrt{2})/8$.
\end{proof}

The $60$-digit values of $S(N)$ give $S(N) - \underline{z}_{\mathbb{Q}}(N) \le 3.1 \times 10^{-8}$ for every $N \in \{4, \dots, 16\}$, and after rounding the values $\underline{z}_{\mathbb{Q}}(N)$ agree with the numerical configurations reported in Figures~3 and~4 of~\cite{paterson2009overhang} for all these sizes.

\section{Computational Results}
\label{sec:Computational}

\subsection{Numerical Certificates and Non-Spinal Exclusion for $N \in \{5, \dots, 16\}$}
\label{subsec:unconditional_certs}

We apply the framework of Section~\ref{sec:Algorithm} at $\eta = 0$; all computations in this section, and the scripts for the chains and the non-spinal branches in the companion repository, concern this case only. All subproblems and branches are solved on a workstation with an Intel Core Ultra~7 265K processor and 50~GB of memory running Linux, with Gurobi~13.0.1 and the settings of Section~\ref{subsec:TwoStagePipeline}. For tight certification, the cutoff is $\gamma_N := \underline{z}_{\mathbb{Q}}(N) + 10^{-5}$ rounded upward at the seventh decimal, except at $N = 15$, where we use the slightly larger $\gamma_{15} = 2.0929064$ instead of $2.0929063$. For non-spinal exclusion, the branches of Corollary~\ref{cor:nonspinal_exclusion} use the relative margin $\rho = 10^{-2}$ and the threshold $\gamma_{\mathrm{ns}}(N) := \lceil 0.99\, S^*(N) \rceil_7$, the value $0.99\, S^*(N)$ rounded upward at the seventh decimal; the solver receives the double-precision value of $0.99\, S^*(N)$, which does not exceed $\gamma_{\mathrm{ns}}(N)$. For $N = 4$, the optimum $(15 - 4\sqrt{2}) / 8 \approx 1.16789322$ is known in closed form~\cite{ainley1979,paterson2009maximum} and is proved again by Theorem~\ref{thm:exact_spinal_main}(v); rounded upward at the seventh decimal, it gives the anchor $\overline{z}(4) = 1.1678933$ of the chains below, and monolithic $\mathcal{M}_0$ at relative gap $10^{-9}$ reproduces it independently.

\paragraph{Main Numerical Result.}
Combining the exact lower bounds of Theorem~\ref{thm:exact_spinal_main} with the subproblems and the non-spinal branches that the solver reports closed yields our main numerical result. It rests on the closures reported by the solver, which we call numerical certificates. Table~\ref{tab:constants} collects the constants that enter it.

\begin{theorem}[Numerical Bounds and Spinal Structure]
\label{thm:main_certified_optimality}
Consider single-root stacks at zero stability margin ($\eta = 0$), with the model and solver settings of Sections~\ref{sec:Algorithm} and~\ref{subsec:unconditional_certs} (Gurobi~13.0.1).
\begin{enumerate}
\item[(a)] \textbf{Numerical upper bounds:} For every $N \in \{5, \dots, 16\}$, the solver reports every subproblem of the spinal-cut layer-profile decomposition closed within its numerical tolerances. Read as exact, these closures give $z^*(N) \le \overline{z}(N)$ with $\overline{z}(N) - \underline{z}_{\mathbb{Q}}(N) \le 1.02 \times 10^{-5}$ and relative optimality gaps between $4.67 \times 10^{-6}$ and $7.68 \times 10^{-6}$ ($4.71 \times 10^{-6}$ to $7.74 \times 10^{-6}$ under seven-decimal downward truncation of $\underline{z}_{\mathbb{Q}}(N)$). In particular, no single-root stack attains $S^*(N)$ at $N \in \{5, 7\}$: $z^*(5) \le 1.3045649 < S^*(5) \approx 1.3171618$ and $z^*(7) \le 1.5300553 < S^*(7) \approx 1.5318708$, which extends the strict inequalities $S(5) < S^*(5)$ and $S(7) < S^*(7)$ of Theorem~\ref{thm:exact_spinal_main} to all single-root stacks (the size $N = 3$ is not treated here). These intervals need not determine the last published digit; at $N = 5$, for instance, the interval contains the rounding boundary $1.304555$.
\item[(b)] \textbf{Numerical spinal structure:} For every $N \in \{5, \dots, 16\}$, the solver reports every non-spinal branch closed within its numerical tolerances. Read as exact, these closures show that every feasible single-root stack whose overhang exceeds $\gamma_{\mathrm{ns}}(N) = \lceil 0.99\, S^*(N) \rceil_7$ is spinal in the sense of Section~\ref{sec:SpinalCert}, which for every stack balanced in the sense of Paterson and Zwick implies spinality under their Definition~3.1~\cite{paterson2009overhang}.
\item[(c)] \textbf{Consequence of part~(b):} If the conclusion of part~(b) holds, then every globally optimal single-root stack is spinal, and $z^*(N) = S(N)$ for every $N \in \{5, \dots, 16\}$: $z^*(N) = S^*(N)$ for $N \in \{6, 8, \dots, 16\}$, $z^*(5) = (12 - \sqrt{30})/5$, and $z^*(7) = \hat{S}_3(7)$.
\end{enumerate}
\end{theorem}

\begin{proof}
Assuming that the reported solver closures are exact, part~(a) follows from Theorem~\ref{thm:spinal_cert_decomp} together with Table~\ref{tab:main_certs} and the base-singleton chain of Table~\ref{tab:constants}, and part~(b) from Corollary~\ref{cor:nonspinal_exclusion} together with Table~\ref{tab:main_certs} and the chain discussed below. For part~(c), let $X$ be a globally optimal single-root stack; since $\underline{z}_{\mathbb{Q}}(N) > \gamma_{\mathrm{ns}}(N)$, part~(b) makes $X$ spinal, so its support set consists of one block $b_t$ in each layer $t \in \{1, \dots, \ell^*\}$, with $b_1 = 1$ and $b_{\ell^*}$ the principal block, where $\ell^* \le H \le N$. For $t \in \{2, \dots, \ell^*\}$, every block on which $b_t$ rests belongs to the support set and lies in layer $t - 1$, so $b_t$ rests only on $b_{t-1}$; since a contact links a block only to a block on which it rests, each $b_t$ with $t \ge 2$ is supported exclusively by $b_{t-1}$. For $N \in \{6, 8, \dots, 16\}$, Proposition~\ref{prop:spinal_bound}(ii) with $p = 2$ and $\ell = \ell^*$ (or with $\ell = 1$ if $\ell^* = 1$) then gives $z^*(N) = X_{b_{\ell^*}} + W \le 2w\, S_{\ell^*}^*(N) \le S^*(N) = S(N)$, because $2w = 1$ at $\eta = 0$ and $S^*(N)$ is the maximum of $S_K^*(N)$ over $K \le N$; for $N \in \{5, 7\}$, Proposition~\ref{prop:exception_sizes} gives $z^*(N) \le S(N)$ directly. Since $S(N) \le z^*(N)$ by Theorem~\ref{thm:exact_spinal_main}, this proves part~(c).
\end{proof}

Parts~(a) and~(b) rest on the solver: the auditor of Section~\ref{sec:Algorithm} checks that every subproblem and branch is present, closed, and solved with the stated threshold, but it does not verify the dual bounds that closed them. The feasibility tolerance of the solver bounds the residuals of the constraints it accepts; it does not translate into an error bound on the reported upper bounds, and neither does the relative gap, which measures the width of the reported interval. The stability of the hardest subproblems under tighter tolerances and the positive reverse controls below are numerical robustness checks of these parts, not proofs. Part~(c) adds to part~(b) only the analytical bounds of Proposition~\ref{prop:spinal_bound}(ii) and Appendix~\ref{app:ExactSpinal}; its values are exact, but their identification with $z^*(N)$ is only as reliable as the floating-point closures of part~(b). The thresholds $\gamma_{\mathrm{ns}}(N)$ lie at least $1.3 \times 10^{-2}$ below $\underline{z}_{\mathbb{Q}}(N)$ for $N \ge 6$ and $5.6 \times 10^{-4}$ below it at $N = 5$, so the conclusion does not hinge on the exact value of the threshold; this margin is not an error bound for the solver either. Dual bounds that can be verified independently, such as safe bounds for the linear relaxations~\cite{neumaier2004} or certificates checked in exact arithmetic~\cite{cheung2017}, applied to the relaxations of the bilinear subproblems together with a proof that they cover the branches, would make parts~(a) to~(c) unconditional; we leave this to future work.

\begin{plaintable}[t]
\centering
\caption{Constants entering Theorems~\ref{thm:exact_spinal_main} and~\ref{thm:main_certified_optimality} and the chains that link consecutive sizes ($\eta = 0$). The exact lower bounds $\underline{z}_{\mathbb{Q}}(N)$ are truncated downward at the tenth decimal; the cutoffs $\gamma_N$ are the values passed to the solver, rounded upward at the seventh decimal, and the non-spinal thresholds $\gamma_{\mathrm{ns}}(N) = \lceil 0.99\, S^*(N) \rceil_7$ round the double-precision solver threshold $0.99\, S^*(N)$ upward at the seventh decimal; $S^*(N)$ is evaluated in $60$-digit arithmetic and rounded upward at the tenth decimal, and the column $S^*(N) - \underline{z}_{\mathbb{Q}}(N)$ is rounded upward. At $N = 4$, Theorem~\ref{thm:exact_spinal_main}(v) proves $z^*(4) = S^*(4) = (15 - 4\sqrt{2})/8$ without a solver. The spinal optimum is $S(N) = S^*(N)$ except at $N \in \{5, 7\}$, where $S(5) \approx 1.3045548850$ and $S(7) \approx 1.5300452371$ exceed $\underline{z}_{\mathbb{Q}}(N)$ by less than $10^{-9}$. The last two columns are the margins of the base-singleton conditions: \textsf{S2+} is $\underline{z}_{\mathbb{Q}}(N) - \overline{z}(N - 1) - 1/(2N)$ for~\eqref{eq:s2p_condition}, with $\overline{z}(4) = 1.1678933$ and $\overline{z}(N - 1) = \gamma_{N - 1}$ otherwise, and ns is $\gamma_{\mathrm{ns}}(N) - \gamma_{\mathrm{ns}}(N - 1) - 1/(2N)$ for~\eqref{eq:s2p_nonspinal_cond}, except at $N = 5$, marked $^{*}$, where it is $\gamma_{\mathrm{ns}}(5) - \overline{z}(4) - 1/10$; both margin columns are rounded to six decimals.}
\label{tab:constants}
\setlength{\tabcolsep}{4pt}
\resizebox{\textwidth}{!}{%
\begin{tabular}{@{}rccccccc@{}}
\toprule
$N$ & $\underline{z}_{\mathbb{Q}}(N)$ & $\gamma_N$ & $\gamma_{\mathrm{ns}}(N)$ & $S^*(N)$ & $S^*(N) - \underline{z}_{\mathbb{Q}}(N)$ & \textsf{S2+} margin & ns margin \\
\midrule
4 & 1.1678932107 & -- & -- & 1.1678932189 & $8.1 \times 10^{-9}$ & -- & -- \\
5 & 1.3045548840 & 1.3045649 & 1.3039902 & 1.3171617922 & $1.3 \times 10^{-2}$ & $0.036662$ & $0.036097^{*}$ \\
6 & 1.4367002109 & 1.4367103 & 1.4223333 & 1.4367002144 & $3.5 \times 10^{-9}$ & $0.048802$ & $0.035010$ \\
7 & 1.5300452363 & 1.5300553 & 1.5165521 & 1.5318707920 & $1.9 \times 10^{-3}$ & $0.021906$ & $0.022790$ \\
8 & 1.6315068053 & 1.6315169 & 1.6151918 & 1.6315068207 & $1.6 \times 10^{-8}$ & $0.038952$ & $0.036140$ \\
9 & 1.7152728131 & 1.7152829 & 1.6981201 & 1.7152728156 & $2.4 \times 10^{-9}$ & $0.028200$ & $0.027373$ \\
10 & 1.7871258598 & 1.7871359 & 1.7692547 & 1.7871258620 & $2.1 \times 10^{-9}$ & $0.021843$ & $0.021135$ \\
11 & 1.8587795195 & 1.8587896 & 1.8401918 & 1.8587795504 & $3.1 \times 10^{-8}$ & $0.026189$ & $0.025483$ \\
12 & 1.9250913613 & 1.9251014 & 1.9058405 & 1.9250913753 & $1.4 \times 10^{-8}$ & $0.024635$ & $0.023982$ \\
13 & 1.9845089017 & 1.9845190 & 1.9646639 & 1.9845089033 & $1.5 \times 10^{-9}$ & $0.020946$ & $0.020362$ \\
14 & 2.0382169409 & 2.0382270 & 2.0178348 & 2.0382169446 & $3.7 \times 10^{-9}$ & $0.017984$ & $0.017457$ \\
15 & 2.0928962987 & 2.0929064 & 2.0719674 & 2.0928963033 & $4.5 \times 10^{-9}$ & $0.021336$ & $0.020799$ \\
16 & 2.1438427861 & 2.1438528 & 2.1224044 & 2.1438427967 & $1.1 \times 10^{-8}$ & $0.019686$ & $0.019187$ \\
\bottomrule
\end{tabular}}
\end{plaintable}

\paragraph{Relation to the Model of Paterson and Zwick.}
Theorems~\ref{thm:exact_spinal_main} and~\ref{thm:main_certified_optimality} are stated for our model, whose closed contact intervals also admit forces at point contacts. Both transfer to the model of Paterson and Zwick (Appendix~\ref{app:pz_relation}); apart from the spinal bounds $S(N)$ of parts~(ii) to~(iv) of Theorem~\ref{thm:exact_spinal_main}, stacks with several blocks on the table are not covered.

\paragraph{Inductive Chains for the Base-Singleton Reduction.}
Restricting the enumeration to $\mathcal{P}_2(N)$ relies on Corollary~\ref{cor:s2p_elimination}, whose condition~\eqref{eq:s2p_condition} links each size to the previous one. For $N = 5$, the anchor is $\overline{z}(4)$, and for $N \ge 6$ it is the certified bound $\overline{z}(N - 1) = \gamma_{N - 1}$. The condition only validates fixing $n_2 = 2$ and never enters a subproblem, and the certificate for $N - 1$ does not depend on $N$, so the chain is acyclic. Its margin exceeds $0.01798$, with the minimum at $N = 14$ (Table~\ref{tab:constants}). For non-spinal exclusion, the general form of condition~\eqref{eq:s2p_nonspinal_cond}, with $\overline{z}(N - 1)$ in place of $\overline{z}_{\mathrm{ns}}(N - 1)$, fails at $N = 14$ and $N = 16$, by $2.40 \times 10^{-3}$ and $1.75 \times 10^{-3}$. The sharper form with $\overline{z}_{\mathrm{ns}}(N - 1) = \gamma_{\mathrm{ns}}(N - 1)$, available inductively, holds for every $N \in \{6, \dots, 16\}$ with margin above $0.01745$, again with the minimum at $N = 14$, and together with Lemma~\ref{lem:wlog_base_shift} it closes the non-spinal chain for $N \in \{5, \dots, 16\}$.

\paragraph{Certification and Non-Spinal Exclusion Results.}
Table~\ref{tab:main_certs} reports the certificates and the exclusions for $N \in \{5, \dots, 16\}$ next to the numerical configurations reported by Paterson and Zwick~\cite[Figures~3--4]{paterson2009overhang}. Every profile in $\mathcal{P}_2(N)$ is closed, either analytically before a model is built or, as reported by the solver, through infeasibility or the objective cutoff; the analytic closures are the energy-pruned profiles of Table~\ref{tab:profile_counts} together with $1$, $1$, and $2$ profiles without a candidate layer at $N = 14$, $15$, and $16$. For $N \le 9$, certification and exclusion each take less than one thread-second per size; from $N = 10$ to $N = 16$, certification grows from $2$ to $72{,}855$ thread-seconds ($17{,}523$~s wall time at $N = 16$). All non-spinal branches close for $N \le 15$. At $N = 16$, $10{,}902$ of the $10{,}905$ branches close within the per-branch limit of $7{,}200$~s; the other three are split by the lowest impure layer (Proposition~\ref{prop:t0_split}) into $14$ sub-branches, of which the three with $t_0 = \ell^*$ need no computation and the other $11$ close by cutoff in $16{,}244$ thread-seconds (Appendix~\ref{app:n16_split}).

\paragraph{Robustness Checks.}
Re-solving the three most expensive profiles at $N = 14$ and at $N = 15$ with feasibility, optimality, and integrality tolerances of $10^{-9}$ and the most conservative numerical setting of the solver closes all six by cutoff again, in $86$ to $870$~s each (against $110$ to $731$~s in the certificate runs); with the same settings, the $11$ sub-branches at $N = 16$ close by cutoff again in $23{,}443$ thread-seconds (against $16{,}244$). Positive reverse controls check the implementation against erroneous pruning. On the $N = 13$ profile $(1, 2, 2, 2, 3, 2, 1)$, which contains the optimal spinal stack, with family~\textsf{SP} active and the cutoff lowered to $1.9844089$, the solver with primal heuristics enabled, integer feasibility tolerance $10^{-6}$, and relative gap $10^{-4}$ recovers a feasible spinal stack of overhang $1.9845088923$, within $1.1 \times 10^{-8}$ of $S^*(13)$, in $36.4$~s; with the certificate settings, it finds no feasible stack within one hour, but its dual bound stops at $S^*(13)$ and the profile stays open, as it must. At $N = 7$, lowering the non-spinal threshold to $1.38$ recovers a feasible stack of overhang $1.3842371$ that is non-spinal under our resting relation; it transmits a small force through a point contact, which our model admits and that of Paterson and Zwick does not. The companion repository contains these stacks with positions, contact forces, and contact points.

\begin{plaintable}[t]
\centering
\caption{Numerical optimality certificates, with solver-reported upper bounds ($\gamma_N = \underline{z}_{\mathbb{Q}}(N) + 10^{-5}$ rounded upward at the seventh decimal, and $\gamma_{15} = 2.0929064$) and numerical non-spinal exclusion ($\rho = 10^{-2}$) for single-root stacks at $\eta = 0$ (Gurobi~13.0.1, single-thread subproblems, integer feasibility tolerance $10^{-9}$), compared with the numerical configurations reported by Paterson and Zwick~\cite[Figures~3--4]{paterson2009overhang}. The exact zero-residual rational lower bounds $\underline{z}_{\mathbb{Q}}(N)$, whose exact values are the fractions of Table~\ref{tab:rational_seeds}, are displayed truncated downward to seven decimal places; the cutoff $\gamma_N$ passed to the solver equals the displayed $\overline{z}(N)$, and relative gaps $(\overline{z}(N) - \underline{z}_{\mathbb{Q}}(N)) / \underline{z}_{\mathbb{Q}}(N)$ are computed conservatively from the displayed bounds (exact untruncated gaps range from $7.68 \times 10^{-6}$ at $N = 5$ to $4.67 \times 10^{-6}$ at $N = 16$). Every profile in $\mathcal{P}_2(N)$ is closed, analytically or as reported by the solver: A counts the profiles closed analytically before any model is built (energy test or no candidate layer), I and C those closed by the solver with a proof of infeasibility or by the cutoff. Every non-spinal branch is closed in the same sense: E counts the branches with $J_{\ell^*} = \emptyset$, which need no model, I those proved infeasible without branching, all of them in presolve, and C those closed by the cutoff. $^{\dagger}$At $N = 16$, three further branches reached the $7{,}200$~s limit and were closed through the $14$ sub-branches of Proposition~\ref{prop:t0_split} ($11$ by the cutoff, $3$ with $t_0 = \ell^*$); the thread-seconds include the time-limited attempts.}
\label{tab:main_certs}
\setlength{\tabcolsep}{3.0pt}
\resizebox{\textwidth}{!}{%
\begin{tabular}{@{}rrrrrrrrrrr@{}}
\toprule
 & & \multicolumn{7}{c}{Numerical optimality certificate ($\gamma_N = \underline{z}_{\mathbb{Q}}(N) + 10^{-5}$)} & \multicolumn{2}{c}{Non-spinal ($\rho = 10^{-2}$)} \\
\cmidrule(lr){3-9} \cmidrule(l){10-11}
$N$ & PZ (2009) & $\underline{z}_{\mathbb{Q}}(N)$, truncated & Reported $\overline{z}(N)$ & Rel. gap & Profiles A/I/C & Nodes & Thread-s & & Branches E/I/C & Thread-s \\
\midrule
 5 & 1.30455 & 1.3045548 & 1.3045649 & $7.74 \times 10^{-6}$ & 0/0/2 & 12 & $< 1$ & & 2/3/0 & $< 1$ \\
 6 & 1.4367 & 1.4367002 & 1.4367103 & $7.03 \times 10^{-6}$ & 1/0/3 & 3 & $< 1$ & & 3/7/0 & $< 1$ \\
 7 & 1.53005 & 1.5300452 & 1.5300553 & $6.60 \times 10^{-6}$ & 1/0/6 & 204 & $< 1$ & & 6/13/4 & $< 1$ \\
 8 & 1.63151 & 1.6315068 & 1.6315169 & $6.19 \times 10^{-6}$ & 1/0/11 & 85 & $< 1$ & & 11/28/9 & $< 1$ \\
 9 & 1.71527 & 1.7152728 & 1.7152829 & $5.89 \times 10^{-6}$ & 4/0/18 & 486 & $< 1$ & & 18/50/24 & $< 1$ \\
10 & 1.78713 & 1.7871258 & 1.7871359 & $5.65 \times 10^{-6}$ & 5/2/32 & 5{,}229 & 2 & & 34/96/61 & 4 \\
11 & 1.85878 & 1.8587795 & 1.8587896 & $5.43 \times 10^{-6}$ & 5/5/60 & 22{,}191 & 8 & & 65/191/140 & 12 \\
12 & 1.92509 & 1.9250913 & 1.9251014 & $5.25 \times 10^{-6}$ & 12/8/106 & 125{,}352 & 37 & & 114/360/292 & 46 \\
13 & 1.98451 & 1.9845089 & 1.9845190 & $5.09 \times 10^{-6}$ & 13/12/200 & 697{,}072 & 247 & & 212/674/643 & 320 \\
14 & 2.03822 & 2.0382169 & 2.0382270 & $4.96 \times 10^{-6}$ & 26/11/366 & 3{,}866{,}654 & 1{,}279 & & 378/1{,}217/1{,}348 & 1{,}344 \\
15 & 2.0929 & 2.0928962 & 2.0929064 & $4.87 \times 10^{-6}$ & 42/18/663 & 26{,}237{,}013 & 8{,}641 & & 682/2{,}255/2{,}745 & 16{,}499 \\
16 & 2.14384 & 2.1438427 & 2.1438528 & $4.71 \times 10^{-6}$ & 72/30/1{,}194 & 232{,}080{,}537 & 72{,}855 & & 1{,}226/4{,}229/5{,}447$^{\dagger}$ & 113{,}000$^{\dagger}$ \\
\bottomrule
\end{tabular}}
\end{plaintable}

\paragraph{Scalability and the Limits of the Spinal Regime.}
The number of profiles $|\mathcal{P}_2(N)|$ grows by a factor of about $1.8$ per block, but the certification effort grows faster, by factors of $5.2$, $6.8$, and $8.4$ in thread-seconds from $N = 13$ to $N = 16$, so the effort per profile grows as well. If the last factor persisted, a certificate for $N = 17$ ($2{,}325$ profiles) would need about $6 \times 10^{5}$ thread-seconds, which appears tractable on a parallel cluster, with the exact spinal stack of Appendix~\ref{app:ExactSpinal} as lower bound. At $N = 20$, the numerical configurations of Paterson and Zwick~\cite[Section~4]{paterson2006overhang}, \cite[Sections~3--4]{paterson2009overhang} include a non-spinal stack of overhang $2.32014$, above their spinal value of about $2.31899$ and our evaluation $S^*(20) \approx 2.3189945$; this comparison is numerical, since we have neither verified that stack in exact arithmetic nor proved $S(20) = S^*(20)$. Once non-spinal stacks exceed $S^*(N)$, as these configurations indicate, the pure-chain bound in~\eqref{eq:sp_cut} no longer lies below the cutoff, so the cuts, which remain valid, prune little, and certification for $N \ge 20$ will need relaxations aimed at non-spinal architectures.

\subsection{Methodological Progression and Cross-Verification Without Spinal Cuts}
\label{subsec:method_progression}

To quantify each methodological layer and to cross-check the certificates without family~\textsf{SP}, Table~\ref{tab:method_progression} compares successive stages at $N = 12$ and $N = 13$. On four threads, monolithic $\mathcal{M}_0$ reaches the relative gap $5.02 \times 10^{-4}$ at $N = 12$ after $6.86$~h ($98{,}749$ thread-seconds) and is still open after $72$~h, and at $N = 13$ its gap is still $2.9 \times 10^{-3}$ after $48$~h. Profile decomposition without family~\textsf{SP} closes all $126$ profiles at $N = 12$ to the same gap in $1{,}068$ thread-seconds, $92.5$ times less than $\mathcal{M}_0$, and all $225$ profiles at $N = 13$ in $13{,}266$ thread-seconds. At $N = 13$, the profile $(1, 2, 2, 2, 3, 2, 1)$ that contains the optimal spinal stack dominates this effort; a separate split of this profile by candidate layer, closed by cutoff in $15{,}587$~s at the critical layer~$5$, confirms its closure independently. The final pipeline, which adds family~\textsf{SP}, the energy test, and single-thread subproblems at the tight cutoff, reduces thread-seconds by factors of $28.9$ at $N = 12$ and $53.7$ at $N = 13$, and nodes by a factor of $83.0$ at $N = 13$, while it tightens the gap by two orders of magnitude, to $5.21 \times 10^{-6}$ and $5.09 \times 10^{-6}$. These factors compare complete configurations; the matched comparison below isolates family~\textsf{SP}.

\paragraph{Matched Comparison of the Spinal Cuts.}
Table~\ref{tab:sp_ablation} repeats the tight certificates at $N \in \{11, 12, 13\}$ without family~\textsf{SP} and leaves everything else unchanged. When $J_\ell = \emptyset$, the cut~\eqref{eq:sp_cut} is a bound on $X_{s_\ell}$, which our implementation imposes as a variable bound, dropping layer~$\ell$ from the candidates when the bound does not exceed $\gamma_N$; the comparison therefore measures family~\textsf{SP} as a whole. At $N = 11$, both configurations close all $70$ profiles, so the run without \textsf{SP} is a second certificate of $\overline{z}(11)$, and \textsf{SP} reduces its effort by factors of $198$ in thread-seconds and $685$ in nodes. At $N = 12$ and $N = 13$, $2$ and $4$ profiles remain open without \textsf{SP} at the time limit, with relative gaps $7.1 \times 10^{-5}$ and $2.1 \times 10^{-3}$, whereas with \textsf{SP} each of them closes in at most $23$~s; on the profiles that close in both configurations, \textsf{SP} reduces thread-seconds by factors of $16$ and $57$ and nodes by factors of $27$ and $109$. The totals without \textsf{SP} at $N = 12$ and $N = 13$ include the time-limited profiles and are therefore only lower bounds on the effort of a complete certificate. With \textsf{SP}, $4$, $6$, and $8$ more profiles end with a proof of infeasibility at $N = 11$, $12$, and $13$; the energy test closes the same profiles in both configurations, so these are proofs of the solver. Since the solver treats the cutoff as a constraint on the objective and the cuts remove no feasible stack (Proposition~\ref{prop:sp_cuts}), such a proof shows, as a closure by the cutoff does, that no stack of the profile exceeds $\gamma_N$.

\begin{plaintable}[t]
\centering
\caption{Methodological progression from monolithic spatial branch-and-bound to spinal-cut layer-profile decomposition at $N = 12$ and $N = 13$ ($\eta = 0$, Gurobi~13.0.1). The upper bound $\overline{z}$ is the solver's final dual bound for $\mathcal{M}_0$ and the cutoff for the decompositions, rounded upward at the seventh decimal; all relative gaps are measured against the same exact lower bound $\underline{z}_{\mathbb{Q}}(N)$, without truncation, so the gap $5.21 \times 10^{-6}$ at $N = 12$ is slightly smaller than the value $5.25 \times 10^{-6}$ of Table~\ref{tab:main_certs}, which uses the truncated bound. For monolithic $\mathcal{M}_0$ at $N = 12$, the $5 \times 10^{-4}$ gap milestone ($6.86$~h) and the complete $72$-hour run are listed separately. Threads are per run for $\mathcal{M}_0$ and per subproblem for the decompositions, which run four (without \textsf{SP}) or three (with \textsf{SP}) subproblems in parallel; at $N = 13$ without \textsf{SP}, the profile $(1, 2, 2, 2, 3, 2, 1)$ was solved separately on three threads. The cutoffs without \textsf{SP} are $\underline{z}(N)(1 + 5 \times 10^{-4})$ with the feasible values available at the time, namely $1.9260577476$ at $N = 12$ and $1.9855055567$ at $N = 13$, displayed rounded to seven decimals, whereas $\overline{z}$ is rounded upward; $\mathcal{M}_0$ uses no cutoff. The runs without \textsf{SP} were made without the energy test, so their infeasible profiles were all proved infeasible by the solver. The stages differ in model, threads, screening, and cutoff at once, so the table shows the progression of complete configurations; Table~\ref{tab:sp_ablation} compares the spinal cuts at the same cutoff.}
\label{tab:method_progression}
\setlength{\tabcolsep}{3.5pt}
\resizebox{\textwidth}{!}{%
\begin{tabular}{@{}clcccrrrl@{}}
\toprule
$N$ & Stage & Threads & Cutoff & Upper bound $\overline{z}$ & Rel. gap & Nodes & Thread-s & Closure \\
\midrule
12 & $\mathcal{M}_0$, $5 \times 10^{-4}$ milestone & 4 & none & 1.9260573 & $5.02 \times 10^{-4}$ & 188{,}122{,}744 & 98{,}749 & at $6.86$\,h \\
12 & $\mathcal{M}_0$, full run & 4 & none & 1.9251508 & $3.09 \times 10^{-5}$ & 1{,}755{,}164{,}331 & 1{,}036{,}800 & open at $72$\,h \\
12 & Profiles, no \textsf{SP} & 1 & 1.9260577 & 1.9260578 & $5.02 \times 10^{-4}$ & 5{,}710{,}202 & 1{,}068 & 126/126 \\
12 & Profiles, \textsf{SP} & 1 & 1.9251014 & 1.9251014 & $5.21 \times 10^{-6}$ & 125{,}352 & 37 & 126/126 (106 C, 20 I) \\
\midrule
13 & $\mathcal{M}_0$ & 4 & none & 1.9902762 & $2.91 \times 10^{-3}$ & 855{,}637{,}535 & 690{,}670 & open at $48$\,h \\
13 & Profiles, no \textsf{SP} & 1 & 1.9855056 & 1.9855056 & $5.02 \times 10^{-4}$ & 57{,}847{,}794 & 13{,}266 & 225/225 (220 C, 5 I) \\
13 & Profiles, \textsf{SP} & 1 & 1.9845190 & 1.9845190 & $5.09 \times 10^{-6}$ & 697{,}072 & 247 & 225/225 (200 C, 25 I) \\
\bottomrule
\end{tabular}}
\end{plaintable}

\begin{plaintable}[t]
\centering
\caption{Matched comparison of the tight certificates with and without family~\textsf{SP} ($\eta = 0$, Gurobi~13.0.1, cutoff $\gamma_N$, one thread per subproblem, otherwise identical options). C and I count the profiles closed by cutoff and proved infeasible (I includes the $5$, $12$, and $13$ profiles closed by the energy test, the same in both configurations); open profiles reached the time limit of $3{,}600$~s ($N \le 12$) or $14{,}400$~s ($N = 13$), and the last column gives the largest dual bound among them. Wall times are for seven parallel workers without \textsf{SP} and three with \textsf{SP} and are not used for the comparison.}
\label{tab:sp_ablation}
\begin{tabular}{@{}ccrrrrrr@{}}
\toprule
$N$ & \textsf{SP} & Closed (C + I) & Open & Nodes & Thread-s & Wall (s) & Open bound \\
\midrule
11 & no & 70 (64 + 6) & 0 & 15{,}194{,}632 & 1{,}578 & 1{,}096 & -- \\
11 & yes & 70 (60 + 10) & 0 & 22{,}191 & 8 & 4 & -- \\
12 & no & 124 (110 + 14) & 2 & 58{,}943{,}212 & 7{,}623 & 3{,}616 & 1.9252284 \\
12 & yes & 126 (106 + 20) & 0 & 125{,}352 & 37 & 15 & -- \\
13 & no & 221 (204 + 17) & 4 & 362{,}115{,}068 & 67{,}606 & 14{,}706 & 1.9886878 \\
13 & yes & 225 (200 + 25) & 0 & 697{,}072 & 247 & 88 & -- \\
\bottomrule
\end{tabular}
\end{plaintable}

\paragraph{Why the Spinal Cuts Help.}
For a candidate layer with $2w\,\hat{S}_\ell(N) \le \gamma_N$, every configuration above the cutoff must satisfy $\Sigma_\ell \ge 1$, so the search is confined to configurations whose rightmost chain rests on a non-rightmost block below the candidate. This holds for every candidate layer at the non-exception sizes, where $\hat{S}(N) = S^*(N) < \gamma_N$; at $N = 5$ (layer~$3$) and $N = 7$ (layer~$4$), the certificate rests on closing the subproblems (Theorem~\ref{thm:spinal_cert_decomp}(ii)). We attribute the reductions in Table~\ref{tab:method_progression} to the removal of pure-chain configurations above the loaded-spine bound, which in profiles that contain near-optimal spinal stacks keep the dual bound close to the cutoff. The region $\Sigma_\ell \ge 1$ also contains spinal stacks whose spine passes through a non-rightmost block, so Corollary~\ref{cor:nonspinal_exclusion} bounds only its non-spinal part, and lower thresholds locate that part more precisely: at $N = 9$, the single-root non-spinal optimum lies in $[1.63188, 1.66382)$, between $3.0\%$ and $4.9\%$ below $z^*(9)$. The lower end is an exactly balanced rational non-spinal stack of overhang $1.6318870324\ldots$, obtained from a solver solution at threshold $1.60$ by fixing its contact forces at rational values and solving the remaining conditions, which are then linear, in exact arithmetic; all its loaded contacts have positive length, the shortest, between blocks~$4$ and~$3$, about $1.66 \times 10^{-7}$ but positive in exact arithmetic, so it is non-spinal also in the sense of Paterson and Zwick. The upper end, valid within solver tolerances, is the threshold $0.97\, S^*(9)$ rounded upward, at which all non-spinal branches close; stacks with $n_2 = 1$ are excluded by condition~\eqref{eq:s2p_nonspinal_cond}, because the branches at $N = 8$ close at threshold $1.60 > \overline{z}(7) + 1/16$ and $1.60 + 1/18 < 0.97\, S^*(9)$.

\section{Conclusion}
\label{sec:Conclusion}

We have developed a global optimization framework for the single-root maximum overhang problem, a nonconvex mixed-integer bilinear program over discrete support topologies and continuous block coordinates. Structural reductions, a base-shift lemma, an energy criterion, an exact layer-profile decomposition, spinal cuts derived from a loaded-spine chain bound, and exact rational lower bounds replace the weak global McCormick boxes of the monolithic model by narrow profile-specific boxes and, in our experiments, remove the stalling of spatial branch-and-bound near spinal optima.

At zero stability margin and without any solver, our exact stacks and the bounds of Appendix~\ref{app:ExactSpinal} prove the equality $S(N) = S^*(N)$ of Conjecture~3.9 of Paterson and Zwick~\cite{paterson2009overhang} for $N \in \{4, 6, 8, \dots, 19\}$ and determine the spinal optima $S(5) = (12 - \sqrt{30})/5$ and $S(7)$ in closed form. Within the single-root class, the solver reports every subproblem and every non-spinal branch closed for $N \in \{5, \dots, 16\}$; read as exact, these numerical certificates bound the optimal overhang within $1.02 \times 10^{-5}$ of these values and show that every stack above $\lceil 0.99\, S^*(N) \rceil_7$ is spinal, and if these solver closures are correct, $z^*(N) = S(N)$ for every $N \in \{5, \dots, 16\}$.

The combination of an analytical equilibrium bound with sparse cuts that isolate an analytically understood regime inside a discrete decomposition may extend to other bilinear programs governed by static equilibrium, such as palletizing and container loading with stability constraints or limit analysis of dry-jointed masonry; whether it does remains to be verified.

\paragraph{Open Problems.}
\begin{enumerate}
\item \textbf{Solver-independent upper bounds:} Theorem~\ref{thm:main_certified_optimality} rests on floating-point closures. Dual bounds that can be checked independently, such as safe bounds for the linear relaxations~\cite{neumaier2004} or certificates verified in exact arithmetic~\cite{cheung2017}, combined with a proof that the relaxations and branches cover each bilinear subproblem, would make the identification $z^*(N) = S(N)$ unconditional, starting with the smallest sizes, where most branches close in presolve.
\item \textbf{Multi-root foundations:} All optimal stacks displayed by Paterson and Zwick~\cite[Figures~3--5]{paterson2009overhang} rest on a single block on the table, and our certificates hold only within this class. Whether several table-contacting blocks can ever do better at finite $N$ remains open.
\item \textbf{Certifying $N \ge 17$ and the transition at $N = 20$:} Appendix~\ref{app:ExactSpinal} determines the spinal optima for $N \le 19$ exactly; certifying that no single-root stack exceeds them at $N = 17$ is the next target, estimated at about $6 \times 10^{5}$ thread-seconds. At $N = 20$, a numerical configuration of Paterson and Zwick exceeds $S^*(20)$ with a non-spinal stack~\cite{paterson2006overhang,paterson2009overhang}. There, the spinal cuts no longer separate the optimum from the rest of a profile, the non-spinal exclusion has no threshold below the optimum, and from $N = 21$ on, the base-singleton chain needs upper bounds at sizes whose optima are themselves non-spinal; which relaxations can replace these mechanisms is open.
\item \textbf{Three-dimensional and frictional models:} Extending global certification to three-dimensional stacks, including rotated arrangements~\cite[Section~6]{paterson2009maximum}, and to frictional contact~\cite[Section~6]{paterson2009overhang}, where the asymptotic order of the overhang changes, is a rich direction for nonconvex global optimization.
\end{enumerate}

\backmatter

\section*{Declarations}

\paragraph{Funding.}
This work was supported in part by JSPS KAKENHI Grant Number 26K02969 and the Logistics Research Grant from the SBS Kamata Foundation.

\paragraph{Competing interests.}
The authors have no relevant financial or non-financial interests to disclose.

\paragraph{Data and code availability.}
All source code (layer-profile decomposition, spinal cuts, non-spinal exclusion, and certificate auditing), the per-profile and per-branch certification summaries, the solver logs and model files of key instances, the exact rational stack data for $N \in \{4, \dots, 16\}$, the standalone exact-arithmetic verification script, the exact re-check of every analytic screening decision, and the two verifiers of the exact spinal optima of Appendix~\ref{app:ExactSpinal} with their outputs and the replay of their box coverings are publicly available at \url{https://github.com/CyouRyuu/Max_Overhang} and archived on Zenodo (version 1.0.6, Git tag v1.0.6, commit 4e7e453) at \url{https://doi.org/10.5281/zenodo.24576588}; the concept DOI \url{https://doi.org/10.5281/zenodo.24518950} refers to all versions.

\paragraph{Author contributions.}
Aolong Zha conceptualized the study, developed the optimization models and the certification framework, implemented the algorithms, conducted the computational experiments, and wrote the manuscript. Qianwei Zhuang contributed to the model validation and data analysis. Qiong Chang supervised the project, acquired funding, and revised the manuscript. All authors read and approved the final manuscript.

\paragraph{Use of artificial intelligence tools.}
During this work, the authors used the large language models Claude Opus (Anthropic), Gemini Neon (Google), and ChatGPT Astra (OpenAI) as assistants for exploring candidate proof steps, drafting and reviewing analysis and verification scripts, cross-checking drafts of the manuscript, and improving English exposition. The models did not generate any result that was not independently established by the authors: every theorem and proof was checked by the authors, every computational certificate was produced by the authors' runs of the published code with Gurobi, the exact rational lower bounds were verified in exact arithmetic by the scripts in the companion repository, and the exact spinal optima were reproduced by the authors with both verifiers, which decide the inequalities by different methods. The authors take full responsibility for the content of this publication.

\clearpage
\begin{appendices}
\renewcommand{\theHtable}{app.\arabic{table}}\renewcommand{\theHfigure}{app.\arabic{figure}}

\section{Exact Rational Stack Seeds and Reconstruction Index}
\label{app:RationalSeeds}

This appendix catalogs the compact rational seed parameters that generate the zero-residual feasible stacks for all $N \in \{4, \dots, 16\}$ at zero stability margin ($\eta = 0$, $W = 1/2$) and gives the closed-form reconstruction of Section~\ref{subsec:RationalConstructor}. Every configuration listed in Table~\ref{tab:rational_seeds} has been reconstructed and audited independently in exact arbitrary-precision rational arithmetic over $\mathbb{Q}$ with zero residual across all force-balance equations, moment-balance equations, contact-interval inclusions, same-layer non-overlap inequalities, and single-root table-clearance conditions. Every contact that carries a positive force has an overlap of positive length, at least $3.5 \times 10^{-3}$, so the stacks are also balanced in the sense of Paterson and Zwick~\cite{paterson2009overhang}, in whose model point contacts carry no force.

\subsection{Reconstruction Procedure}
\label{subsec:ReconIndex}

For $N \ne 5$, the shield block $B_r'$ ($r \in \{1, \dots, k - 1\}$) lies in layer $k - r + 1$ beside $B_r$, the bridge block $B_0'$ in layer $k + 1$, and the tower blocks $T_1, \dots, T_{v_1}$ rest on $B_1'$ in layers $k + 1, \dots, k + v_1$. Each shield block $B_r'$ has center $X(B_r') = y_r + 1/2$ and applies downward force $\omega_{r+1}$ and torque $\omega_{r+1} x_{r+1}$ to spine block $B_{r+1}$ at its left edge $x_{r+1}$. Initializing $u_k := 0$ and $\tau_k^{\mathrm{sh}} := 0$, vertical force and moment balance on $B_r'$ propagate bottom-up for $r = k - 1, \dots, 2$ to give the exact rational force $u_r$ and torque $\tau_r^{\mathrm{sh}}$ exerted by $B_{r-1}'$ on $B_r'$:
\begin{equation}
\label{eq:rat_shield_prop}
u_r = u_{r+1} + \omega_{r+1} - 1,
\qquad
\tau_r^{\mathrm{sh}} = \tau_{r+1}^{\mathrm{sh}} + \omega_{r+1} x_{r+1} - \left(y_r + \frac{1}{2}\right)
\qquad (r = k - 1, \dots, 2).
\end{equation}
At the top shield level $r = 1$, the bridge block $B_0'$ applies force $\omega_1$ and torque $\omega_1 x_1$ to $B_1$ at $x_1$ and transmits the remaining weight $u_1 := 1 - \omega_1 > 0$ onto $B_1'$. Defining the net upward moment demand on $B_1'$ by
\begin{equation}
\label{eq:rat_m1_demand}
M_1 := \tau_2^{\mathrm{sh}} + \omega_2 x_2 - \left(y_1 + \frac{1}{2}\right),
\end{equation}
exact moment equilibrium of $B_1'$ and $B_0'$ yields the closed-form rational bridge and tower coordinates:
\begin{align}
\text{if } v_1 = 0:\quad &
\tau_1^{\mathrm{sh}} = M_1, \qquad
z_1 = \frac{M_1}{u_1}, \qquad
X(B_0') = M_1 + \omega_1 x_1; \label{eq:rat_top_v0} \\
\text{if } v_1 \ge 1:\quad &
\tau_1^{\mathrm{sh}} = u_1 z_1, \qquad
X(B_0') = u_1 z_1 + \omega_1 x_1, \qquad
p_{v_1} = \frac{M_1 - u_1 z_1}{v_1}. \label{eq:rat_top_vpos}
\end{align}
When $v_1 \ge 1$, all $v_1$ tower blocks $T_1, \dots, T_{v_1}$ share the exact rational center $X(T_s) = p_{v_1}$ ($s = 1, \dots, v_1$), with bottom tower contact $(F_{T_1, B_1'}, \tau_{T_1, B_1'}) = (v_1, v_1 p_{v_1})$ and vertical inter-tower contacts $(F_{T_s, T_{s-1}}, \tau_{T_s, T_{s-1}}) = (v_1 - s + 1,\, (v_1 - s + 1) p_{v_1})$ for $s = 2, \dots, v_1$.

For $N = 5$, the spine has $d_1 = 1/2$ ($\omega_1 = 0$, $t_1 = 1$), and the counterweights $C_2$ in layer~$2$ and $C_1$ in layer~$3$ supply $\omega_2 = t_2 - 2$ to $B_2$ and $\omega_3 = 4 - t_2$ to $B_3$. Fixing rational seeds $t_2 \in \mathbb{Q}$ (approximating $\sqrt{15/2}$) and $X(C_2) \in \mathbb{Q}$ determines all remaining quantities in closed form over $\mathbb{Q}$:
\begin{equation}
\label{eq:rat_n5_closed}
\begin{gathered}
\omega_2 = t_2 - 2, \quad
\omega_3 = 4 - t_2, \quad
u_2 = 1 - \omega_2, \\
\tau_{C_1, C_2} = \omega_3 x_3 - X(C_2), \quad
X(C_1) = \omega_2 x_2 + \tau_{C_1, C_2},
\end{gathered}
\end{equation}
with $d_1 = 1/2$, $d_2 = (\omega_2 + 1/2)/t_2$, $d_3 = (\omega_3 + 1/2)/5$, spine coordinates from~\eqref{eq:rat_spine_steps}, and contacts $(C_2, B_3)$ carrying $(\omega_3, \omega_3 x_3)$, $(C_1, B_2)$ carrying $(\omega_2, \omega_2 x_2)$, and $(C_1, C_2)$ carrying $(u_2, \tau_{C_1, C_2})$.

In summary, given a seed tuple from Table~\ref{tab:rational_seeds}, every coordinate, contact force, and contact torque is obtained in $\mathbb{Q}$ as follows. For $N \ne 5$, compute $\omega_r$, $d_r$, $x_r$, and the spine contacts from~\eqref{eq:rat_spine_steps} and~\eqref{eq:rat_overhang} using the loads $(t_1, \dots, t_{k-1})$; place the shield blocks at $X(B_r') = y_r + 1/2$ and propagate their contacts by~\eqref{eq:rat_shield_prop}; finally, obtain the bridge and the tower from~\eqref{eq:rat_m1_demand} together with~\eqref{eq:rat_top_v0} when $v_1 = 0$, or with~\eqref{eq:rat_top_vpos} and the seed $z_1$ when $v_1 \ge 1$. For $N = 5$, use $t_1 = 1$, the seeds $(t_2, X(C_2))$, and~\eqref{eq:rat_n5_closed}. Every resulting stack is sorted canonically (Proposition~\ref{prop:single_root}) and passes the exact checks (i) to (v) of Section~\ref{subsec:RationalConstructor}, including same-layer separation of at least one unit.

\subsection{Complete Rational Seed Table for $N \in \{4, \dots, 16\}$}
\label{subsec:SeedTable}

Table~\ref{tab:rational_seeds} lists the architecture $(k, v_1)$, exact rational load vector $(t_1, \dots, t_{k-1})$, shield left-edge vector $(y_1, \dots, y_{k-1})$, top bridge contact seed $z_1$, and resulting verified overhang $\underline{z}_{\mathbb{Q}}(N)$ for all $N \in \{4, \dots, 16\}$. Across all non-exception sizes $N \in \{4, 6, 8, \dots, 16\}$, the exact rational overhangs lie within $1.5 \times 10^{-9}$ to $3.1 \times 10^{-8}$ of $S^*(N)$, and since these stacks are spinal, both under the resting relation of Section~\ref{sec:SpinalCert} and in the sense of Paterson and Zwick, and $S^*(N)$ bounds every spinal stack~\cite[Section~3]{paterson2009overhang} (Proposition~\ref{prop:spinal_bound} for single-root stacks), $0 \le S^*(N) - S(N) \le 3.1 \times 10^{-8}$ already follows from these stacks; Appendix~\ref{app:ExactSpinal} shows that $S(N) = S^*(N)$ by replacing their loads with exact algebraic ones. At the exception sizes $N = 5$ and $N = 7$, the rational overhangs lie within $9.3 \times 10^{-10}$ of $S(5) = \frac{12}{5} - 2\sqrt{3/10}$ and within $7.1 \times 10^{-10}$ of $S(7) = S_3^*(7)$, respectively (Proposition~\ref{prop:exception_sizes}).

\begin{plaintable}[t]
\centering
\caption{Exact rational seed parameters generating zero-residual single-root feasible stacks for $N \in \{4, \dots, 16\}$ at $\eta = 0$. For $N \ne 5$, $(k, v_1)$ denotes spine height $k$ and counterweight tower size $v_1 = N - 2k$; when $v_1 = 0$, $z_1$ is determined automatically by~\eqref{eq:rat_top_v0}. For the bridgeless exception $N = 5$ (profile $(1, 2, 2)$), the listed seeds are $t_2$ and $X(C_2)$ for~\eqref{eq:rat_n5_closed}. The final column reports the verified exact rational overhang truncated downward at the tenth decimal (exact reduced fractions are listed in Section~\ref{subsec:SeedTable}).}
\label{tab:rational_seeds}
\footnotesize
\setlength{\tabcolsep}{3.5pt}
\resizebox{\textwidth}{!}{%
\begin{tabular}{@{}cclllr@{}}
\toprule
$N$ & $(k, v_1)$ & Intermediate loads $(t_1, \dots, t_{k-1})$ & Shield edges $(y_1, \dots, y_{k-1})$ & $z_1$ & Verified $\underline{z}_{\mathbb{Q}}(N)$ \\
\midrule
 4 & $(2, 0)$ & $\bigl(\frac{707}{500}\bigr)$ & $(-1)$ & via~\eqref{eq:rat_top_v0} & 1.1678932107 \\[3pt]
 5 & $(3, \text{br.})$ & $t_2 = \frac{5477}{2000}$ ($\omega_1 = 0$) & $X(C_2) = -\frac{18}{25}$ & N/A & 1.3045548840 \\[3pt]
 6 & $(3, 0)$ & $\bigl(\frac{639}{500},\, \frac{1633}{500}\bigr)$ & $\bigl(-1,\, -\frac{6}{5}\bigr)$ & via~\eqref{eq:rat_top_v0} & 1.4367002109 \\[3pt]
 7 & $(3, 1)$ & $\bigl(\frac{2679}{2000},\, \frac{7177}{2000}\bigr)$ & $\bigl(-\frac{3}{4},\, -\frac{6}{5}\bigr)$ & $\frac{1}{5}$ & 1.5300452363 \\[3pt]
 8 & $(4, 0)$ & $\bigl(\frac{154}{125},\, \frac{6069}{2000},\, \frac{2127}{400}\bigr)$ & $\bigl(-\frac{17}{20},\, -1,\, -\frac{7}{5}\bigr)$ & via~\eqref{eq:rat_top_v0} & 1.6315068053 \\[3pt]
 9 & $(4, 1)$ & $\bigl(\frac{631}{500},\, \frac{6371}{2000},\, \frac{11519}{2000}\bigr)$ & $\bigl(-\frac{13}{20},\, -1,\, -\frac{7}{5}\bigr)$ & $\frac{11}{40}$ & 1.7152728131 \\[3pt]
10 & $(4, 2)$ & $\bigl(\frac{129}{100},\, \frac{416}{125},\, \frac{6187}{1000}\bigr)$ & $\bigl(-\frac{11}{20},\, -1,\, -\frac{7}{5}\bigr)$ & $\frac{2}{5}$ & 1.7871258598 \\[3pt]
11 & $(5, 1)$ & $\bigl(\frac{154}{125},\, \frac{607}{200},\, \frac{5317}{1000},\, 8\bigr)$ & $\bigl(-\frac{11}{20},\, -\frac{9}{10},\, -\frac{6}{5},\, -\frac{303}{200}\bigr)$ & $\frac{2}{5}$ & 1.8587795195 \\[3pt]
12 & $(5, 2)$ & $\bigl(\frac{156}{125},\, \frac{623}{200},\, \frac{111}{20},\, \frac{213}{25}\bigr)$ & $\bigl(-\frac{1}{2},\, -\frac{4}{5},\, -\frac{6}{5},\, -\frac{3}{2}\bigr)$ & $\frac{17}{40}$ & 1.9250913613 \\[3pt]
13 & $(5, 3)$ & $\bigl(\frac{1263}{1000}, \frac{319}{100}, \frac{2309}{400}, \frac{903}{100}\bigr)$ & $\bigl(-\frac{3}{5}, -\frac{31}{50}, -\frac{11}{10}, -\frac{29}{20}\bigr)$ & $\frac{2}{5}$ & 1.9845089017 \\[3pt]
14 & $(5, 4)$ & $\bigl(\frac{1277}{1000}, \frac{1631}{500}, \frac{1497}{250}, \frac{19061}{2000}\bigr)$ & $\bigl(-\frac{1}{2}, -\frac{3}{5}, -\frac{21}{20}, -\frac{7}{5}\bigr)$ & $\frac{1}{2}$ & 2.0382169409 \\[3pt]
15 & $(6, 3)$ & $\bigl(\frac{1239}{1000}, \frac{6141}{2000}, \frac{2711}{500}, \frac{3293}{400}, \frac{2289}{200}\bigr)$ & $\bigl(-\frac{1}{2}, -\frac{101}{200}, -\frac{19}{20}, -\frac{13}{10}, -\frac{8}{5}\bigr)$ & $\frac{1}{2}$ & 2.0928962987 \\[3pt]
16 & $(6, 4)$ & $\bigl(\frac{156}{125}, \frac{779}{250}, \frac{2777}{500}, \frac{17061}{2000}, \frac{24041}{2000}\bigr)$ & $\bigl(-\frac{9}{20}, -\frac{1}{2}, -\frac{9}{10}, -\frac{5}{4}, -\frac{31}{20}\bigr)$ & $\frac{11}{20}$ & 2.1438427861 \\
\bottomrule
\end{tabular}}
\end{plaintable}

For complete traceability, the reduced exact rational overhangs produced by Table~\ref{tab:rational_seeds} for $N = 4$ through $N = 16$ are
\[
\underline{z}_{\mathbb{Q}}(4) = \frac{1{,}651{,}401}{1{,}414{,}000}, \quad
\underline{z}_{\mathbb{Q}}(5) = \frac{71{,}450{,}471}{54{,}770{,}000}, \quad
\underline{z}_{\mathbb{Q}}(6) = \frac{21{,}115{,}183}{14{,}697{,}000},
\]
\[
\underline{z}_{\mathbb{Q}}(7) = \frac{411{,}858{,}436{,}609}{269{,}180{,}562{,}000}, \quad
\underline{z}_{\mathbb{Q}}(8) = \frac{247{,}112{,}620{,}283}{151{,}462{,}819{,}200}, \quad
\underline{z}_{\mathbb{Q}}(9) = \frac{158{,}860{,}140{,}542{,}171}{92{,}615{,}086{,}838{,}000},
\]
\[
\underline{z}_{\mathbb{Q}}(10) = \frac{1{,}483{,}400{,}102{,}771}{830{,}047{,}920{,}000}, \quad
\underline{z}_{\mathbb{Q}}(11) = \frac{335{,}947{,}378{,}939}{180{,}735{,}464{,}000}, \quad
\underline{z}_{\mathbb{Q}}(12) = \frac{12{,}287{,}515{,}301}{6{,}382{,}821{,}900},
\]
and
\[
\underline{z}_{\mathbb{Q}}(13) = \frac{120{,}400{,}867{,}309{,}071}{60{,}670{,}358{,}899{,}150}, \quad
\underline{z}_{\mathbb{Q}}(14) = \frac{23{,}072{,}966{,}927{,}875{,}223}{11{,}320{,}172{,}285{,}996{,}000},
\]
\[
\underline{z}_{\mathbb{Q}}(15) = \frac{58{,}035{,}616{,}013{,}404{,}873}{27{,}729{,}809{,}664{,}879{,}000}, \quad
\underline{z}_{\mathbb{Q}}(16) = \frac{41{,}648{,}849{,}852{,}234{,}795{,}719}{19{,}427{,}194{,}065{,}076{,}512{,}000}.
\]
Per-block coordinates and per-contact forces of all thirteen stacks are provided in the companion repository.

\section{Technical Proofs and Supplementary Remarks}
\label{app:proofs}

This appendix collects, in their order of appearance, the proofs omitted from Sections~\ref{sec:Formulation} to~\ref{sec:SpinalCert}, together with two supplementary remarks, the comparison with the model of Paterson and Zwick, the solver settings of the reported runs, and the sub-branch details of the $N = 16$ non-spinal exclusion.

\subsection{Proof of Lemma~\ref{lem:upper_dag_base}}
\label{app:upper_dag_proof}

\begin{proof}[Proof of Lemma~\ref{lem:upper_dag_base}]
For part~(i), if $T_{i,j} = 1$ with $i, j \in \mathcal{N}$, then~\eqref{eq:M0_layer_link} requires $L_i = L_j + 1 > L_j$, whereas canonical ordering gives $L_i \le L_j$ whenever $j \ge i$, a contradiction. Thus $T_{i,j} = 0$ for all $j \ge i$, and~\eqref{eq:M0_force_cap}--\eqref{eq:M0_tau_cap} force $F_{i,j} = P_{i,j} = \tau_{i,j} = 0$. The table equations for $i \ge 2$ follow directly from~\eqref{eq:M0_upper_root}.

For part~(ii), because every block $i \ge 2$ requires at least one support in layer $L_i - 1$ by~\eqref{eq:M0_support_min} and~\eqref{eq:M0_layer_link}, no empty layer can exist below the highest occupied layer; combined with $L_1 = 1$ and $L_1 \le \dots \le L_N$, this implies $L_{k+1} - L_k \in \{0, 1\}$ for all $k \in \{1, \dots, N - 1\}$ and hence $L_i \le 1 + (i - 1) = i$. By part~(i), only blocks in $\{i + 1, \dots, N\}$ can exert downward force on block~$i$; summing the vertical force balance equations~\eqref{eq:force_bal} over the subset $\{i, i + 1, \dots, N\}$ cancels all internal forces within $\{i, \dots, N\}$ and shows that the total downward force exiting block~$i$ toward $\{1, \dots, i - 1\}$ is at most the combined weight $N - i + 1$ of blocks $\{i, \dots, N\}$. Since all forces are nonnegative, each individual force $F_{i,j}$ satisfies~\eqref{eq:upper_force_cap}. The tightened layer-linking bounds follow from $1 \le L_j \le j$ and $2 \le L_i \le i$, which imply $1 - j \le L_i - L_j - 1 \le i - 1$ whenever $T_{i,j} = 0$.

For part~(iii), summing the vertical force balance equations~\eqref{eq:force_bal} over all $i \in \mathcal{N}$ cancels every internal contact force $F_{i,j}$ ($1 \le j < i \le N$) and yields $\sum_{i=1}^N F_{i,0} = N$. Because $F_{i,0} = 0$ for all $i \ge 2$ by part~(i), we obtain $F_{1,0} = N$ identically, so $\tau_{1,0} = F_{1,0} P_{1,0} = N P_{1,0}$, establishing~\eqref{eq:base_linearization}. Summing the moment balance equations~\eqref{eq:moment_bal} over all $i \in \mathcal{N}$ similarly cancels every internal torque $\tau_{i,j}$ ($1 \le j < i \le N$) and gives $\sum_{i=1}^N X_i = \tau_{1,0} = N P_{1,0}$, proving~\eqref{eq:global_torque_identity}.
\end{proof}

\subsection{Proof of Lemma~\ref{lem:support_bound}}
\label{app:support_bound_proof}

\begin{proof}[Proof of Lemma~\ref{lem:support_bound}]
For part~(i), an active contact $T_{ij} = 1$ requires the effective cores of blocks $i$ and $j$ to intersect, so $|X_i - X_j| \le 2w$. The centers of all supports of block~$i$ therefore lie in an interval of length $4w = 2(1 - \eta)$, and, lying in one layer, they are pairwise at least $1$ apart. Such an interval contains at most $\lfloor 4w \rfloor + 1$ points with pairwise distances at least $1$, that is, three points when $\eta = 0$ and two points when $\eta > 0$. The same argument bounds the number of blocks resting on block~$j$.

For part~(ii), every block of layer~$2$ rests only on block~$1$, the loads satisfy $\sum_{k \in \mathcal{B}_2} F_{k,1} = N - 1$, and every contact point satisfies $P_{k,1} \ge X_1 - w$. Rotational equilibrium of block~$1$ together with $P_{1,0} \le -\epsilon$ gives $-N\epsilon \ge N P_{1,0} = X_1 + \sum_{k \in \mathcal{B}_2} F_{k,1} P_{k,1} \ge N X_1 - (N - 1) w$.

For part~(iii), every block in a layer $L_i \ge 2$ rests on at least one block of layer $L_i - 1$. Following such supports downward produces a path of $L_i - 1$ active contacts from block~$i$ to block~$1$, and each contact changes the center coordinate by at most $2w$.
\end{proof}

\subsection{Proof of Lemma~\ref{lem:s1_translation}}
\label{app:s1_translation_proof}

\begin{proof}[Proof of Lemma~\ref{lem:s1_translation}]
Consider an optimal stack and suppose for contradiction that $P_{1,0} < -\epsilon$. Translate all block centers and all contact points, including $P_{1,0}$, rightward by $\delta := -\epsilon - P_{1,0} > 0$, leaving all layer indices, contact indicators $T_{ij}$, and forces $F_{ij}$ unchanged. Every relative coordinate $P_{ij} - X_i$, $P_{ij} - X_j$, and $X_i - X_k$ is invariant, so same-layer non-overlap and all contact-interval inclusions are preserved, and the translated table reaction equals $-\epsilon$. Vertical force balance is unaffected, and for the rotational equilibrium of any block $i$, the force-balance identity $\sum_{j < i} F_{ij} - \sum_{k > i} F_{ki} = 1$ (with $j = 0$ included for $i = 1$) gives
\[
\sum_{j < i} F_{ij} (P_{ij} + \delta) - \sum_{k > i} F_{ki} (P_{ki} + \delta)
= X_i + \delta.
\]
Hence the translated stack is feasible with overhang $z + \delta > z$, a contradiction. Therefore $P_{1,0} = -\epsilon$, $\tau_{1,0} = F_{1,0} P_{1,0} = -N\epsilon$, and~\eqref{eq:global_torque_identity} gives $\sum_{i=1}^{N} X_i = -N\epsilon$.
\end{proof}

\subsection{Proof of Theorem~\ref{thm:layer_doubling}}
\label{app:layer_doubling_proof}

\begin{proof}[Proof of Theorem~\ref{thm:layer_doubling}]
Fix any layer $\ell \in \{1, \dots, H - 1\}$, and list the block centers in $\mathcal{B}_\ell$ in strictly increasing order as $y_1 < y_2 < \dots < y_m$, where $m = n_\ell$. Partition $\{y_1, \dots, y_m\}$ into maximal contiguous clusters $C_1, \dots, C_p$ such that consecutive elements within the same cluster satisfy $y_{r+1} - y_r \le 4w$, while consecutive clusters are separated by a gap exceeding $4w$.

Consider any block in layer $\ell + 1$ with center $z$. Because contact with a block at $y_r$ requires $|z - y_r| \le 2w$, any block $z$ contacting both $y_b$ and $y_{b+1}$ would imply $y_{b+1} - y_b \le (y_{b+1} - z) + (z - y_b) \le 4w$. Consequently, all lower supports of block $z$ must belong to a single cluster $C_s$.

Take any cluster $C_s = \{y_a, \dots, y_b\}$ of cardinality $c := b - a + 1$, and let $z_1 < z_2 < \dots < z_k$ be the centers of all blocks in layer $\ell + 1$ whose supports lie in $C_s$. Because $z_1$ must contact at least one block in $C_s$, we have $z_1 \ge y_a - 2w$. Suppose for contradiction that $z_1 = y_a - 2w$. For any other block $y_r \in C_s$ with $r > a$, same-layer non-overlap implies $y_r - z_1 \ge (y_a + 1) - (y_a - 2w) = 1 + 2w > 2w$, so $z_1$ cannot contact $y_r$. Hence $y_a$ is the unique lower support of $z_1$, and the intersection of their effective cores degenerates to the single point $x_0 := z_1 + w = y_a - w$. Taking rotational equilibrium of block $z_1$ about $x_0$, the upward reaction from $y_a$ acts at $x_0$ and produces zero moment, every downward force $F_{u, z_1} \ge 0$ from layer $\ell + 2$ acts at $P_{u, z_1} \le z_1 + w = x_0$ and produces non-positive moment $F_{u, z_1}(P_{u, z_1} - x_0) \le 0$, whereas the unit self-weight of block $z_1$ acts at $z_1 = x_0 - w$ and produces moment $-w < 0$. The resulting moment balance reads
\[
0 = (z_1 - x_0) + \sum_{u \in \mathcal{B}_{\ell+2}} F_{u, z_1}(P_{u, z_1} - x_0) \le -w < 0,
\]
a contradiction. Thus strictly $z_1 > y_a - 2w$. By the symmetric argument taking moments of block $z_k$ about $y_b + w$, we obtain strictly $z_k < y_b + 2w$. Applying this strict inequality to the leftmost cluster $C_1$ (where $y_a = y_{\min}(\ell)$) and the rightmost cluster $C_p$ (where $y_b = y_{\max}(\ell)$) proves~\eqref{eq:endpoint_exclusion} in part~(i).

Next, within cluster $C_s$, adjacent centers satisfy $y_{r+1} - y_r \le 4w$, so $y_b - y_a \le 4w(c - 1)$. Using same-layer non-overlap $z_{r+1} - z_r \ge 2W = 1$ in layer $\ell + 1$ together with $2w = 1 - \eta \le 1$, we obtain
\[
k - 1 \le z_k - z_1 < (y_b + 2w) - (y_a - 2w) = (y_b - y_a) + 4w \le 4wc \le 2c.
\]
Since $k$ and $2c$ are integers, $k - 1 < 2c$ implies $k \le 2c$. Summing over all clusters $s = 1, \dots, p$ yields
\[
n_{\ell + 1} = \sum_{s=1}^{p} k_s \le \sum_{s=1}^{p} 2 c_s = 2 n_\ell,
\]
establishing the first inequality in~\eqref{eq:layer_doubling}.

For the linear expansion bound, starting from $y_{\min}(1) = y_{\max}(1) = X_1$ and applying~\eqref{eq:endpoint_exclusion} inductively gives
\[
y_{\max}(\ell) - y_{\min}(\ell) < 4w(\ell - 1) \le 2(\ell - 1) \qquad (\forall \ell \in \{2, \dots, H\}).
\]
By same-layer non-overlap in layer $\ell$, we have $n_\ell - 1 \le y_{\max}(\ell) - y_{\min}(\ell) < 2(\ell - 1)$, and integrality yields $n_\ell \le 2\ell - 2$, proving part~(ii). Finally, if $n_2 = 1$, then $y_{\min}(2) = y_{\max}(2) = X_2$; repeating the same induction from layer~$2$ yields $y_{\max}(\ell) - y_{\min}(\ell) < 4w(\ell - 2) \le 2(\ell - 2)$ for all $\ell \in \{3, \dots, H\}$, and hence $n_\ell \le 2\ell - 4$, proving part~(iii).
\end{proof}

\subsection{Proof of Theorem~\ref{thm:base_singleton}}
\label{app:base_singleton_proof}

\begin{proof}[Proof of Theorem~\ref{thm:base_singleton}]
Because $n_2 = 1$, block~$2$ is the sole block in layer~$2$ and carries the entire upper sub-stack $\{2, \dots, N\}$, so $F_{2,1} = N - 1$ and $F_{1,0} = N$. Rotational equilibrium of block~$1$ combined with the table clearance constraint $P_{1,0} \le -\epsilon$ gives
\begin{equation}
\label{eq:base_singleton_mom1}
X_1 + (N - 1) P_{2,1} = N P_{1,0} \le -N\epsilon.
\end{equation}
Since the contact point $P_{2,1}$ lies within the effective core $[X_1 - w, X_1 + w]$ of block~$1$, we have $-X_1 \le w - P_{2,1}$. Adding this inequality to~\eqref{eq:base_singleton_mom1} cancels $X_1$ and yields
\[
N P_{2,1} \le w - N\epsilon,
\qquad\text{or equivalently,}\qquad
\Delta := P_{2,1} + \epsilon \le \frac{w}{N}.
\]

We next construct an $(N - 1)$-block single-root $\eta$-stable stack from the upper blocks $\{2, \dots, N\}$ by translating them horizontally by $-\Delta$ and shifting all layer indices down by one. Specifically, set $L'_i := L_i - 1$ and $X'_i := X_i - \Delta$ for all $i \in \{2, \dots, N\}$; keep all contact indicators $T'_{ij} := T_{ij}$ and vertical forces $F'_{ij} := F_{ij}$ for $2 \le j < i \le N$, with translated contact points $P'_{ij} := P_{ij} - \Delta$; and define the new base reaction on block~$2$ by $F'_{2,0} := F_{2,1} = N - 1$ acting at $P'_{2,0} := P_{2,1} - \Delta = -\epsilon$.

We verify that this $(N - 1)$-block configuration is feasible. Block~$2$ is the unique block with $L'_2 = 1$, while $L'_i \ge 2$ for all $i \ge 3$. At the table, $P'_{2,0} = -\epsilon$ satisfies the table clearance constraint, and $P'_{2,0} - X'_2 = P_{2,1} - X_2 \in [-w, w]$, so $P'_{2,0} \in [X'_2 - w, \min(-\epsilon, X'_2 + w)]$. Uniform horizontal translation preserves all same-layer center separations $X'_i - X'_j = X_i - X_j \ge 2W$ and all internal contact core inclusions $P'_{ij} - X'_i = P_{ij} - X_i \in [-w, w]$ and $P'_{ij} - X'_j = P_{ij} - X_j \in [-w, w]$. Vertical force equilibrium holds identically for every $i \in \{2, \dots, N\}$. For rotational equilibrium of any block $i \in \{2, \dots, N\}$, using vertical force balance $\sum_{j} F'_{ij} - \sum_{k} F'_{ki} = 1$ with $j = 0$ included for $i = 2$, we obtain
\[
\sum_{j} F'_{ij} P'_{ij} - \sum_{k} F'_{ki} P'_{ki}
= \left(\sum_{j} F_{ij} P_{ij} - \sum_{k} F_{ki} P_{ki}\right)
  - \Delta \left(\sum_{j} F_{ij} - \sum_{k} F_{ki}\right)
= X_i - \Delta = X'_i.
\]
Hence the shifted configuration is a feasible single-root $(N - 1)$-block $\eta$-stable stack, and therefore $\max_{2 \le i \le N}(X'_i + W) \le z^*(N - 1, \eta)$. Reversing the shift gives
\begin{equation}
\label{eq:base_singleton_upper}
\max_{2 \le i \le N}(X_i + W) = \Delta + \max_{2 \le i \le N}(X'_i + W) \le z^*(N - 1, \eta) + \frac{w}{N}.
\end{equation}

It remains to bound the overhang of block~$1$ itself. Using $P_{2,1} \ge X_1 - w$ in~\eqref{eq:base_singleton_mom1} yields $N X_1 - (N - 1) w \le -N\epsilon$, so
\begin{equation}
\label{eq:base_singleton_b1}
X_1 + W \le W - \epsilon + \frac{N - 1}{N} w.
\end{equation}
On the other hand, consider the classical single-column harmonic stack of $N - 1$ blocks with base reaction at $-\epsilon$, base center $-\epsilon - \frac{N - 2}{N - 1} w$, and inter-block center increments $w / r$ for $r = 1, \dots, N - 2$ from top to bottom (empty for $N = 2$). Its top block achieves the feasible overhang $W - \epsilon - \frac{N - 2}{N - 1} w + w \sum_{r=1}^{N-2} 1/r = W - \epsilon + w \sum_{k=2}^{N-1} 1/k$, so
\[
z^*(N - 1, \eta) + \frac{w}{N} \ge W - \epsilon + w \sum_{k=2}^{N} \frac{1}{k}.
\]
Because $\sum_{k=2}^{N} 1/k = 1/2 = (N - 1)/N$ for $N = 2$ and $\sum_{k=2}^{N} 1/k > (N - 1)/N$ for all $N \ge 3$ (the increment from $N$ to $N + 1$ satisfies $1/(N + 1) > 1/(N(N + 1))$), and $w > 0$, inequality~\eqref{eq:base_singleton_b1} implies $X_1 + W \le z^*(N - 1, \eta) + w / N$. Combining this with~\eqref{eq:base_singleton_upper} proves~\eqref{eq:base_singleton_bound}.
\end{proof}

\subsection{Proof of Corollary~\ref{cor:s2p_elimination}}
\label{app:s2p_elimination_proof}

\begin{proof}[Proof of Corollary~\ref{cor:s2p_elimination}]
By Theorem~\ref{thm:base_singleton} and~\eqref{eq:s2p_condition}, any feasible stack with $n_2 = 1$ has overhang at most $z^*(N - 1, \eta) + w / N \le \overline{z}(N - 1, \eta) + w / N < \underline{z}(N, \eta) \le z^*(N, \eta)$, and is therefore strictly suboptimal. Since $n_2 \in \{1, 2\}$ by Theorem~\ref{thm:layer_doubling}(ii), every optimal stack must have $n_2 = 2$, equivalently $b_2 = 0$ and $L_2 = L_3 = 2$. Because each block index $j \ge 4$ can increase the layer index by at most one per block, $L_j \le L_3 + (j - 3) = j - 1$ holds for all $j \ge 3$, proving part~(i). By Theorem~\ref{thm:layer_doubling}(ii), $n_3 \le 2 n_2 = 4$, so at most $n_1 + n_2 + n_3 \le 7$ blocks occupy layers~$1$ through~$3$; hence every block $i \ge 8$ satisfies $L_i \ge 4$, whereas $L_3 = 2$, so no contact $(i, 3)$ is possible, establishing part~(ii). Finally, blocks~$2$ and~$3$ both lie in layer~$2$ and rest only on block~$1$, so $T_{2,1} = T_{3,1} = 1$, $T_{3,2} = 0$, and $F_{2,1} + F_{3,1} = N - 1$ with $F_{2,1}, F_{3,1} \ge 1$; strict endpoint exclusion~\eqref{eq:endpoint_exclusion} at layer~$1$ together with $X_2 - X_3 \ge 1$ gives $X_2 > X_1 > X_3$, proving part~(iii).
\end{proof}

\subsection{Proof of Lemma~\ref{lem:wlog_base_shift}}
\label{app:wlog_base_shift_proof}

\begin{proof}[Proof of Lemma~\ref{lem:wlog_base_shift}]
We construct $\tilde{X}$ in three steps that keep the support topology $T$ and all contact forces $F$ unchanged.

\emph{Step~1 (global translation).}
Translate all block centers and all contact points, including $P_{1,0}$, rightward by $\Delta_0 \ge 0$, leaving all forces unchanged, and denote the result by $X^{(1)}, P^{(1)}$. Every relative coordinate $P_{i,j} - X_i$, $P_{i,j} - X_j$, and $X_i - X_k$ is invariant, so force balance, moment balance (by vertical force balance, as in the proof of Lemma~\ref{lem:s1_translation}), contact-interval inclusions, and same-layer non-overlap remain satisfied, and $P_{1,0}^{(1)} = -\epsilon$. The quantity $\delta_1 = X_1^{(1)} + w - P_{2,1}^{(1)}$ is unchanged, and $\delta_1 \ge 0$ by the contact interval of $(2, 1)$.

\emph{Step~2 (leftward shift of the base block).}
Block~$1$ is the only block of layer~$1$, so every block of $\mathcal{B}_2$ rests only on block~$1$; vertical force balance gives $F_{k,1} \ge 1$ for $k \in \mathcal{B}_2$, $\sum_{k \in \mathcal{B}_2} F_{k,1} = N - 1$, and $F_{1,0} = N$. If $n_2 = 2$, then $P_{2,1}^{(1)} \ge X_2^{(1)} - w$ and $P_{3,1}^{(1)} \le X_3^{(1)} + w \le X_2^{(1)} - 1 + w$, so
\begin{equation}
\label{eq:p21_ge_p31}
P_{2,1}^{(1)} - P_{3,1}^{(1)} \ge 1 - 2w = \eta \ge 0.
\end{equation}
Hence $P_{k,1}^{(1)} \le P_{2,1}^{(1)} = X_1^{(1)} + w - \delta_1$ for every $k \in \mathcal{B}_2$ (trivially when $n_2 = 1$). Now hold blocks $2, \dots, N$, all inter-block contact points, and all forces fixed, and move only block~$1$ leftward by $\delta \in [0, \delta_1]$ to $X_1^{(2)} := X_1^{(1)} - \delta$. Layer~$1$ contains no other block, and nothing changes for blocks $2, \dots, N$. For block~$1$, vertical force balance is unchanged, and rotational equilibrium
\[
N P_{1,0}^{(2)} = X_1^{(2)} + \sum_{k \in \mathcal{B}_2} F_{k,1} P_{k,1}^{(1)} = N P_{1,0}^{(1)} - \delta
\]
moves the table reaction to $P_{1,0}^{(2)} = -\epsilon - \delta / N$. The contact conditions of block~$1$ hold for every $\delta \in [0, \delta_1]$:
\begin{enumerate}
\item[(a)] \emph{table clearance:} $P_{1,0}^{(2)} \le -\epsilon$;
\item[(b)] \emph{upper contacts, right edge:} $P_{k,1}^{(1)} \le X_1^{(1)} + w - \delta_1 \le X_1^{(2)} + w$ for $k \in \mathcal{B}_2$;
\item[(c)] \emph{upper contacts, left edge:} $P_{k,1}^{(1)} \ge X_1^{(1)} - w \ge X_1^{(2)} - w$ for $k \in \mathcal{B}_2$;
\item[(d)] \emph{table reaction, right edge:} by~(b) and $\sum_{k \in \mathcal{B}_2} F_{k,1} = N - 1$, $N P_{1,0}^{(2)} \le X_1^{(2)} + (N - 1)(X_1^{(2)} + w) = N(X_1^{(2)} + w) - w$, so $P_{1,0}^{(2)} \le X_1^{(2)} + w - w / N$;
\item[(e)] \emph{table reaction, left edge:} $P_{1,0}^{(2)} = P_{1,0}^{(1)} - \delta / N \ge X_1^{(1)} - w - \delta = X_1^{(2)} - w$.
\end{enumerate}
We take $\delta = \delta_1$, which gives $P_{2,1}^{(1)} = X_1^{(2)} + w$ and $P_{1,0}^{(2)} = -\epsilon - \delta_1 / N$.

\emph{Step~3 (final translation).}
Translate the configuration of Step~2 rightward by $\delta_1 / N$, which preserves feasibility as in Step~1. The result $\tilde{X}$ satisfies $\tilde{P}_{1,0} = -\epsilon$, $\sum_{i=1}^N \tilde{X}_i = N \tilde{P}_{1,0} = -N\epsilon$ by~\eqref{eq:global_torque_identity}, and $\tilde{P}_{2,1} = P_{2,1}^{(1)} + \delta_1 / N = \tilde{X}_1 + w$. Every block $i \ge 2$ has moved by $\Delta_0 + \delta_1 / N = \Delta$, and block~$1$ by $\Delta_0 - \delta_1 + \delta_1 / N = \Delta - \delta_1$. This proves~(i) and~(ii).

For~(iii), $\tilde{X}_1 \ge -W$ follows from $\tilde{P}_{1,0} = -\epsilon \le \tilde{X}_1 + w$. If $n_2 = 2$ and $N \ge 3$, substituting $\tilde{P}_{1,0} = -\epsilon$, $\tilde{P}_{2,1} = \tilde{X}_1 + w$, $\tilde{P}_{3,1} \ge \tilde{X}_1 - w$, $F_{3,1} = N - 1 - F_{2,1}$, and $F_{2,1} \ge 1$ into the rotational equilibrium of block~$1$ gives
\[
-N\epsilon = \tilde{X}_1 + F_{2,1}(\tilde{X}_1 + w) + (N - 1 - F_{2,1})\tilde{P}_{3,1}
\ge N\tilde{X}_1 + (2 F_{2,1} - N + 1) w
\ge N\tilde{X}_1 - (N - 3) w,
\]
that is, $\tilde{X}_1 \le \frac{N - 3}{N} w - \epsilon$.

For~(iv), every block $i \ge 2$ gains exactly $\Delta$. If block~$1$ does not attain the maximum overhang of $X$, then $X_1 < \max_{i \ge 2} X_i$ and $\tilde{X}_1 = X_1 + \Delta - \delta_1 \le X_1 + \Delta < \max_{i \ge 2} \tilde{X}_i$, so block~$1$ does not attain the maximum overhang of $\tilde{X}$ either. Since all other blocks move by the same amount, the set of maximum-overhang blocks, and hence the principal block, is the same in both stacks, $z(\tilde{X}) = z(X) + \Delta$. The spinal or non-spinal status is determined by the geometric resting relation, not by $T$, which may omit geometric contacts without force; it is preserved because blocks $2, \dots, N$ move by a common shift, so the resting relations among them are those of $X$, and every block $k$ of layer~$2$ rests on block~$1$ in both stacks: block~$1$ is the only block of layer~$1$, so in both $X$ and $\tilde{X}$ it carries block~$k$ with a force $F_{k,1} \ge 1$ through a contact point in the intersection of their effective cores, which gives $|X_k - X_1| \le 2w$ and $|\tilde{X}_k - \tilde{X}_1| \le 2w$. Hence the principal block has the same support set in both stacks. For a maximum-overhang layer $\ell^* \ge 2$, the inequalities $\tilde{X}_{s_l} \le \tilde{X}_{s_{\ell^*}}$ hold for $l \ge 2$ because they hold in $X$ and both sides move by $\Delta$, and for $l = 1$ because $\tilde{X}_1 < \tilde{X}_{s_{\ell^*}}$.
\end{proof}

\subsection{Proof of Theorem~\ref{thm:base_right_corner}}
\label{app:base_right_corner_proof}

\begin{proof}[Proof of Theorem~\ref{thm:base_right_corner}]
For $N = 2$, the harmonic placement $X_1 = -\epsilon - w / 2$, $X_2 = -\epsilon + w / 2$ with $P_{1,0} = -\epsilon$ and $P_{2,1} = X_1 + w$ is feasible and achieves overhang $X_2 + W = W - \epsilon + w / 2 = 3w / 2$, which matches the upper bound of Theorem~\ref{thm:base_singleton} ($z^*(1, \eta) + w / 2 = w + w / 2 = 3w / 2$).

Now fix $N \ge 3$ and consider any globally optimal single-root $N$-block $\eta$-stable stack with overhang $z = z^*(N, \eta)$. We proceed in two steps.

\emph{Step 1: Block~$1$ cannot be a rightmost block.}
Because every block in layer~$2$ rests solely on block~$1$ and transfers the total upper weight $\sum_{k \in \mathcal{B}_2} F_{k,1} = N - 1$ at contact points $P_{k,1} \ge X_1 - w$, rotational equilibrium of block~$1$ combined with $P_{1,0} \le -\epsilon$ gives
\[
-N\epsilon \ge N P_{1,0} = X_1 + \sum_{k \in \mathcal{B}_2} F_{k,1} P_{k,1}
\ge X_1 + (N - 1)(X_1 - w) = N X_1 - (N - 1) w,
\]
and hence
\[
X_1 + W \le W - \epsilon + \frac{N - 1}{N} w.
\]
On the other hand, the $N$-block single-column harmonic stack with $P_{1,0} = -\epsilon$ achieves the feasible overhang $W - \epsilon + w \sum_{k=2}^{N} 1/k$. As shown in the proof of Theorem~\ref{thm:base_singleton} (Appendix~\ref{app:base_singleton_proof}), $\sum_{k=2}^{N} 1/k > (N - 1)/N$ for all $N \ge 3$ and $w > 0$, whence $z^*(N, \eta) \ge W - \epsilon + w \sum_{k=2}^{N} 1/k > X_1 + W$. Therefore every block attaining the optimal overhang belongs to $\{2, \dots, N\}$.

\emph{Step 2: Base-block shift.}
Suppose for contradiction that $P_{2,1} < X_1 + w$, so that $\delta_1 = X_1 + w - P_{2,1} > 0$ in Lemma~\ref{lem:wlog_base_shift}. By Lemma~\ref{lem:s1_translation}, $P_{1,0} = -\epsilon$, so $\Delta_0 = 0$ and $\Delta = \delta_1 / N > 0$. The stack $\tilde{X}$ of Lemma~\ref{lem:wlog_base_shift} is feasible and shifts every block $i \in \{2, \dots, N\}$ rightward by $\Delta$. By Step~1, the optimal overhang $z$ is attained by some block $i^* \ge 2$, so $\tilde{X}$ attains overhang $z + \delta_1 / N > z$, contradicting optimality. Therefore $P_{2,1} = X_1 + w$.
\end{proof}

\subsection{Scope of the Base Right-Corner Reduction}
\label{app:s4_scope}

\begin{remark}[Scope Boundaries of Corner Tightness]
\label{rem:s4_limits}
Two natural extensions of Theorem~\ref{thm:base_right_corner} are deliberately excluded from our unconditional formulations:
\begin{enumerate}
\item[(i)] \emph{Base left corner $P_{3,1} = X_1 - w$ when $n_2 = 2$.} Although $P_{3,1} = X_1 - w$ holds empirically across all computed optimal solutions for $N \in \{6, \dots, 13\}$, a local perturbation proof fails once $P_{2,1} = X_1 + w$ is active: shifting block~$1$ rightward violates $P_{2,1} \le X_1 + w$ and $P_{1,0} \le -\epsilon$, while shifting block~$3$ leftward can drag bridging blocks in layer~$3$ that rest jointly on blocks~$2$ and~$3$ and alter the load balance of block~$2$. Whether $P_{3,1} = X_1 - w$ holds at every optimum remains open.
\item[(ii)] \emph{Intermediate singleton layers.} Theorem~\ref{thm:base_right_corner} does not extend to an intermediate singleton layer $\mathcal{B}_\ell = \{u\}$ ($\ell \ge 2$) supporting two blocks $v_1, v_2$ in layer $\ell + 1$, because block~$u$ rests on finite blocks in layer $\ell - 1$ rather than on the half-infinite table $(-\infty, 0]$, so an active lower left-edge contact can block leftward translation of $u$. Indeed, for $N = 7$ and $\eta = 0$ on the fixed layer profile $(1, 2, 1, 2, 1)$, the profile-optimal overhang $1.394049$ (attained by block~$2$) is realized by the balanced stack \[\begin{aligned}(X_1, \dots, X_4) &= (0.161310,\ 0.894049,\ -0.480380,\ -0.087544), \\ (X_5, X_6, X_7) &= (0.258462,\ -0.743727,\ -0.002170),\end{aligned}\] in which block~$4$ has both upper contacts strictly interior to its core ($X_4 + W - P_{5,4} = 0.156734 > 0$, with fixed-geometry minimum non-corner force $0.183847 > 0$, and $P_{6,4} - (X_4 - W) = 0.012250 > 0$) while its lower contact $(4, 2)$ is pinned at the left edge $X_2 - W = 0.394049$. Although forcing $P_{v_1, u} = X_u + w$ did not decrease the profile optimum across numerical tests on $N \in \{6, \dots, 9\}$, proving existence of a corner-tight optimum at intermediate singleton layers remains open. Accordingly, \textsf{S4} is applied exclusively at the base block~$1$.
\end{enumerate}
\end{remark}

\subsection{Proofs for the Cut-Set Inequalities}
\label{app:c2_proofs}

\begin{proof}[Proof of Proposition~\ref{prop:c2_singleton}]
Because $u$ is the unique block in layer~$\ell$, both $v_1$ and $v_2$ rest exclusively on $u$, so $T_{v_1, u} = T_{v_2, u} = 1$, $F_{v_1, u} \ge 1$, and $F_{v_2, u} \ge 1$. Summing the vertical force equilibrium equations over all blocks in $\mathcal{B}_{>\ell} := \bigcup_{t = \ell + 1}^{H} \mathcal{B}_t$ cancels all internal forces within $\mathcal{B}_{>\ell}$ and yields
\[
F_{v_1, u} + F_{v_2, u} = M_\ell,
\]
whence $F_{v_1, u} \in [1, M_\ell - 1]$ and $F_{v_2, u} = M_\ell - F_{v_1, u}$.

Rotational equilibrium of block~$u$ reads $\tau_u^{\mathrm{out}} - (\tau_{v_1, u} + \tau_{v_2, u}) = X_u$. Substituting $\tau_{v_r, u} = F_{v_r, u} P_{v_r, u}$ ($r \in \{1, 2\}$) and $F_{v_2, u} = M_\ell - F_{v_1, u}$ yields the exact algebraic identity
\begin{align*}
\tau_u^{\mathrm{out}} - X_u
&= F_{v_1, u} P_{v_1, u} + (M_\ell - F_{v_1, u}) P_{v_2, u} \\
&= M_\ell P_{v_2, u} + F_{v_1, u} (P_{v_1, u} - P_{v_2, u}),
\end{align*}
so $Q_u = F_{v_1, u} \delta_u$.

We next bound $\delta_u = P_{v_1, u} - P_{v_2, u}$. By canonical ordering in layer $\ell + 1$, same-layer non-overlap gives $X_{v_1} - X_{v_2} \ge 2W = 1$. Since $P_{v_1, u} \ge X_{v_1} - w$ and $P_{v_2, u} \le X_{v_2} + w$, we obtain the lower bound
\[
\delta_u \ge (X_{v_1} - w) - (X_{v_2} + w) = (X_{v_1} - X_{v_2}) - 2w \ge 1 - 2w = d_L.
\]
On the other hand, both $P_{v_1, u}$ and $P_{v_2, u}$ lie in the effective core $[X_u - w, X_u + w]$ of block~$u$, which gives the upper bound
\[
\delta_u \le (X_u + w) - (X_u - w) = 2w = d_U.
\]
Moreover, strict endpoint exclusion in Theorem~\ref{thm:layer_doubling}(i) applied at layer~$\ell$ gives $X_{v_1} - X_{v_2} < 4w$; combined with $X_{v_1} - X_{v_2} \ge 1$, this forces $4w > 1$, so $d_L = 1 - 2w < 2w = d_U$ whenever such a configuration is feasible.

Because $(F_{v_1, u}, \delta_u)$ belongs to the rectangle $[1, M_\ell - 1] \times [d_L, d_U]$, the four standard McCormick inequalities for the product $Q_u = F_{v_1, u} \delta_u$ yield~\eqref{eq:c2_a}--\eqref{eq:c2_d}.
\end{proof}

\begin{proof}[Proof of Proposition~\ref{prop:c2_layer}]
Let $F_{\ell, 2} := \sum_{k \in \mathcal{B}_{\ell + 1}} F_{k, u_2}$. Summing the vertical force balance equations over all blocks $b \in \mathcal{B}_{>\ell}$ cancels all internal contact forces among blocks above layer~$\ell$, leaving
\begin{equation}
\label{eq:c2l_force}
F_{\ell, 1} + F_{\ell, 2} = M_\ell,
\end{equation}
which together with non-negativity of contact forces implies $F_{\ell, 1} \in [0, M_\ell]$.

Because $n_1 = 1$, the index $s_0(\ell) \in \{1, \dots, \ell - 1\}$ is well defined and satisfies $y_{\min}(s_0(\ell)) = y_{\max}(s_0(\ell))$. Applying Theorem~\ref{thm:layer_doubling}(i) inductively from layer $s_0(\ell)$ to layer~$\ell$ yields
\[
\Delta_\ell = X_{u_1} - X_{u_2} = y_{\max}(\ell) - y_{\min}(\ell) < 4w\big(\ell - s_0(\ell)\big).
\]
Combined with the valid box bound $\Delta_\ell \le x_{u_1}^{\mathrm{ub}} - x_{u_2}^{\mathrm{lb}}$ and same-layer non-overlap $\Delta_\ell \ge 2W = 1$, we obtain $\Delta_\ell \in [1, \overline{\Delta}_\ell]$.

Next, summing the rotational equilibrium equations $\sum_{j < b} F_{bj} P_{bj} - \sum_{k > b} F_{kb} P_{kb} = X_b$ over all $b \in \mathcal{B}_{>\ell}$ cancels all internal moments within $\mathcal{B}_{>\ell}$ and yields the cut-set moment identity
\[
\sum_{b \in \mathcal{B}_{>\ell}} X_b
= \sum_{k \in \mathcal{B}_{\ell + 1}} \big(F_{k, u_1} P_{k, u_1} + F_{k, u_2} P_{k, u_2}\big).
\]
Writing $P_{k, u_r} = X_{u_r} + (P_{k, u_r} - X_{u_r})$ for $r \in \{1, 2\}$ and substituting $F_{\ell, 2} = M_\ell - F_{\ell, 1}$ from~\eqref{eq:c2l_force} gives
\begin{align}
\sum_{b \in \mathcal{B}_{>\ell}} X_b
&= F_{\ell, 1} X_{u_1} + (M_\ell - F_{\ell, 1}) X_{u_2}
   + \sum_{k \in \mathcal{B}_{\ell + 1}} \Big[ F_{k, u_1}(P_{k, u_1} - X_{u_1}) + F_{k, u_2}(P_{k, u_2} - X_{u_2}) \Big] \nonumber \\
&= M_\ell X_{u_2} + F_{\ell, 1} \Delta_\ell + R_\ell, \label{eq:c2l_exact}
\end{align}
where $R_\ell$ denotes the sum in brackets. Because every active contact satisfies $|P_{k, u_r} - X_{u_r}| \le w$ and $F_{k, u_r} \ge 0$, the triangle inequality and~\eqref{eq:c2l_force} imply
\[
|R_\ell| \le w \sum_{k \in \mathcal{B}_{\ell + 1}} (F_{k, u_1} + F_{k, u_2}) = M_\ell w.
\]
Rearranging~\eqref{eq:c2l_exact} therefore yields the bilateral enclosure
\begin{equation}
\label{eq:c2l_sandwich}
F_{\ell, 1} \Delta_\ell - M_\ell w \le W_\ell \le F_{\ell, 1} \Delta_\ell + M_\ell w.
\end{equation}

Finally, since $(F_{\ell, 1}, \Delta_\ell) \in [0, M_\ell] \times [1, \overline{\Delta}_\ell]$, the McCormick convex hull of $F_{\ell, 1} \Delta_\ell$ over this rectangle is characterized by
\[
\max\big(F_{\ell, 1},\; M_\ell \Delta_\ell + \overline{\Delta}_\ell F_{\ell, 1} - M_\ell \overline{\Delta}_\ell\big)
\le F_{\ell, 1} \Delta_\ell
\le \min\big(M_\ell \Delta_\ell + F_{\ell, 1} - M_\ell,\; \overline{\Delta}_\ell F_{\ell, 1}\big).
\]
Combining these four bounds with~\eqref{eq:c2l_sandwich} gives~\eqref{eq:c2l_a}--\eqref{eq:c2l_d}.
\end{proof}

\begin{remark}[Structure and Scope of \textsf{C2} and \textsf{C2-Layer}]
\label{rem:c2_discussion}
A downward cut-set aggregation below a two-block layer~$\ell$ requires layer $\ell - 1$ to contain two blocks and then coincides with the upward cut of layer $\ell - 1$; the only exception, $\ell = 2$ with $n_1 = 1$, is covered by Proposition~\ref{prop:c2_singleton} at the base without the $\pm M_\ell w$ residual. We apply Proposition~\ref{prop:c2_layer} at every non-top two-block layer, which gives one aggregation per such layer; other aggregations, for instance of the downward forces of two upper blocks onto a wider lower layer, are possible and are not used here. Each such layer adds only four sparse inequalities that link all $2 n_{\ell + 1}$ bilinear contact terms crossing the cut; \eqref{eq:c2l_a} and~\eqref{eq:c2l_c} use only $\Delta_\ell \ge 1$, whereas~\eqref{eq:c2l_b} and~\eqref{eq:c2l_d} are tightest when the nearest singleton layer lies directly below layer~$\ell$, where $\overline{\Delta}_\ell \le 4w$.
\end{remark}

\subsection{Proof of Theorem~\ref{thm:second_moment_energy}}
\label{app:second_moment_proof}

\begin{proof}[Proof of Theorem~\ref{thm:second_moment_energy}]
Define shifted contact points $p_{ij} := P_{ij} + \epsilon$ for $1 \le j < i \le N$ and $p_{1,0} := P_{1,0} + \epsilon = 0$. For each block $b \in \{1, \dots, N\}$, adopt the convention that the lower sum over $j < b$ includes $j = 0$ (with $F_{1,0} = N$) when $b = 1$. Vertical force equilibrium and rotational equilibrium about $-\epsilon$ read
\[
\sum_{j < b} F_{bj} - \sum_{k > b} F_{kb} = 1,
\qquad
\sum_{j < b} F_{bj} p_{bj} - \sum_{k > b} F_{kb} p_{kb} = x_b.
\]
Subtracting $x_b$ times the first equation from the second gives the center-of-mass moment balance
\begin{equation}
\label{eq:energy_local_mom}
\sum_{j < b} F_{bj} (p_{bj} - x_b) - \sum_{k > b} F_{kb} (p_{kb} - x_b) = 0.
\end{equation}
For any scalar $p$, write $p^2 = x_b^2 + 2 x_b (p - x_b) + (p - x_b)^2$, and note that $p_{rs} - x_b = P_{rs} - X_b$. Forming the weighted difference of squared shifted contact points around block~$b$ and eliminating the linear cross-term via~\eqref{eq:energy_local_mom} yields
\[
\sum_{j < b} F_{bj} p_{bj}^2 - \sum_{k > b} F_{kb} p_{kb}^2
= x_b^2 + \sum_{j < b} F_{bj} (P_{bj} - X_b)^2 - \sum_{k > b} F_{kb} (P_{kb} - X_b)^2.
\]
Summing this identity over all $b \in \{1, \dots, N\}$, every internal contact term $F_{ij} p_{ij}^2$ ($1 \le j < i \le N$) on the left-hand side appears once with a positive sign (for $b = i$) and once with a negative sign (for $b = j$) and cancels identically, leaving only the base term $F_{1,0} p_{1,0}^2 = N \cdot 0^2 = 0$. On the right-hand side, the base contact contributes $F_{1,0}(P_{1,0} - X_1)^2 = N(-\epsilon - X_1)^2 = N x_1^2$, while each internal pair $(i, j)$ contributes $F_{ij}(P_{ij} - X_i)^2 - F_{ij}(P_{ij} - X_j)^2$. Rearranging yields~\eqref{eq:second_moment_identity}.

Next, multiply the vertical force balance equation of each block $b \in \{1, \dots, N\}$ by $L_b - 1$ and sum over $b$. Because $L_1 - 1 = 0$, the base contact force $F_{1,0}$ vanishes, and because every active internal contact $(i, j)$ connects adjacent layers ($L_i - L_j = 1$), the net coefficient of $F_{ij}$ equals $(L_i - 1) - (L_j - 1) = 1$. Hence
\[
\sum_{1 \le j < i \le N} F_{ij} = \sum_{b=1}^{N} (L_b - 1) = \sum_{\ell=2}^{H} (\ell - 1) n_\ell.
\]
Since $F_{ij} \ge 0$, $F_{ij}(P_{ij} - X_i)^2 \ge 0$, and $|P_{ij} - X_j| \le w$ for every active contact, applying these bounds to~\eqref{eq:second_moment_identity} gives
\[
\sum_{b=1}^{N} x_b^2 + N x_1^2
\le \sum_{1 \le j < i \le N} F_{ij} (P_{ij} - X_j)^2
\le w^2 \sum_{1 \le j < i \le N} F_{ij}
= w^2 \sum_{\ell=2}^{H} (\ell - 1) n_\ell = E(\mathbf{n}),
\]
which establishes~\eqref{eq:energy_quad}.
\end{proof}

\subsection{Proof of Corollary~\ref{cor:energy_infeasibility}}
\label{app:energy_infeasibility_proof}

\begin{proof}[Proof of Corollary~\ref{cor:energy_infeasibility}]
Suppose a feasible single-root $\eta$-stable stack with profile $\mathbf{n}$ exists. Translating the entire stack rightward by $-\epsilon - P_{1,0} \ge 0$ produces a feasible stack with the same profile satisfying $P_{1,0} = -\epsilon$. Let $x_b := X_b + \epsilon$, and for each layer $\ell \in \{1, \dots, H\}$ let $\bar{x}_\ell := \frac{1}{n_\ell} \sum_{b \in \mathcal{B}_\ell} x_b$ denote the shifted layer mean. Because blocks in $\mathcal{B}_\ell = \{i_{\ell, 0} < \dots < i_{\ell, n_\ell - 1}\}$ are indexed from right to left and spaced by at least $2W = 1$, we have $x_{i_{\ell, r}} - x_{i_{\ell, s}} \ge s - r$ for all $0 \le r < s \le n_\ell - 1$. By Lagrange's identity,
\begin{align}
\sum_{b \in \mathcal{B}_\ell} x_b^2
&= n_\ell \bar{x}_\ell^2 + \frac{1}{n_\ell} \sum_{0 \le r < s \le n_\ell - 1} (x_{i_{\ell, r}} - x_{i_{\ell, s}})^2 \nonumber \\
&\ge n_\ell \bar{x}_\ell^2 + \frac{1}{n_\ell} \sum_{0 \le r < s \le n_\ell - 1} (s - r)^2
= n_\ell \bar{x}_\ell^2 + V_{\min}(n_\ell), \label{eq:layer_var_lb}
\end{align}
where the final equality uses the variance of $n_\ell$ consecutive integers, $\sum_{r=0}^{n_\ell - 1} \big(r - (n_\ell - 1)/2\big)^2 = n_\ell(n_\ell^2 - 1)/12$. Summing~\eqref{eq:layer_var_lb} over $\ell \in \{1, \dots, H\}$ and substituting into~\eqref{eq:energy_quad} yields
\begin{equation}
\label{eq:mean_energy_budget}
N x_1^2 + \sum_{\ell=1}^{H} n_\ell \bar{x}_\ell^2 \le \Delta E(\mathbf{n}).
\end{equation}
Since the left-hand side is non-negative, feasibility requires $\Delta E(\mathbf{n}) \ge 0$.
\end{proof}

\subsection{Proof of the Layer-Wise Energy Bound}
\label{app:energy_proof}

\begin{proof}[Proof of Corollary~\ref{cor:energy_layer_ub}]
Fix $\ell \in \{2, \dots, H\}$ and write $x_b := X_b + \epsilon$, $\bar{x} := \bar{x}_\ell$, and $m := n_\ell$. Note that $X_{i_{\ell, 0}} + W = x_{i_{\ell, 0}} + w$. By~\eqref{eq:global_torque_identity} and $P_{1,0} = -\epsilon$, $x_1 + \sum_{t \notin \{1, \ell\}} n_t \bar{x}_t = -m \bar{x}$. Applying the Cauchy--Schwarz inequality to the weighted sum $\big(\sqrt{N + 1}\, x_1\big)\frac{1}{\sqrt{N + 1}} + \sum_{t \notin \{1, \ell\}} \big(\sqrt{n_t}\, \bar{x}_t\big)\sqrt{n_t} = -m \bar{x}$ and using $\sum_{t \notin \{1, \ell\}} n_t = N - m - 1$ gives
\[
m^2 \bar{x}^2
\le \left((N + 1) x_1^2 + \sum_{t \notin \{1, \ell\}} n_t \bar{x}_t^2\right)
    \left(\frac{1}{N + 1} + N - m - 1\right),
\]
which holds for all $N \ge 2$ (including the degenerate case $N - m - 1 = 0$). Adding $m \bar{x}^2$ yields
\begin{equation}
\label{eq:energy_mean_lb}
N x_1^2 + \sum_{t=1}^{H} n_t \bar{x}_t^2 \ge A_m \bar{x}^2.
\end{equation}

If $m = 1$, then $x_{i_{\ell, 0}} = \bar{x}$, and combining~\eqref{eq:mean_energy_budget} with~\eqref{eq:energy_mean_lb} immediately yields $A_1 \bar{x}^2 \le \Delta E(\mathbf{n})$, hence $x_{i_{\ell, 0}} + w \le w + \sqrt{\Delta E(\mathbf{n}) / A_1}$, establishing the $m = 1$ case of~\eqref{eq:energy_layer_ub}.

Now suppose $m \ge 2$. Set $c_r := (m - 1)/2 - r$ and $u_r := (x_{i_{\ell, r}} - \bar{x}) - c_r$ for $r \in \{0, \dots, m - 1\}$, so that $\sum_{r=0}^{m-1} c_r = \sum_{r=0}^{m-1} u_r = 0$, $\sum_{r=0}^{m-1} c_r^2 = V_{\min}(m)$, and $x_{i_{\ell, 0}} = \bar{x} + (m - 1)/2 + u_0$. Because $x_{i_{\ell, r}} - x_{i_{\ell, r+1}} \ge 1 = c_r - c_{r+1}$, the sequence $u_0 \ge u_1 \ge \dots \ge u_{m-1}$ is non-increasing, and zero sum forces $u_0 \ge 0$. Expanding the excess within-layer variance gives
\[
\sum_{r=0}^{m-1} (x_{i_{\ell, r}} - \bar{x})^2 - V_{\min}(m)
= 2 \sum_{r=0}^{m-1} c_r u_r + \sum_{r=0}^{m-1} u_r^2.
\]
For the linear term, since both $(c_1, \dots, c_{m-1})$ and $(u_1, \dots, u_{m-1})$ are non-increasing, Chebyshev's sum inequality gives
\[
\sum_{r=1}^{m-1} c_r u_r
\ge \frac{1}{m - 1} \left(\sum_{r=1}^{m-1} c_r\right) \left(\sum_{r=1}^{m-1} u_r\right)
= \frac{1}{m - 1}\left(-\frac{m - 1}{2}\right)(-u_0) = \frac{1}{2} u_0,
\]
so $2 \sum_{r=0}^{m-1} c_r u_r = (m - 1) u_0 + 2 \sum_{r=1}^{m-1} c_r u_r \ge m u_0$. For the quadratic term, Cauchy--Schwarz gives $\sum_{r=1}^{m-1} u_r^2 \ge \frac{1}{m - 1}(\sum_{r=1}^{m-1} u_r)^2 = \frac{u_0^2}{m - 1}$, whence $\sum_{r=0}^{m-1} u_r^2 \ge \frac{m}{m - 1} u_0^2$. Combining this excess variance with~\eqref{eq:layer_var_lb} on all layers $t \neq \ell$, \eqref{eq:energy_quad}, and~\eqref{eq:energy_mean_lb} yields
\[
A_m \bar{x}^2 + m u_0 + \frac{m}{m - 1} u_0^2 \le \Delta E(\mathbf{n}), \qquad u_0 \ge 0.
\]
Completing the square via $v := u_0 + (m - 1)/2 \ge (m - 1)/2$ transforms the objective to $x_{i_{\ell, 0}} + w = w + \bar{x} + v$ subject to the elliptic constraint
\begin{equation}
\label{eq:energy_ellipse}
A_m \bar{x}^2 + \frac{m}{m - 1} v^2 \le K := \Delta E(\mathbf{n}) + \frac{m(m - 1)}{4}, \qquad v \ge \frac{m - 1}{2}.
\end{equation}
By the Cauchy--Schwarz inequality, the unconstrained maximum of $\bar{x} + v$ over the ellipse~\eqref{eq:energy_ellipse} is $\sqrt{\big(A_m^{-1} + (m - 1)/m\big) K}$, attained at
\[
v^* = \frac{m - 1}{m} \sqrt{\frac{K}{A_m^{-1} + (m - 1)/m}}.
\]
A direct algebraic check shows that $v^* > (m - 1)/2$ if and only if $2\sqrt{A_m \Delta E(\mathbf{n})} > m$, which gives the second branch of~\eqref{eq:energy_layer_ub}. When $2\sqrt{A_m \Delta E(\mathbf{n})} \le m$, the concave upper boundary $\bar{x}(v) = \sqrt{A_m^{-1}\big(K - \frac{m}{m - 1} v^2\big)}$ has slope at $v = (m - 1)/2$ satisfying $\bar{x}'\big((m - 1)/2\big) = -m / \big(2\sqrt{A_m \Delta E(\mathbf{n})}\big) \le -1$ (or $\Delta E(\mathbf{n}) = 0$, forcing $v = (m - 1)/2$ and $\bar{x} = 0$). Hence $v + \bar{x}(v)$ is non-increasing for $v \ge (m - 1)/2$ and attains its maximum at the boundary $v = (m - 1)/2$, giving $\bar{x} + v \le (m - 1)/2 + \sqrt{\Delta E(\mathbf{n}) / A_m}$ and proving the first branch of~\eqref{eq:energy_layer_ub}.
\end{proof}

\subsection{Proof of Proposition~\ref{prop:profile_tightening}}
\label{app:profile_tightening_proof}

\begin{proof}[Proof of Proposition~\ref{prop:profile_tightening}]
For part~(i), repeated application of Theorem~\ref{thm:layer_doubling}(i) from layer~$1$ (where $y_{\min}(1) = y_{\max}(1) = X_1$) to layer~$\ell$ yields
\[
X_1 - 2w(\ell - 1) < X_{i_{\ell, n_\ell - 1}} \le X_{i_{\ell, 0}} < X_1 + 2w(\ell - 1).
\]
By canonical ordering and same-layer non-overlap ($X_{i_{\ell, s}} - X_{i_{\ell, s+1}} \ge 2W = 1$), block $i = i_{\ell, r}$ satisfies $X_{i_{\ell, 0}} - X_i \ge r$ and $X_i - X_{i_{\ell, n_\ell - 1}} \ge n_\ell - 1 - r$. Combining these inequalities yields~\eqref{eq:rank_box_rel}, and substituting $x_1^{\mathrm{lb}} \le X_1 \le x_1^{\mathrm{ub}}$ yields~\eqref{eq:rank_box_abs}. When $T_{ij} = 1$, contact core inclusion implies $\max(X_i, X_j) - w \le P_{ij} \le \min(X_i, X_j) + w$, whence $P_{ij} \in [p_{ij}^{\mathrm{lb}}, p_{ij}^{\mathrm{ub}}]$; when $T_{ij} = 0$, the formulation forces $P_{ij} = \tau_{ij} = 0$, so including $0$ yields valid bounds for $P_{ij}$ and $\tau_{ij} = F_{ij} P_{ij}$.

For part~(ii), because blocks within each layer are indexed consecutively in decreasing order of horizontal centers with unit minimum separation, $X_{i_1} - X_{i_2} \ge i_2 - i_1$ and $X_{j_1} - X_{j_2} \ge j_2 - j_1$. Summing these two inequalities gives
\[
(X_{i_1} - X_{j_2}) + (X_{j_1} - X_{i_2})
= (X_{i_1} - X_{i_2}) + (X_{j_1} - X_{j_2})
\ge (i_2 - i_1) + (j_2 - j_1).
\]
If $T_{i_1, j_2} = T_{i_2, j_1} = 1$, then $|X_{i_1} - X_{j_2}| \le 2w$ and $|X_{j_1} - X_{i_2}| \le 2w$, so the left-hand side is at most $4w = 2(1 - \eta)$. When $\eta = 0$, $4w = 2 < 3$, so $(i_2 - i_1) + (j_2 - j_1) \ge 3$ is impossible. When $\eta > 0$, $4w = 2(1 - \eta) < 2$, so $(i_2 - i_1) + (j_2 - j_1) \ge 2$ is impossible. This proves~\eqref{eq:planar_noncrossing}.

For part~(iii), every block in $\mathcal{B}_{\ell+1}$ requires at least one support in $\mathcal{B}_\ell = \{u\}$, so $T_{v, u} = 1$. Summing vertical force balance over $\mathcal{B}_{>\ell}$ gives $\sum_{v \in \mathcal{B}_{\ell+1}} F_{v, u} = M_\ell$, which fixes $F_{v, u} = M_\ell$ when $n_{\ell+1} = 1$. When $n_{\ell+1} = 2$, strict endpoint exclusion~\eqref{eq:endpoint_exclusion} at layer~$\ell$ gives $X_{v_2} > X_u - 2w$ and $X_{v_1} < X_u + 2w$; combining these with $X_{v_1} - X_{v_2} \ge 1$ and $2w \le 1$ yields $X_{v_1} > X_u > X_{v_2}$, and Proposition~\ref{prop:c2_singleton} supplies~\eqref{eq:c2_facets}.
\end{proof}

\subsection{Proof of Theorem~\ref{thm:profile_decomp}}
\label{app:profile_decomp_proof}

\begin{proof}[Proof of Theorem~\ref{thm:profile_decomp}]
Let $X^*$ be any globally optimal single-root $N$-block $\eta$-stable stack with overhang $z^*(N, \eta)$ and layer profile $\mathbf{n}^* = (n_1^*, \dots, n_{H^*}^*)$. By Theorem~\ref{thm:layer_doubling}, $\mathbf{n}^* \in \mathcal{P}(N)$; by Corollary~\ref{cor:s2p_elimination} under~\eqref{eq:s2p_condition}, $n_2^* = 2$; and by Corollary~\ref{cor:energy_infeasibility}, $\Delta E(\mathbf{n}^*) \ge 0$. Hence $\mathbf{n}^* \in \mathcal{P}_2^{\mathrm{E}}(N)$. Moreover, by Lemma~\ref{lem:support_bound}, Lemma~\ref{lem:s1_translation}, Corollary~\ref{cor:s2p_elimination}, Theorem~\ref{thm:base_right_corner}, Corollary~\ref{cor:base_ub_tight}, and Propositions~\ref{prop:c2_singleton}, \ref{prop:c2_layer}, and~\ref{prop:profile_tightening}, $X^*$ is feasible for $\mathcal{F}(\mathbf{n}^*)$.

By Step~1 of the proof of Theorem~\ref{thm:base_right_corner} (Appendix~\ref{app:base_right_corner_proof}), $z^*(N, \eta)$ is attained by some block $i^*$ in a layer $\ell^* \ge 2$. Within layer $\ell^*$, canonical ordering implies $X^*_{i_{\ell^*, 0}} \ge X^*_{i^*}$, so without loss of generality $i^* = i_{\ell^*, 0}$. Combining the layer-rank upper bound~\eqref{eq:rank_box_abs} at $r = 0$ with the closed-form energy upper bound~\eqref{eq:energy_layer_ub} of Corollary~\ref{cor:energy_layer_ub} yields
\[
z^*(N, \eta) = X^*_{i_{\ell^*, 0}} + W
\le \min\!\big(x_1^{\mathrm{ub}} + W + 2w(\ell^* - 1),\; \overline{z}^{\mathrm{E}}(\mathbf{n}^*, \ell^*)\big).
\]
Every case gives $\overline{u}(\mathbf{n}) \ge \gamma$, so if $z^*(N, \eta) \le \gamma$, then $z^*(N, \eta) \le \overline{u}(\mathbf{n}^*)$. If instead $z^*(N, \eta) > \gamma$, then the minimum above strictly exceeds $\gamma$, so $i_{\ell^*, 0} \in \mathcal{C}(\mathbf{n}^*, \gamma)$. Thus $\mathcal{C}(\mathbf{n}^*, \gamma) \neq \emptyset$, ruling out case~(a), and $X^*$ is a feasible point of $\mathcal{S}(\mathbf{n}^*, \gamma)$ with objective $z^*(N, \eta) > \gamma$, so $\mathcal{S}(\mathbf{n}^*, \gamma)$ is not closed, ruling out case~(b). Hence case~(c) occurs for $\mathbf{n}^*$, and the validity of $\overline{z}(\mathbf{n}^*)$ gives $z^*(N, \eta) \le \overline{z}(\mathbf{n}^*) \le \overline{u}(\mathbf{n}^*)$. In both cases,
\[
z^*(N, \eta) \le \overline{u}(\mathbf{n}^*) \le \max_{\mathbf{n} \in \mathcal{P}_2^{\mathrm{E}}(N)} \overline{u}(\mathbf{n}),
\]
which proves~\eqref{eq:pd_global_bound}.
\end{proof}

\subsection{Proof of Proposition~\ref{prop:spinal_bound}}
\label{app:spinal_i_proof}

\begin{proof}[Proof of Proposition~\ref{prop:spinal_bound}(i)]
Fix $K \ge 2$ and $M > 0$, and introduce logarithmic coordinates $y_s := \ln u_s \in \mathbb{R}$ for $s \in \{1, \dots, K - 1\}$ with $y_K := \ln M$ fixed. Maximizing the objective in~\eqref{eq:pz_unconstrained} is equivalent to minimizing
\[
\phi_K(y; M) := \frac{1}{2} e^{-y_1} + \sum_{s=2}^{K} \left(e^{y_{s-1} - y_s} + \frac{1}{2} e^{-y_s}\right)
\]
over $y \in \mathbb{R}^{K-1}$. Because every summand is an exponential of an affine function, $\nabla^2 \phi_K(y; M) \succeq 0$. For any $v \in \mathbb{R}^{K-1}$, the quadratic form $v^\top \nabla^2 \phi_K(y; M) v$ vanishes only if $v_1 = 0$ (from $\frac{1}{2} e^{-y_1}$) and $v_{s-1} - v_s = 0$ for all $s \in \{2, \dots, K - 1\}$ (from $e^{y_{s-1} - y_s}$), forcing $v = 0$. Hence $\phi_K(\cdot; M)$ is strictly convex on $\mathbb{R}^{K-1}$. Furthermore, $\phi_K(y; M)$ is coercive: if $y_s \to -\infty$ for some $s \le K - 1$, then $\frac{1}{2} e^{-y_s} \to +\infty$; if $y_r \to +\infty$ for some $r \le K - 1$, then telescoping $\sum_{s=r+1}^{K}(y_{s-1} - y_s) = y_r - \ln M$ implies $\max_{s > r}(y_{s-1} - y_s) \ge (y_r - \ln M) / (K - r) \to +\infty$, so $\sum_{s=2}^{K} e^{y_{s-1} - y_s} \to +\infty$. Consequently, $\phi_K(\cdot; M)$ has a unique global minimizer, characterized by the first-order stationarity conditions $\partial \phi_K / \partial y_s = 0$ ($s = 1, \dots, K - 1$), which rewrite in the original variables $u_s = e^{y_s}$ directly as~\eqref{eq:pz_recurrence}.

To verify monotonicity of $u_K(u_1)$, define $\rho_s(u_1) := u_s(u_1) / (u_{s-1}(u_1) + 1/2)$ so that $u_{s+1} = u_s \rho_s$. Both $u_1$ and $\rho_1(u_1) = 2u_1$ are positive and strictly increasing in $u_1 > 0$. By induction, if $u_s$ and $\rho_s$ are positive and strictly increasing in $u_1$, then $u_{s+1} = u_s \rho_s$ and $\rho_{s+1} = \frac{u_s}{u_s + 1/2} \rho_s$ are also positive and strictly increasing in $u_1$. As $u_1 \to 0^+$, $u_2 = 2u_1^2 \to 0^+$ and inductively $u_{s+1} \le 2u_s^2 \to 0^+$. Conversely, for any $u_1 > 0$ and any $s \ge 1$,
\[
\rho_s - \rho_{s+1} = \rho_s - \frac{u_s \rho_s}{u_s + 1/2} = \frac{u_s}{(u_{s-1} + 1/2)(2u_s + 1)} < 1,
\]
so $\rho_s > 2u_1 - (s - 1)$. Thus for $u_1 \ge K$ we have $\rho_s > 1$ for all $s \le K - 1$ and $u_K(u_1) = u_1 \prod_{s=1}^{K-1} \rho_s > u_1 \to +\infty$ as $u_1 \to +\infty$. Continuity and the intermediate value theorem yield a unique root $u_1^*$ of $u_K(u_1) = M$. The equality $S_K^*(M) = \hat{S}_K(M)$ when $u_s^* - u_{s-1}^* \ge 1$ is immediate from~\eqref{eq:pz_unconstrained} and~\eqref{eq:pz_constrained}.
\end{proof}

\begin{proof}[Proof of Proposition~\ref{prop:spinal_bound}(ii)]
For part~(ii), if $\ell = 1$ (so $b_1 = 1$), rotational equilibrium of block~$1$ with $P_{1,0} \le -\epsilon$, $F_{1,0} = N$, and $P_{k,1} - X_1 \ge -w$ for all $k \in \mathcal{B}_2$ gives
\[
N(-\epsilon - X_1) \ge N(P_{1,0} - X_1) = \sum_{k \in \mathcal{B}_2} F_{k,1}(P_{k,1} - X_1) \ge -w(N - 1),
\]
whence $X_1 + W = X_1 + w + \epsilon \le w(2 - 1/N) = 2w \, S_1^*(N) = 2w \, \hat{S}_1(N)$.

Now assume $2 \le p \le \ell \le H$. For each $t \in \{p, \dots, \ell\}$, set $F_t := F_{b_t, b_{t-1}}$, $d_t := X_{b_t} - X_{b_{t-1}}$, and $r_t := \ell - t + 1$. We first prove by reverse induction on $j = \ell, \ell - 1, \dots, p$ that
\begin{equation}
\label{eq:spinal_ind}
F_j (P_{b_j, b_{j-1}} - X_{b_j}) \ge \sum_{t=j+1}^{\ell} F_t d_t - w(F_j - r_j).
\end{equation}
At the base step $j = \ell$ ($r_\ell = 1$), block $b_\ell$ is supported exclusively by $b_{\ell-1}$, so vertical force balance yields $F_\ell = 1 + \sum_{k \in \mathcal{B}_{\ell+1}} F_{k, b_\ell} \ge 1$, and moment balance of $b_\ell$ about $X_{b_\ell}$ together with $P_{k, b_\ell} - X_{b_\ell} \ge -w$ gives
\[
F_\ell (P_{b_\ell, b_{\ell-1}} - X_{b_\ell})
= \sum_{k \in \mathcal{B}_{\ell+1}} F_{k, b_\ell} (P_{k, b_\ell} - X_{b_\ell})
\ge -w(F_\ell - 1),
\]
confirming~\eqref{eq:spinal_ind} for $j = \ell$. For the inductive step $j + 1 \to j$ ($p \le j \le \ell - 1$), unique support of $b_j$ on $b_{j-1}$ gives vertical force balance
\[
F_j = 1 + F_{j+1} + \sum_{k \in \mathcal{B}_{j+1} \setminus \{b_{j+1}\}} F_{k, b_j} \ge F_{j+1} + 1.
\]
Writing the moment balance of $b_j$ about $X_{b_j}$, separating the $k = b_{j+1}$ term via $P_{b_{j+1}, b_j} - X_{b_j} = (P_{b_{j+1}, b_j} - X_{b_{j+1}}) + d_{j+1}$, and applying the induction hypothesis together with $P_{k, b_j} - X_{b_j} \ge -w$ yields
\begin{align*}
F_j (P_{b_j, b_{j-1}} - X_{b_j})
&= F_{j+1}(P_{b_{j+1}, b_j} - X_{b_{j+1}}) + F_{j+1} d_{j+1}
   + \sum_{k \in \mathcal{B}_{j+1} \setminus \{b_{j+1}\}} F_{k, b_j}(P_{k, b_j} - X_{b_j}) \\
&\ge \left[\sum_{t=j+2}^{\ell} F_t d_t - w(F_{j+1} - r_{j+1})\right]
   + F_{j+1} d_{j+1} - w(F_j - F_{j+1} - 1) \\
&= \sum_{t=j+1}^{\ell} F_t d_t - w\big(F_j - (r_{j+1} + 1)\big)
 = \sum_{t=j+1}^{\ell} F_t d_t - w(F_j - r_j),
\end{align*}
completing the induction. Because $P_{b_j, b_{j-1}} \le X_{b_{j-1}} + w$, we have $P_{b_j, b_{j-1}} - X_{b_j} \le w - d_j$, and substituting this into~\eqref{eq:spinal_ind} yields
\begin{equation}
\label{eq:spinal_prefix_mom}
S_j := \sum_{t=j}^{\ell} F_t d_t \le w(2 F_j - r_j) \qquad (\forall j \in \{p, \dots, \ell\}).
\end{equation}

When $p = 2$, block $b_1 = 1$ carries $F_{1,0} = N$ with $P_{1,0} \le -\epsilon$ and $\sum_{k \in \mathcal{B}_2 \setminus \{b_2\}} F_{k, 1} = N - 1 - F_2 \ge 0$. Applying the same moment-balance decomposition to block~$1$ using~\eqref{eq:spinal_ind} at $j = 2$ ($r_2 = \ell - 1$) gives
\begin{align*}
N(-\epsilon - X_1) &\ge N(P_{1,0} - X_1) \\
&\ge \left[\sum_{t=3}^{\ell} F_t d_t - w\big(F_2 - (\ell - 1)\big)\right] + F_2 d_2 - w(N - 1 - F_2) \\
&= \sum_{t=2}^{\ell} F_t d_t - w(N - \ell).
\end{align*}
Setting $F_1 := N \ge F_2 + 1$ and $d_1 := X_1 + W = X_1 + w + \epsilon$ (so $-\epsilon - X_1 = w - d_1$) extends~\eqref{eq:spinal_prefix_mom} to $j = 1$ with $r_1 = \ell$: $S_1 = \sum_{t=1}^{\ell} F_t d_t \le w(2N - \ell)$.

Finally, let $q \in \{p, 1\}$ (with $q = 1$ when $p = 2$ for~\eqref{eq:spinal_full_bound}, and $q = p$ for~\eqref{eq:spinal_suffix_bound}) and set $K := \ell - q + 1$. Re-index in reverse order by $u_s := F_{\ell - s + 1}$, $\tilde{d}_s := d_{\ell - s + 1}$, and $\tilde{S}_k := \sum_{s=1}^{k} u_s \tilde{d}_s = S_{\ell - k + 1} \le w(2 u_k - k)$ for $s, k \in \{1, \dots, K\}$, with $u_0 := 0$. Because $u_1 = F_\ell \ge 1$ and $u_s - u_{s-1} \ge 1$ for every $s \in \{1, \dots, K\}$, we have $0 < u_1 < \dots < u_K$ and $(u_1, \dots, u_K)$ is feasible for~\eqref{eq:pz_constrained} with total load $M = u_K$. By Abel's summation-by-parts identity and~\eqref{eq:spinal_prefix_mom},
\begin{align*}
\sum_{s=1}^{K} \tilde{d}_s
&= \frac{\tilde{S}_K}{u_K} + \sum_{k=1}^{K-1} \left(\frac{1}{u_k} - \frac{1}{u_{k+1}}\right) \tilde{S}_k \\
&\le w \frac{2 u_K - K}{u_K} + w \sum_{k=1}^{K-1} \left(\frac{1}{u_k} - \frac{1}{u_{k+1}}\right) (2 u_k - k) \\
&= 2w \sum_{s=1}^{K} \left(1 - \frac{u_{s-1}}{u_s}\right) - w \sum_{s=1}^{K} \frac{1}{u_s}
= 2w \sum_{s=1}^{K} \left(1 - \frac{u_{s-1} + 1/2}{u_s}\right) \\
&\le 2w \, S_K^*(u_K) \le 2w \, \hat{S}_K(u_K).
\end{align*}
Setting $q = p$ (where $\sum_{s=1}^{r} \tilde{d}_s = X_{b_\ell} - X_{b_{p-1}}$ and $u_r = F_p$) proves~\eqref{eq:spinal_suffix_bound}, and setting $q = 1$ (where $\sum_{s=1}^{\ell} \tilde{d}_s = X_{b_\ell} + W$ and $u_\ell = N$) proves~\eqref{eq:spinal_full_bound}.
\end{proof}

\subsection{Proof of Proposition~\ref{prop:sp_cuts}}
\label{app:sp_cuts_proof}

\begin{proof}[Proof of Proposition~\ref{prop:sp_cuts}]
Consider any feasible single-root $\eta$-stable stack with layer profile $\mathbf{n}$ satisfying $X_{s_\ell} \le \overline{x}_{s_\ell}$, and denote the integer double sum in~\eqref{eq:sp_cut} by $\Sigma_\ell \in \mathbb{Z}_{\ge 0}$. If $\Sigma_\ell = 0$, then $T_{s_t, i_{t-1, q}} = 0$ for every $t \in \{3, \dots, \ell\}$ and every $q \in \{1, \dots, n_{t-1} - 1\}$ (vacuously when $n_{t-1} = 1$). Because $n_1 = 1$, block $s_2$ is supported exclusively by $s_1 = 1$; because every block in layer $t \ge 3$ requires at least one support in layer $t - 1$, every block $s_t$ ($t = 3, \dots, \ell$) is supported exclusively by $s_{t-1}$. Hence $(s_1, s_2, \dots, s_\ell)$ is a base-rooted unique-support chain of length $\ell$, and~\eqref{eq:spinal_full_bound} of Proposition~\ref{prop:spinal_bound}(ii) yields $X_{s_\ell} + W \le 2w \, \hat{S}_\ell(N)$, confirming~\eqref{eq:sp_cut}.

If instead $\Sigma_\ell \ge 1$, then $M_\ell \ge 0$ and~\eqref{eq:sp_bigm} imply that the right-hand side of~\eqref{eq:sp_cut} is at least
\[
2w \, \hat{S}_\ell(N) + M_\ell \ge \overline{x}_{s_\ell} + W \ge X_{s_\ell} + W,
\]
so~\eqref{eq:sp_cut} holds in this case as well.
\end{proof}

\subsection{Relation to the Model of Paterson and Zwick}
\label{app:pz_relation}

For $\eta = 0$, Paterson and Zwick~\cite[Section~2]{paterson2009overhang} let a block rest on another when $|X_i - X_j| < 1$; our relation of Section~\ref{subsec:spinal_bound}, with $2w = 1$, adds the point contacts $|X_i - X_j| = 1$, through which no force can act in their model. For every stack balanced in their sense, their support set is contained in ours, because their relation is contained in ours. Moreover, every block of their support set above layer~$1$ carries a positive load, at least its own weight, so it rests in their sense on some block of the layer below; hence their support set meets every layer $1, \dots, \ell^*$. Consequently, if our support set has exactly one block in each of these layers, so has theirs: a stack balanced in their sense that is spinal under our relation is spinal in their sense. The converse can fail when a block of their spine touches a second block of the layer below at a single point. Every stack balanced in their sense is feasible in our model, because a nonnegative force distribution over a contact interval has its resultant in that interval. Hence excluding non-spinal stacks above a threshold under our relation excludes them under the definition of~\cite{paterson2009overhang} as well.

For Theorems~\ref{thm:exact_spinal_main} and~\ref{thm:main_certified_optimality}, which are stated for our model, this has the following consequences. The rational stacks and the stacks of Appendix~\ref{app:ExactSpinal} carry forces only through contacts of positive length and are spinal in their sense, so Theorem~\ref{thm:exact_spinal_main} holds in their model too, and its values $S(N)$ refer to their class of spinal stacks, which includes stacks with several blocks on the table. Every single-root stack balanced in the sense of Paterson and Zwick is feasible in our model, so the upper bounds of parts~(a) and~(c) of Theorem~\ref{thm:main_certified_optimality} bound their single-root optimum as well, and part~(b) implies that every single-root stack balanced in their sense whose overhang exceeds $\gamma_{\mathrm{ns}}(N)$ is spinal in their sense. Apart from the spinal bounds $S(N)$ of parts~(ii) to~(iv) of Theorem~\ref{thm:exact_spinal_main}, stacks with several blocks on the table are not covered.

\subsection{Proof of Theorem~\ref{thm:spinal_cert_decomp}}
\label{app:spinal_cert_proof}

\begin{proof}[Proof of Theorem~\ref{thm:spinal_cert_decomp}]
Let $X^*$ be any globally optimal single-root $N$-block $\eta$-stable stack with overhang $z^*(N, \eta)$ and layer profile $\mathbf{n}^*$. As shown in the proof of Theorem~\ref{thm:profile_decomp}, $\mathbf{n}^* \in \mathcal{P}_2^{\mathrm{E}}(N)$ and $X^* \in \mathcal{F}(\mathbf{n}^*)$. By Proposition~\ref{prop:sp_cuts}, every cut~\eqref{eq:sp_cut} holds for $X^*$, so $X^* \in \mathcal{F}_{\mathrm{SP}}(\mathbf{n}^*)$.

For part~(i), if $X^*$ is spinal, then, because every pair linked by its contact representation is a resting pair, some block $b^*$ in layer $\ell^* \in \{1, \dots, N\}$ attaining $z^*(N, \eta)$ is supported by a base-rooted unique-support chain of length $\ell^*$, and~\eqref{eq:spinal_full_bound} yields $z^*(N, \eta) = X^*_{b^*} + W \le 2w \, \hat{S}_{\ell^*}(N) \le 2w \max_{1 \le \ell \le N} \hat{S}_\ell(N)$. If $X^*$ is non-spinal, then by definition of $\overline{u}_{\mathrm{ns}}(\mathbf{n}^*)$ we have $z^*(N, \eta) \le \overline{u}_{\mathrm{ns}}(\mathbf{n}^*) \le \max_{\mathbf{n} \in \mathcal{P}_2^{\mathrm{E}}(N)} \overline{u}_{\mathrm{ns}}(\mathbf{n})$. Moreover, every non-spinal stack $X \in \mathcal{F}(\mathbf{n})$ belongs to $\mathcal{F}_{\mathrm{SP}}(\mathbf{n})$, and its principal block lies in a layer $\ell^* \ge 2$ (Section~\ref{subsec:spinal_bound}, after Lemma~\ref{lem:maximal_contacts}). If $z(X) > \gamma$, then $i_{\ell^*, 0} \in \mathcal{C}(\mathbf{n}, \gamma)$ as in the proof of Theorem~\ref{thm:profile_decomp}, so $\mathcal{S}_{\mathrm{SP}}(\mathbf{n}, \gamma)$ is not closed, case~(c) applies, and $z(X) \le \overline{z}(\mathbf{n}) \le \overline{u}_{\mathrm{SP}}(\mathbf{n})$; if $z(X) \le \gamma$, then $z(X) \le \gamma \le \overline{u}_{\mathrm{SP}}(\mathbf{n})$ in every case. Hence the bound $\overline{u}_{\mathrm{SP}}(\mathbf{n})$ defined by cases~(a) to~(c) of Theorem~\ref{thm:profile_decomp} is a valid choice for $\overline{u}_{\mathrm{ns}}(\mathbf{n})$. Combining the two cases establishes~\eqref{eq:spinal_decomp_bound}.

For part~(ii), because $\mathcal{F}_{\mathrm{SP}}(\mathbf{n}^*)$ still contains $X^*$, whether $X^*$ is spinal or not, repeating the proof of Theorem~\ref{thm:profile_decomp} with $\mathcal{F}_{\mathrm{SP}}(\mathbf{n})$ in place of $\mathcal{F}(\mathbf{n})$ directly yields $z^*(N, \eta) \le \max_{\mathbf{n} \in \mathcal{P}_2^{\mathrm{E}}(N)} \overline{u}_{\mathrm{SP}}(\mathbf{n}) \le \gamma$. When $\eta = 0$ and $\gamma = S^*(N) + \delta > 2w \, \hat{S}_\ell(N)$ for all $\ell \le N$ (whenever $S^*(N) = \hat{S}(N)$), any feasible solution to $\mathcal{S}_{\mathrm{SP}}(\mathbf{n}, \gamma)$ attaining $X_{s_\ell} + W > \gamma$ violates $X_{s_\ell} + W \le 2w \, \hat{S}_\ell(N)$ and must therefore satisfy $\Sigma_\ell \ge 1$ in~\eqref{eq:sp_cut}, preventing the solver from branching on right-chain spinal configurations.
\end{proof}

\subsection{Proof of Lemma~\ref{lem:maximal_contacts}}
\label{app:maximal_contacts_proof}

\begin{proof}[Proof of Lemma~\ref{lem:maximal_contacts}]
Every pair linked by $T$ has $L_i = L_j + 1$ and intersecting effective cores, so every pair linked by $\hat{T}$ has these two properties, and $\hat{T}$ links exactly the resting pairs. The positions $X$ and all forces and moments are unchanged, and every new contact contributes zero force and zero moment, so the equilibrium equations, the coordinate and rank boxes, the dominance and reachability inequalities, \textsf{S1}, \textsf{S4}, and the cut-set families \textsf{C2} and \textsf{C2-Layer}, which involve only positions and the forces, contact points, and moments of contacts that carry force, remain satisfied. For a new pair, $|\hat{P}_{i,j} - X_i| = |\hat{P}_{i,j} - X_j| = |X_i - X_j|/2 \le w$, so $\hat{P}_{i,j}$ lies in both effective cores and in the contact-point box of the pair; together with $L_i - L_j = 1$ and $\hat{F}_{i,j} = \hat{\tau}_{i,j} = 0$, this satisfies the contact-interval, layer-linking, and force, point, and moment capacity rows of the pair. The remaining rows that involve $T$ are of three kinds. The upper support-degree bounds of Lemma~\ref{lem:support_bound} and the non-crossing rows of Proposition~\ref{prop:profile_tightening} are proved from the two properties above and the unit spacing within layers alone, so they hold for every contact representation, in particular for $\hat{T}$. The table-contact fixings $T_{1,0} = 1$ and $T_{i,0} = 0$ for $i \ge 2$ and the fixings $T_{i,j} = 0$ between non-adjacent layers are unaffected, and the fixings $T_{i,j} = 1$ above singleton layers and in~\textsf{S2+} remain true. Finally, the lower support bounds $\sum_{j < i} T_{i,j} \ge 1$ and the right-hand side of each spinal cut~\eqref{eq:sp_cut} are nondecreasing in~$T$.
\end{proof}

\subsection{Proof of Corollary~\ref{cor:nonspinal_exclusion}}
\label{app:nonspinal_exclusion_proof}

\begin{proof}[Proof of Corollary~\ref{cor:nonspinal_exclusion}]
Suppose for contradiction that there exists a feasible non-spinal single-root $N$-block $\eta$-stable stack $X$ with overhang $z(X) > \gamma_{\mathrm{ns}}$ and layer profile $\mathbf{n} = (n_1, \dots, n_H)$. Let $\ell^* := \min\{l \in \{1, \dots, H\} : X_{s_l} + W = z(X)\}$ be the layer of its principal block $s_{\ell^*}$. Because $X$ is non-spinal, block~$1$ cannot attain $z(X)$, since otherwise $s_1 = 1$ would be the principal block with support set $\{1\}$ and $X$ would be spinal. Hence $\ell^* \ge 2$ and every block attaining $z(X)$ lies in layers $2, \dots, H$.

We first prove that $n_2 = 2$. By layer-width doubling (Theorem~\ref{thm:layer_doubling}), $n_2 \le 2 n_1 = 2$, so it suffices to exclude $n_2 = 1$. Suppose $n_2 = 1$, so layer~$2$ consists solely of block~$2$. Rotational equilibrium of block~$1$ with $F_{1,0} = N$, $F_{2,1} = N - 1$, $P_{1,0} \le -\epsilon$, and $X_1 \ge P_{2,1} - w$ yields
\[
-N\epsilon \ge N P_{1,0} = X_1 + (N - 1) P_{2,1} \ge N P_{2,1} - w,
\]
whence $\Delta' := P_{2,1} + \epsilon \le w / N$. Deleting block~$1$ and translating blocks $2, \dots, N$ leftward by $\Delta'$ yields a feasible single-root $(N - 1)$-block $\eta$-stable stack $X'$ rooted at block~$2$ (with table contact $P'_{2,0} = P_{2,1} - \Delta' = -\epsilon$) whose overhang satisfies $z(X') = z(X) - \Delta' \ge z(X) - w / N$. Because block~$1$ did not attain $z(X)$ and all blocks $2, \dots, N$ shift uniformly by $-\Delta'$ while dropping one layer, $X'$ has the exact same set of maximum-overhang blocks and the exact same principal block $s_{\ell^*}$ as $X$. Moreover, $X'$ must be non-spinal: if the support set of $s_{\ell^*}$ in $X'$ contained only one block per layer (ending at block~$2$), then prepending block~$1$ (which uniquely supports block~$2$) would show that the support set of $s_{\ell^*}$ in $X$ also contains one block per layer, contradicting that $X$ is non-spinal. Thus $z(X) \le z(X') + w / N \le z_{\mathrm{ns}}^*(N - 1, \eta) + w / N \le \overline{z}_{\mathrm{ns}}(N - 1, \eta) + w / N \le \gamma_{\mathrm{ns}}$ by~\eqref{eq:s2p_nonspinal_cond}, contradicting $z(X) > \gamma_{\mathrm{ns}}$. Hence $n_2 = 2$ and $\mathbf{n} \in \mathcal{P}_2(N)$.

Next, apply Lemma~\ref{lem:wlog_base_shift} to $X$. Because $n_1 = 1$, $n_2 = 2$, and block~$1$ does not attain $z(X)$, we obtain a feasible non-spinal single-root $N$-block $\eta$-stable stack $\tilde{X}$ with the same profile $\mathbf{n} \in \mathcal{P}_2(N)$, the same support topology $T$, the same principal block $s_{\ell^*}$ ($\ell^* \ge 2$), and overhang $z(\tilde{X}) = z(X) + \Delta \ge z(X) > \gamma_{\mathrm{ns}}$. By Lemma~\ref{lem:wlog_base_shift}(ii)--(iii), $\tilde{X}$ satisfies \textsf{S1}, \textsf{S4}, Corollary~\ref{cor:base_ub_tight}, and the bound $\tilde{X}_1 \ge -W$ (and hence the absolute boxes~\eqref{eq:rank_box_abs}), as well as every unconditional constraint of the profile subproblem, including the spinal cuts~\eqref{eq:sp_cut}, which Proposition~\ref{prop:sp_cuts} makes valid at every layer, so the energy layer bounds that select the layers receiving them play no role here. By Corollary~\ref{cor:energy_infeasibility}, $\mathbf{n} \in \mathcal{P}_2^{\mathrm{E}}(N)$, and $\ell^* \ge 2$ indexes one of its branches. Furthermore, by definition of the principal layer $\ell^*$, the dominance inequalities $\tilde{X}_{s_l} \le \tilde{X}_{s_{\ell^*}}$ hold for all $l \in \{1, \dots, H\}$; because they are weak, they remain valid when higher layers tie with $\ell^*$ for the maximum overhang. Following active supports downward from $s_{\ell^*}$ to any layer $l < \ell^*$ changes the center coordinate by at most $2w(\ell^* - l)$, and $s_l$ is the rightmost block of layer~$l$, so $\tilde{X}_{s_l} \ge \tilde{X}_{s_{\ell^*}} - 2w(\ell^* - l) > \gamma_{\mathrm{ns}} - W - 2w(\ell^* - l)$, which confirms the downward reachability bounds. Finally, replace the contact representation of $\tilde{X}$ by its maximal contact representation; by Lemma~\ref{lem:maximal_contacts}, all constraints verified above remain satisfied. Because Lemma~\ref{lem:wlog_base_shift} moves blocks $2, \dots, N$ by a common shift, the resting relations among these blocks, and hence the support set of $s_{\ell^*}$ above layer~$1$, are those of $X$, so $\tilde{X}$ is non-spinal, and by the implication following Lemma~\ref{lem:maximal_contacts} in Section~\ref{subsec:spinal_bound} the double sum in~\eqref{eq:ns_branch_constraints} is at least~$1$ (and $J_{\ell^*} \neq \emptyset$). Therefore $\tilde{X}$ is a feasible solution to branch $\mathcal{S}_{\mathrm{ns}}(\mathbf{n}, \ell^*, \gamma_{\mathrm{ns}})$ with objective $z(\tilde{X}) > \gamma_{\mathrm{ns}}$, contradicting the closure of every branch at threshold $\gamma_{\mathrm{ns}}$. This proves part~(i), and part~(ii) follows immediately since $\underline{z}_{\mathbb{Q}}(N) > \gamma_{\mathrm{ns}}$.
\end{proof}

\subsection{Proof of Proposition~\ref{prop:t0_split}}
\label{app:t0_split_proof}

\begin{proof}[Proof of Proposition~\ref{prop:t0_split}]
Suppose that a feasible non-spinal stack $X$ has profile $\mathbf{n}$, principal layer $\ell^*$, and overhang $z(X) > \gamma_{\mathrm{ns}}$. As shown in the proof of Corollary~\ref{cor:nonspinal_exclusion}, its normalized image $\tilde{X}$ from Lemma~\ref{lem:wlog_base_shift}, equipped with the maximal contact representation $\hat{T}$ of Lemma~\ref{lem:maximal_contacts}, is feasible for $\mathcal{S}_{\mathrm{ns}}(\mathbf{n}, \ell^*, \gamma_{\mathrm{ns}})$, has principal layer $\ell^*$, and has objective value $z(\tilde{X}) > \gamma_{\mathrm{ns}}$. The maximal contact representation is formed for the branch, before the constraints~\eqref{eq:t0_partition} and~\eqref{eq:t0_cuts} of any sub-branch are imposed. Because $(\tilde{X}, \hat{T})$ satisfies the non-spinal support constraint in~\eqref{eq:ns_branch_constraints}, the lowest impure layer $t_0 := \min\{t \in J_{\ell^*} : \sum_{q \ge 1} \hat{T}_{s_t, i_{t-1, q}} \ge 1\}$ exists, and~\eqref{eq:t0_partition} holds for $\hat{T}$. The block $s_{t_0}$ rests on some block $i_{t_0 - 1, q}$ with $q \ge 1$, which satisfies $\tilde{X}_{i_{t_0 - 1, q}} \le \tilde{X}_{s_{t_0 - 1}} - 1$, so contact-core overlap gives $\tilde{X}_{s_{t_0}} \le \tilde{X}_{i_{t_0 - 1, q}} + 2w \le \tilde{X}_{s_{t_0 - 1}} - 1 + 2w$. Moreover, by~\eqref{eq:t0_partition} the chain $s_1, \dots, s_{t_0 - 1}$ is a base-rooted unique-support chain, so Proposition~\ref{prop:spinal_bound}(ii) gives $\tilde{X}_{s_{t_0 - 1}} + W \le 2w \, \hat{S}_{t_0 - 1}(N)$. Hence~\eqref{eq:t0_cuts} holds. Finally, $t_0 \neq \ell^*$: otherwise the first inequality of~\eqref{eq:t0_cuts} and $1 - 2w = \eta \ge 0$ would give $\tilde{X}_{s_{\ell^*}} \le \tilde{X}_{s_{\ell^* - 1}}$, whereas the principal layer $\ell^*$ is the lowest layer attaining the maximum overhang, so $\tilde{X}_{s_{\ell^* - 1}} < \tilde{X}_{s_{\ell^*}}$. Thus $\tilde{X}$ is feasible for a sub-branch with $t_0 \in J_{\ell^*} \setminus \{\ell^*\}$ and objective value exceeding $\gamma_{\mathrm{ns}}$, contradicting its closure.
\end{proof}

\subsection{Solver Settings of the Reported Runs}
\label{app:run_settings}

Table~\ref{tab:run_settings} lists the settings of every run reported in Section~\ref{sec:Computational}; Section~\ref{subsec:TwoStagePipeline} explains the common settings and the measure of effort.

\begin{plaintable}[h]
\centering
\caption{Solver settings of the reported runs (Gurobi~13.0.1, nonconvex handling on, primal heuristics off, integer feasibility tolerance $10^{-9}$, feasibility and optimality tolerances $10^{-6}$, unless stated otherwise). The runs are those of Tables~\ref{tab:main_certs}, \ref{tab:method_progression}, \ref{tab:sp_ablation}, and~\ref{tab:n16_split} and of the tolerance checks of Section~\ref{subsec:unconditional_certs}. Thr.: threads per subproblem; Par.: subproblems solved in parallel; Limit: time limit per subproblem in seconds; gap: relative gap parameter; int.\ tol.: integer feasibility tolerance. Thread-seconds add, over all subproblems or branches, the run time multiplied by the threads of the subproblem; a run stopped by its time limit counts with its full time. The gap parameter never decides a closure, because only the statuses cutoff and infeasible, and the analytic closures, are accepted.}
\label{tab:run_settings}
\footnotesize
\setlength{\tabcolsep}{2.5pt}
\begin{tabular*}{\textwidth}{@{\extracolsep{\fill}}lccccl@{}}
\toprule
Run & Thr. & Par. & Limit & Seed & Other settings \\
\midrule
Certificates, $N \le 9$ & $1$ & $1$ & $3{,}600$ & $42$ & gap $5 \times 10^{-4}$ \\
Certificates, $N = 10, \dots, 13$ & $1$ & $3$ & $86{,}400$ & $42$ & gap $5 \times 10^{-4}$ \\
Certificates, $N = 14, 15, 16$ & $1$ & $7$ & $259{,}200$ & $42$ & gap $5 \times 10^{-4}$ \\
Non-spinal exclusion & $1$ & ${\le}\,7$ & $7{,}200$ & default & gap $10^{-4}$ \\
Sub-branches, $N = 16$ & $1$ & $11$ & $86{,}400$ & default & gap $10^{-4}$ \\
No \textsf{SP}, matched & $1$ & $7$ & $3{,}600$$^{c}$ & $42$ & gap $5 \times 10^{-4}$ \\
No \textsf{SP}, progression & $1$ & $4$ & $14{,}400$$^{c}$ & $42$ & gap $5 \times 10^{-4}$, int.\ tol.\ $10^{-5}$, no energy test$^{a}$ \\
Monolithic $\mathcal{M}_0$ & $4$ & $1$ & $259{,}200$$^{b}$ & $42$ & gap $10^{-6}$, heuristics on, int.\ tol.\ $10^{-5}$ \\
Monolithic $\mathcal{M}_0$, $N = 4$ & $1$ & $1$ & $600$ & $42$ & gap $10^{-9}$, heuristics on, int.\ tol.\ $10^{-5}$ \\
Tolerance checks & $1$ & $3$$^{d}$ & $172{,}800$ & as above & tolerances $10^{-9}$, numeric focus $3$ \\
\bottomrule
\end{tabular*}

\smallskip
\parbox{\textwidth}{\footnotesize $^{a}$The profile $(1, 2, 2, 2, 3, 2, 1)$ at $N = 13$ was closed by a separate run with three threads and seed~$1$, and independently by a split of its candidate principal layers. $^{b}$At $N = 13$, the state after $48$~h is reported. $^{c}$At $N = 13$: $14{,}400$ (matched) and $259{,}200$ (progression). $^{d}$$11$ for the sub-branches at $N = 16$; the gap is that of the runs checked.}
\end{plaintable}

\subsection{Sub-Branch Closure for $N = 16$}
\label{app:n16_split}

Table~\ref{tab:n16_split} lists the $14$ sub-branches of Proposition~\ref{prop:t0_split} for the three branches at $N = 16$ that reached the per-branch time limit of $7{,}200$~s in the main run. Every sub-branch was solved on one thread, so its run time equals its thread-seconds.

\begin{plaintable}[!htbp]
\centering
\small
\caption{Sub-branches of the three time-limited non-spinal branches at $N = 16$ ($\gamma_{\mathrm{ns}} = 2.1224044$). The solver used the double-precision value $2.12240436\ldots$ of $0.99\, S^*(16)$, which does not exceed $\gamma_{\mathrm{ns}}$. The sub-branches with $t_0 = \ell^*$ need no computation by Proposition~\ref{prop:t0_split}.}
\label{tab:n16_split}
\begin{tabular}{@{}llllrr@{}}
\toprule
Profile $\mathbf{n}$ & $\ell^*$ & $t_0$ & Status & Nodes & Time (s) \\
\midrule
$(1, 2, 2, 2, 2, 2, 3, 2)$ & $7$ & $7$ & no computation ($t_0 = \ell^*$) & & \\
 &  & $6$ & cutoff & 67{,}506 & 24.8 \\
 &  & $5$ & cutoff & 576{,}382 & 187.2 \\
 &  & $4$ & cutoff & 6{,}878{,}427 & 1{,}841.8 \\
 &  & $3$ & cutoff & 18{,}493{,}009 & 4{,}778.7 \\
\midrule
$(1, 2, 2, 2, 2, 3, 3, 1)$ & $7$ & $7$ & no computation ($t_0 = \ell^*$) & & \\
 &  & $6$ & cutoff & 38 & 0.1 \\
 &  & $5$ & cutoff & 816{,}054 & 247.4 \\
 &  & $4$ & cutoff & 7{,}617{,}894 & 2{,}179.8 \\
 &  & $3$ & cutoff & 12{,}835{,}344 & 3{,}809.5 \\
\midrule
$(1, 2, 2, 2, 3, 3, 2, 1)$ & $6$ & $6$ & no computation ($t_0 = \ell^*$) & & \\
 &  & $5$ & cutoff & 1{,}649 & 2.7 \\
 &  & $4$ & cutoff & 249{,}379 & 123.4 \\
 &  & $3$ & cutoff & 10{,}765{,}003 & 3{,}048.1 \\
\midrule
Total & & & 11 closed by cutoff & 58{,}300{,}685 & 16{,}244 \\
\bottomrule
\end{tabular}
\end{plaintable}

\section{Exact Spinal Optima}
\label{app:ExactSpinal}

This appendix proves the statements on the maximum overhang $S(N)$ of an $N$-block spinal stack in Theorem~\ref{thm:exact_spinal_main}. We keep the notation of Section~\ref{subsec:RationalConstructor}: $\eta = 0$, the spine $B_1, \dots, B_k$ is indexed from top to bottom, $B_1$ is the principal block, $t_r$ is the load that $B_r$ transmits to the block below it or to the table ($t_0 := 0$), and $\omega_r = t_r - t_{r-1} - 1$. The upper bounds below use only force and moment balance with forces in closed contact intervals, so they hold both for spinal stacks in the sense of Paterson and Zwick~\cite{paterson2009overhang} and for spinal single-root stacks of our model. The stacks constructed below are balanced in both models, because every contact that carries a positive force has an overlap of positive length.

\subsection{Algebraic Loads and the Optimal Height}
\label{subsec:AlgebraicLoads}

For $K \ge 2$, let $t_0 := 0$ and $t_1 := x$, and generate $t_2, \dots, t_K \in \mathbb{Q}(x)$ by the recurrence~\eqref{eq:pz_recurrence}. By Proposition~\ref{prop:spinal_bound}(i), the equation $t_K(x) = N$ has a unique positive root $\xi_K(N)$, and the maximizer of~\eqref{eq:pz_unconstrained} with $M = N$ consists of the loads $t_s(\xi_K(N))$. All loads are therefore rational functions of one algebraic number. For every $N \in \{4, \dots, 19\}$ and every $K \in \{2, \dots, N\}$, we isolate $\xi_K(N)$ in a short interval with rational endpoints at which $t_K - N$ has rigorously determined opposite signs, and we enclose $\hat{S}_K(N)$ by interval arithmetic with outward rounding; $\hat{S}_1(N) = 1 - 1/(2N)$ is exact. For the height $k$ of Table~\ref{tab:exact_heights}, the lower end of the enclosure of $\hat{S}_k(N)$ exceeds the upper end of every other enclosure, and all $\omega_r$ are positive at the maximizer. Hence $S_k^*(N) = \hat{S}_k(N)$ by Proposition~\ref{prop:spinal_bound}(i), and since $S_K^*(N) \le \hat{S}_K(N)$ for every $K$,
\begin{equation}
\label{eq:exact_sstar}
S^*(N) = \hat{S}_k(N) = \sum_{r=1}^{k} \left(1 - \frac{t_{r-1} + 1/2}{t_r}\right), \qquad t_r = t_r(\xi_k(N)) \ (r < k), \quad t_k = N.
\end{equation}
At these heights, the numerator of $t_k(x) - N$ is irreducible over $\mathbb{Q}$ of degree $2^{k-1}$, which is therefore the degree of $\xi_k(N)$.

\begin{plaintable}[t]
\centering
\caption{Optimal loaded-spine heights at $\eta = 0$. The column ``next'' gives the height with the second largest value $\hat{S}_K(N)$, and ``separation'' a lower bound, truncated downward, on $\hat{S}_k(N) - \max_{K \ne k} \hat{S}_K(N)$. The values $S^*(N)$ are truncated downward at the twelfth decimal. At $N \in \{5, 7\}$, the spinal optimum $S(N)$ is smaller (Proposition~\ref{prop:exception_sizes}).}
\label{tab:exact_heights}
\footnotesize
\setlength{\tabcolsep}{4pt}
\resizebox{\textwidth}{!}{%
\begin{tabular}{@{}rcccccc@{\qquad}rcccccc@{}}
\toprule
$N$ & $k$ & degree & next & separation & $S^*(N)$ & & $N$ & $k$ & degree & next & separation & $S^*(N)$ \\
\midrule
4 & 2 & 2 & 3 & 0.01208 & 1.167893218813 & & 12 & 5 & 16 & 6 & 0.01547 & 1.925091375262 \\
5 & 3 & 4 & 2 & 0.04961 & 1.317161792196 & & 13 & 5 & 16 & 6 & 0.00783 & 1.984508903252 \\
6 & 3 & 4 & 4 & 0.02647 & 1.436700214371 & & 14 & 5 & 16 & 6 & 0.00078 & 2.038216944562 \\
7 & 4 & 8 & 3 & 0.00182 & 1.531870791978 & & 15 & 6 & 32 & 5 & 0.00576 & 2.092896303215 \\
8 & 4 & 8 & 3 & 0.02582 & 1.631506820608 & & 16 & 6 & 32 & 7 & 0.01079 & 2.143842796630 \\
9 & 4 & 8 & 5 & 0.01666 & 1.715272815527 & & 17 & 6 & 32 & 7 & 0.00633 & 2.190898459596 \\
10 & 4 & 8 & 5 & 0.00312 & 1.787125861965 & & 18 & 6 & 32 & 7 & 0.00212 & 2.234570823586 \\
11 & 5 & 16 & 4 & 0.00902 & 1.858779550311 & & 19 & 7 & 64 & 6 & 0.00185 & 2.277130167643 \\
\bottomrule
\end{tabular}}
\end{plaintable}

\subsection{Realization by Ordinary Blocks}
\label{subsec:AlgebraicRealization}

For $N \in \{4, 6, 8, \dots, 19\}$, we take the shielded single-tower construction of Section~\ref{subsec:RationalConstructor} with the height $k$ of Table~\ref{tab:exact_heights} and $v_1 = N - 2k$, replace its rational loads by the exact loads of~\eqref{eq:exact_sstar}, and keep the shield seeds $y_1, \dots, y_{k-1}$ and, when $v_1 \ge 1$, the bridge seed $z_1$ of Table~\ref{tab:rational_seeds} for $N \le 16$. For $N = 17, 18, 19$, we use the seeds $(y_1, \dots, y_{k-1}; z_1)$ given by $(-\frac{2}{5}, -\frac{4}{5}, -\frac{6}{5}, -\frac{13}{10}, -\frac{8}{5}; \frac{3}{5})$, $(-\frac{3}{10}, -\frac{4}{5}, -\frac{11}{10}, -\frac{7}{5}, -\frac{8}{5}; \frac{3}{5})$, and $(-\frac{9}{40}, -\frac{29}{40}, -\frac{43}{40}, -\frac{27}{20}, -\frac{3}{2}, -\frac{13}{8}; \frac{13}{20})$. Equations~\eqref{eq:rat_spine_steps} and~\eqref{eq:rat_overhang} of Section~\ref{subsec:RationalConstructor} and~\eqref{eq:rat_shield_prop} to~\eqref{eq:rat_top_vpos} of Appendix~\ref{subsec:ReconIndex} then express every coordinate, contact force, and contact torque as an element of $\mathbb{Q}(x)$, evaluated at $x = \xi_k(N)$, and the overhang is $x_1 + 1 = \sum_r d_r = S^*(N)$.

The verification has two parts. First, the force residual and the moment residual of every block, recomputed from its individual contact forces and torques, are the zero element of $\mathbb{Q}(x)$, and so is $\sum_i X_i$; balance therefore holds for every $x$, in particular at $\xi_k(N)$. Second, every inequality has the form $g(\xi_k(N)) > 0$ or $g(\xi_k(N)) \ge 0$ with $g \in \mathbb{Q}(x)$: positive contact forces and added loads, contact points inside the overlap of the two blocks, an overlap of positive length for every contact that carries a force, a separation of at least one unit between blocks of the same layer, the table contact, a unique principal block, and, for each spine block, a horizontal distance greater than one to every block of the layer below except the next spine block, so that the support set has exactly one block per layer. Some contact points lie at an end of their overlap by construction; the corresponding differences are identically zero. Section~\ref{subsec:ExactSoftware} describes how the inequalities are decided, and Table~\ref{tab:exact_margins} lists the smallest margins.

\begin{plaintable}[t]
\centering
\caption{Smallest margins of the exact stacks, truncated downward at the sixth decimal: positive contact force, overlap length of a contact that carries a force, positive distance from a contact point to an end of its overlap, same-layer clearance $|X_i - X_j| - 1$, and the clearance $|X_i - X_j| - 1$ between a spine block and the blocks of the layer below other than the next spine block. At $N \in \{5, 7\}$, the rows refer to the stacks of Proposition~\ref{prop:exception_sizes}; at $N = 4$, layer~$1$ contains only the base block. In every stack, the principal block extends at least $0.49$ beyond every other block.}
\label{tab:exact_margins}
\footnotesize
\setlength{\tabcolsep}{4pt}
\resizebox{\textwidth}{!}{%
\begin{tabular}{@{}rccccc@{\qquad}rccccc@{}}
\toprule
$N$ & force & overlap & end & layer & support & $N$ & force & overlap & end & layer & support \\
\midrule
4 & 0.414213 & 0.142766 & 0.142766 & 0.167893 & -- & 12 & 0.247881 & 0.123872 & 0.013269 & 0.035693 & 0.386993 \\
5 & 0.261387 & 0.427722 & 0.151207 & 0.024554 & 0.524554 & 13 & 0.262949 & 0.019592 & 0.019592 & 0.000407 & 0.383041 \\
6 & 0.277909 & 0.126320 & 0.080797 & 0.027964 & 0.636700 & 14 & 0.277120 & 0.070277 & 0.018019 & 0.002789 & 0.374504 \\
7 & 0.339468 & 0.281994 & 0.038637 & 0.023402 & 0.730045 & 15 & 0.239081 & 0.003579 & 0.001049 & 0.001420 & 0.321355 \\
8 & 0.231795 & 0.068676 & 0.068676 & 0.008094 & 0.437418 & 16 & 0.248198 & 0.005579 & 0.004670 & 0.005459 & 0.306524 \\
9 & 0.262078 & 0.175113 & 0.013821 & 0.026564 & 0.511444 & 17 & 0.256859 & 0.011130 & 0.011130 & 0.028735 & 0.388869 \\
10 & 0.289915 & 0.225107 & 0.015980 & 0.046253 & 0.574748 & 18 & 0.265113 & 0.018786 & 0.018786 & 0.049156 & 0.418767 \\
11 & 0.231777 & 0.147555 & 0.015488 & 0.015101 & 0.350383 & 19 & 0.240037 & 0.021078 & 0.006386 & 0.008964 & 0.247211 \\
\bottomrule
\end{tabular}}
\end{plaintable}

\begin{proposition}[Exact Spinal Optima]
\label{prop:exact_spinal}
For every $N \in \{4, 6, 8, \dots, 19\}$, $S(N) = S^*(N)$, and the value is attained by a single-root spinal stack of $N$ ordinary blocks that is balanced in the sense of Paterson and Zwick and feasible in our model. In particular, Conjecture~3.9 of~\cite{paterson2009overhang} holds for these $N$.
\end{proposition}

\begin{proof}
The stack constructed above is balanced, has exactly one support block in each layer, and has overhang $S^*(N)$, so $S(N) \ge S^*(N)$. Conversely, no spinal stack has overhang above $S^*(N)$~\cite[Section~3]{paterson2009overhang}, and Proposition~\ref{prop:spinal_bound}(ii) gives the same bound for spinal single-root stacks of our model.
\end{proof}

\subsection{The Exception Sizes $N = 5$ and $N = 7$}
\label{subsec:ExceptionSizes}

At $N \in \{5, 7\}$, the optimal height $k$ of Table~\ref{tab:exact_heights} satisfies $2k > N$. The following bound shows that $S^*(N)$ is then not attained and determines $S(N)$.

\begin{lemma}[Fixed-Load Bound]
\label{lem:fixed_load}
Let $X$ be a balanced spinal $N$-block stack with spine $B_1, \dots, B_k$, and let $\omega_r \ge 0$ be the total downward force exerted on $B_r$ by blocks outside the spine. Then $t_r = t_{r-1} + 1 + \omega_r$ for every $r$, and
\begin{equation}
\label{eq:fixed_load_bound}
z(X) \le \sum_{r=1}^{k} \left(1 - \frac{t_{r-1} + 1/2}{t_r}\right).
\end{equation}
At most $N - k$ of the loads $\omega_r$ are positive, and if a block other than $B_k$ rests on the table, then $t_k \le N - 1$.
\end{lemma}

\begin{proof}
Let $x_r$ be the left edge of $B_r$, let $p_r$ be the point at which $B_r$ transmits $t_r$ downward, with $p_k$ the resultant point of the table forces, let $a_r \in [x_r, x_r + 1]$ be the resultant point of the forces that make up $\omega_r$ (arbitrary if $\omega_r = 0$), and set $p_0 := x_1 + 1$. Because each $B_r$ with $r < k$ rests only on $B_{r+1}$ and $B_k$ rests only on the table, vertical balance of $B_r$ gives $t_r = t_{r-1} + 1 + \omega_r$, and moment balance gives $t_r p_r = t_{r-1} p_{r-1} + (x_r + 1/2) + \omega_r a_r$. Therefore
\[
t_r (p_{r-1} - p_r) = (1 + \omega_r) p_{r-1} - (x_r + 1/2) - \omega_r a_r \le \omega_r + \frac{1}{2},
\]
because $p_{r-1} \le x_r + 1$ (for $r \ge 2$, the contact point lies on $B_r$) and $a_r \ge x_r$. Dividing by $t_r \ge r$, summing over $r$, and using $p_k \le 0$ gives $z(X) = x_1 + 1 = p_0 \le \sum_r (\omega_r + 1/2)/t_r$, which is~\eqref{eq:fixed_load_bound}. A positive $\omega_r$ requires a block outside the spine in the layer directly above $B_r$; these layers are distinct for distinct $r$, so at most $N - k$ of the $\omega_r$ are positive. A block other than $B_k$ that rests on the table receives at least its own weight from the table, so $t_k \le N - 1$.
\end{proof}

\begin{proposition}[Spinal Optima at the Exception Sizes]
\label{prop:exception_sizes}
The maximum overhang of a $5$-block spinal stack is
\[
S(5) = \frac{12 - \sqrt{30}}{5} \approx 1.304554884990,
\]
and the maximum overhang of a $7$-block spinal stack is
\[
S(7) = \hat{S}_3(7) = 3 - \frac{1}{2\xi} - \frac{\xi + 1/2}{2\xi^2} - \frac{2\xi^2 + 1/2}{7} \approx 1.530045237061,
\]
where $\xi \approx 1.339468862388$ is the positive root of $8\xi^4 - 14\xi - 7 = 0$. Both values are attained by single-root spinal stacks of ordinary blocks that are balanced in the sense of Paterson and Zwick, and they bound every spinal stack of our model as well.
\end{proposition}

\begin{proof}
Let $X$ be a spinal stack with spine height $k$. If a block other than the base of the spine rests on the table, Lemma~\ref{lem:fixed_load} bounds $z(X)$ by the right-hand side of~\eqref{eq:fixed_load_bound} with $k \le t_k \le N - 1$; raising $t_k$ to $N - 1$ keeps $t_k - t_{k-1} \ge 1$ and increases this bound, so $z(X) \le S_k^*(N - 1) \le S^*(N - 1)$, and $S^*(4) < 1.1679$ and $S^*(6) < 1.4368$ lie below the claimed values. Let $X$ therefore be single-root, so $t_k = N$ and $z(X) \le S_k^*(N) \le \hat{S}_k(N)$. Interval enclosures give $\hat{S}_K(5) < 1.30455$ for $K \ne 3$ and $\hat{S}_K(7) < 1.53004$ for $K \notin \{3, 4\}$, both below the claimed values.

For $N = 5$ and $k = 3$, Lemma~\ref{lem:fixed_load} allows at most two positive $\omega_r$, so some $\omega_j$ vanishes. If $\omega_1 = 0$, then $t_1 = 1$, and~\eqref{eq:fixed_load_bound} becomes $12/5 - 3/(2t_2) - t_2/5$ with $2 \le t_2 \le 4$, a concave function with maximum $(12 - \sqrt{30})/5$ at $t_2 = \sqrt{15/2}$. On the faces $\omega_2 = 0$ and $\omega_3 = 0$, the bound stays below $1.30455$ (Section~\ref{subsec:ExactSoftware}). The construction~\eqref{eq:rat_n5_closed} with $t_2 = \sqrt{15/2}$ and the seed $X(C_2) = -18/25$ of Table~\ref{tab:rational_seeds}, verified as in Section~\ref{subsec:AlgebraicRealization}, attains $(12 - \sqrt{30})/5$.

For $N = 7$ and $k = 3$, the maximizer of $\hat{S}_3(7)$ has loads $(0, \xi, 2\xi^2, 7)$ with $\omega_1, \omega_2, \omega_3 > 0$, so every spinal stack of height $3$ has overhang at most $S_3^*(7) = \hat{S}_3(7)$, and the shielded construction with $(k, v_1) = (3, 1)$, the seeds $(-\frac{3}{4}, -\frac{6}{5}; \frac{1}{5})$ of Table~\ref{tab:rational_seeds}, and these loads attains it. For $k = 4$, at most three of the four $\omega_r$ are positive, and on each of the four faces $\omega_j = 0$ the bound~\eqref{eq:fixed_load_bound} stays below $1.53004 < \hat{S}_3(7)$ (Section~\ref{subsec:ExactSoftware}).
\end{proof}

Table~\ref{tab:exception_cases} collects the cases of the proof.

\begin{plaintable}[t]
\centering
\caption{Cases in the proof of Proposition~\ref{prop:exception_sizes}. Upper bounds are rounded upward; the face bounds come from the box covering of Section~\ref{subsec:ExactSoftware}.}
\label{tab:exception_cases}
\footnotesize
\resizebox{\textwidth}{!}{%
\begin{tabular}{@{}cllll@{}}
\toprule
$N$ & Spine height and case & Bound used & Upper bound & Attained by \\
\midrule
$5$ & several blocks on the table & $S^*(4)$ & $1.1679$ & -- \\
$5$ & $k \ne 3$ & $\hat{S}_k(5)$ & $1.2676$ & -- \\
$5$ & $k = 3$, face $\omega_1 = 0$ & \eqref{eq:fixed_load_bound} & $(12 - \sqrt{30})/5$ & construction~\eqref{eq:rat_n5_closed}, $t_2 = \sqrt{15/2}$ \\
$5$ & $k = 3$, faces $\omega_2 = 0$ and $\omega_3 = 0$ & \eqref{eq:fixed_load_bound} & $1.3041667$ & -- \\
\midrule
$7$ & several blocks on the table & $S^*(6)$ & $1.4368$ & -- \\
$7$ & $k \notin \{3, 4\}$ & $\hat{S}_k(7)$ & $1.4825$ & -- \\
$7$ & $k = 3$ & $S_3^*(7) = \hat{S}_3(7)$ & $\hat{S}_3(7)$ & shielded construction, $(k, v_1) = (3, 1)$ \\
$7$ & $k = 4$, faces $\omega_j = 0$, $j = 1, \dots, 4$ & \eqref{eq:fixed_load_bound} & $1.5300399$ & -- \\
\bottomrule
\end{tabular}}
\end{plaintable}

\subsection{Verification Software}
\label{subsec:ExactSoftware}

On a face $\omega_j = 0$ of height $k$ at $N \in \{5, 7\}$, the remaining added loads are nonnegative with sum $N - k$, so the loads $t_1, \dots, t_{k-1}$ lie in a polytope contained in the rational box $s \le t_s \le N - (k - s)$, $s = 1, \dots, k - 1$. We cover this box by bisection. On a box with lower limits $\underline{t}_r$ and upper limits $\overline{t}_r$ (with $t_0 = 0$ and $t_k = N$ fixed), the bound~\eqref{eq:fixed_load_bound} is at most $k - \sum_r (\underline{t}_{r-1} + 1/2)/\overline{t}_r$, because each term decreases in $t_{r-1}$ and increases in $t_r$. A box is discarded as empty when it satisfies one of the linear conditions $\overline{t}_s < \underline{t}_{s-1} + 1$ for some $s$, $\underline{t}_j > \overline{t}_{j-1} + 1$, or $\overline{t}_j < \underline{t}_{j-1} + 1$; it is closed when its bound, computed in exact rational arithmetic, falls below the threshold $1.30455$ at $N = 5$ or $1.53004$ at $N = 7$; otherwise it is bisected along its longest edge. At $N = 5$, the faces $\omega_2 = 0$ and $\omega_3 = 0$ close with $55$ and $17$ boxes and largest bounds $1.3041667$ and $1.3031250$. At $N = 7$, the faces $\omega_1 = 0, \dots, \omega_4 = 0$ of height $4$ close with $4{,}311$, $1{,}541$, $1{,}584$, and $644$ boxes and largest bounds $1.5300399$, $1.5300262$, $1.5300246$, and $1.5299587$. Every leaf of these bisection trees, empty or closed, is archived with its path from the initial box. A separate program, which shares no code with the verifiers below, replays the trees in exact rational arithmetic: it checks that the leaves of each face form a complete bisection tree of the initial box, that every empty leaf satisfies one of the linear conditions, and that every other leaf has the archived bound and that this bound lies below the threshold.

Two independent implementations perform all other checks of this appendix; they share only the input data, namely the seeds of Table~\ref{tab:rational_seeds} and, for $N \in \{17, 18, 19\}$, the seeds of Section~\ref{subsec:AlgebraicRealization}. The first uses only the Python standard library: it isolates each root $t_1 = x$ in an interval of width $3 \times 10^{-26}$ whose endpoints are checked in exact rational arithmetic, encloses every quantity in intervals with rational endpoints rounded outward to multiples of $10^{-60}$, decides every inequality from these enclosures, and covers the exception faces with its own box subdivision; an optional SymPy script checks the balance identities in $\mathbb{Q}(x)$ and the irreducibility of the polynomials. The second recomputes every construction in the field $\mathbb{Q}(x)$ with SymPy and checks the balance identities there; it decides every inequality from interval enclosures of the numerator and the denominator over the isolating interval of the root and, where these are inconclusive, by exact Sturm sign counting, encloses the values $\hat{S}_K(N)$ with the interval arithmetic of mpmath, and performs the box covering whose leaves are replayed as described above. Root isolation and sign decisions are thus made twice, by different methods, and the identities in $\mathbb{Q}(x)$ are checked by the second program and by the optional script of the first. Both programs reach the conclusions of this appendix for every $N \in \{4, \dots, 19\}$. Their outputs, including the per-block balance identities in $\mathbb{Q}(x)$, the irreducible polynomials of the roots, all margins, and the leaves of the box covering, are archived in the companion repository together with a script that reruns all three programs and compares the results with the archive.

\end{appendices}

\clearpage
\bibliographystyle{sn-mathphys-num}
\bibliography{references}

\end{document}